\documentclass[12pt]{iopart}
\usepackage{lipsum}
\usepackage{amsfonts}
\usepackage{graphicx}
\usepackage{amsmath}
\usepackage{amsthm}
\usepackage{multirow}
\usepackage{algorithmic}
\usepackage{caption}
\usepackage{afterpage}
\usepackage{subcaption}
\usepackage{lineno}
\usepackage{hyperref}
\usepackage{epstopdf}
\usepackage{mathtools}
\usepackage{algorithm}
\ifpdf
  \DeclareGraphicsExtensions{.eps,.pdf,.png,.jpg}
\else
  \DeclareGraphicsExtensions{.eps}
\fi

\DeclareMathOperator*{\argmin}{arg\,min}
\DeclareMathOperator*{\argmax}{arg\,max}
\DeclareMathOperator{\sign}{sgn}
\newcommand{\WDNote}[1]           
{\textcolor{red}{#1}\marginpar{\textcolor{red}{WD $\longleftarrow$}}}
\usepackage{rotating}
\newcommand{\rottext}[1]{\rotatebox{90}{\hbox to 27mm{\hss #1\hss}}}
\newcommand{\rottextt}[1]{\rotatebox{90}{\hbox to 4mm{\hss #1\hss}}}
\newcommand{\rottextttt}[1]{\rotatebox{90}{\hbox to 40mm{\hss #1\hss}}}
\usepackage{float}
\newtheorem{thm}{Theorem}
\newtheorem{lemma}[thm]{Lemma}
\newtheorem{assumption}[thm]{Assumption}
\newtheorem{theorem}[thm]{Theorem}
\newtheorem{proposition}[thm]{Proposition}
\newtheorem{definition}[thm]{Definition}
\newtheorem{remark}{Remark}
\counterwithin{thm}{section}
\usepackage{booktabs}
\usepackage{multirow}
\usepackage{siunitx}
\usepackage{makecell}
\usepackage{afterpage}
\usepackage{color,soul}
\usepackage{enumerate}
\usepackage{enumitem}
\usepackage{comment}
\allowdisplaybreaks
\usepackage{placeins}

\renewcommand{\appendix}{\par
  \setcounter{section}{0}
  \setcounter{subsection}{0}
  \gdef\thesection{\Alph{section}}
}

\makeatletter
\def\maketag@@@#1{\hbox{\m@th\normalfont\normalsize#1}}
\makeatother

\begin{document}

\title[Stochastic ADMM for Image Ptychography]{A Stochastic ADMM Algorithm for Large-Scale Ptychography with Weighted
Difference of Anisotropic and Isotropic Total Variation}

\author{Kevin Bui$^1$ and Zichao (Wendy) Di$^2$}

\address{$^1$ Department of Mathematics, University of California at Irvine, Irvine, CA 92697 USA\\
$^2$ Mathematics and Computer Science Division, Argonne National Laboratory, Lemont, IL 60439 USA}
\ead{kevinb3@uci.edu}
\vspace{10pt}

\begin{abstract}
Ptychography, a prevalent imaging technique in fields such as biology and optics, poses substantial challenges in its reconstruction process, characterized by nonconvexity and large-scale requirements. This paper presents a novel approach by introducing a class of variational models that incorporate the weighted difference of anisotropic--isotropic total variation. This formulation enables the handling of measurements corrupted by Gaussian or Poisson noise, effectively addressing the nonconvex challenge. To tackle the large-scale nature of the problem, we propose an efficient stochastic alternating direction method of multipliers, which guarantees convergence under mild conditions. Numerical experiments validate the superiority of our approach by demonstrating its capability to successfully reconstruct complex-valued images, especially in recovering the phase components even in the presence of highly corrupted measurements.
\end{abstract}

%
\vspace{2pc}
\submitto{\IP}
%
%
%

\section{Introduction}
\begin{figure}
    \centering
    \includegraphics[scale=0.30]{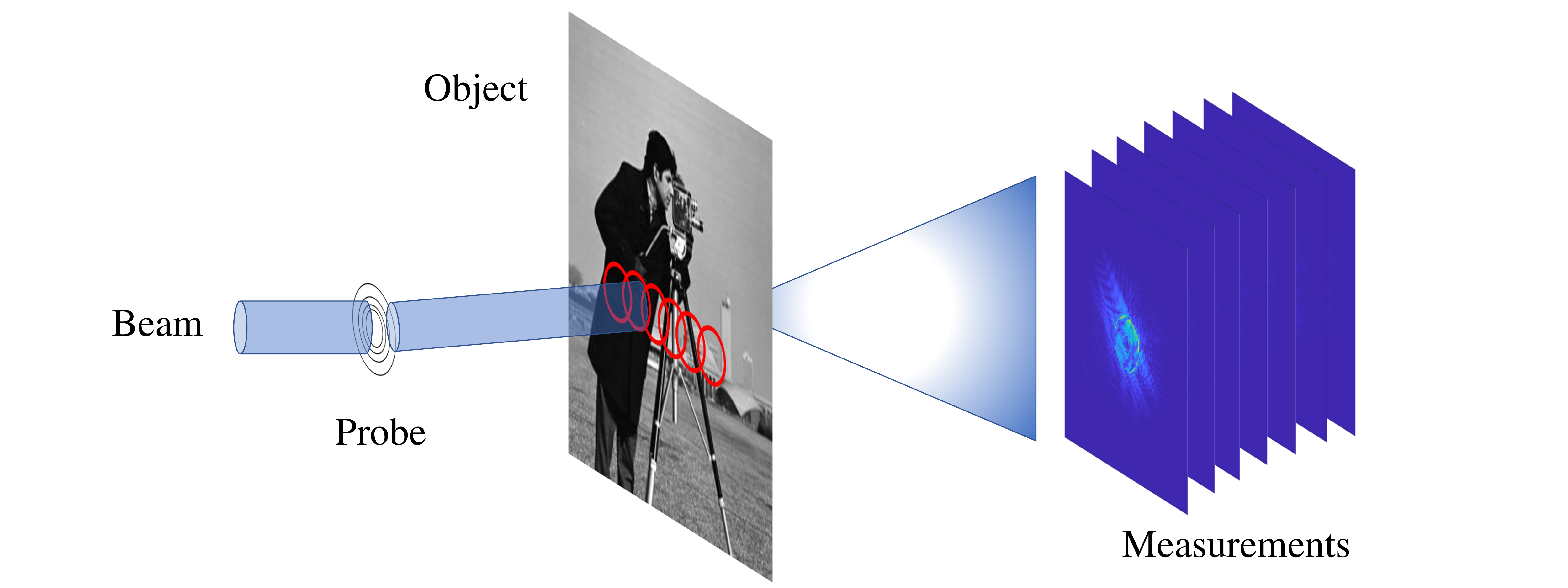}
    \caption{Schematic of a ptychography experiment.}
    \label{fig:ptycho_experiment}
\end{figure}

Ptychography is a popular imaging technique that combines both coherent diffractive imaging and scanning transmission microscopy. It is used in various industrial and scientific applications, including biology \cite{marrison2013ptychography, suzuki2016dark, zhou2020low}, crystallography \cite{de2016ptychographic}, and optics \cite{shechtman2015phase,walther1963question}. To perform a ptychographic experiment (see Figure \ref{fig:ptycho_experiment}), a coherent beam is scanned across the object of interest, where each scan may have overlapping positions with another. The scanning procedure provides a set of phaseless measurements that can be used to reconstruct an image of the object of interest.

Algorithms have been developed to solve the nonblind ptychography problem, where the probe is known, or, more generally, the phase retrieval problem. One of the most popular methods is the ptychographical iterative engine (PIE) \cite{rodenburg2004phase}, which applies gradient descent to the object using each measurement sequentially. Convergence analysis of PIE is provided in \cite{melnyk2022stochastic}.  Other gradient-based methods for phase retrieval include Wirtinger flow \cite{candes2015phase} and its variants \cite{chen2017solving, yuan2019phase, zhang2017nonconvex}, which use adaptive step sizes and careful initialization based on spectral method. PIE belongs to the class of projection methods for phase retrieval, where the object update is projected onto a nonconvex modulus constraint set. Other projection-based algorithms are hybrid projection-reflection \cite{bauschke2003hybrid}, Douglas-Rachford splitting \cite{thibault2009probe}, and the relaxed averaged alternating reflections \cite{luke2004relaxed}. For Douglas-Rachford splitting, fixed-point analysis and local, linear convergence were studied extensively \cite{chen2018coded, chen2018fourier, fannjiang2020fixed}. The phase retrieval problem can be formulated as variational models \cite{chang2018variational, wen2012alternating} that are solved by alternating direction method of multipliers (ADMM) \cite{boyd2011distributed} or proximal alternating linearized minimization \cite{bolte2014proximal}.  Incorporating total variation regularization \cite{chang2018total, chang2016phase, chang2018variational} improves the performance and robustness of these models, especially when the magnitude measurements are corrupted by noise.

Some algorithms for the nonblind ptychography problems were extended to jointly solve for the probe and object in the more challenging, blind ptychography problem, where the probe is unknown. Later refinements to PIE led to ePIE \cite{maiden2009improved} and rPIE \cite{maiden2017further}, which apply gradient descent updates with adaptive, stable step sizes to both the object and the probe. A convergence analysis of ePIE is provided in \cite{melnyk2023convergence}.  A Wirtinger flow algorithm was extended to blind polychromatic ptychography via alternating minimization and was proven to sublinearly converge to a stationary point \cite{filbir2023image}.  A semi-implicit relaxed version of DR \cite{pham2019semi} was proposed as a robust algorithm to solve for both probe and object, especially with sparse data. Based on Kurdyka-Łojasiewicz assumptions \cite{bolte2014proximal}, globally convergent algorithms, such as ADMM \cite{chang2019blind} and proximal alternating linearized minimization \cite{hesse2015proximal}, were proposed for blind ptychography. 

Large-scale ptychography, characterized by a substantial number of scans, presents challenges in terms of memory usage and computational cost for existing algorithms. Parallel algorithms have been developed to address these challenges. For example, an asynchronous, parallel version of ePIE was implemented on GPUs, where sub-images from partitioned measurement sets are fused together \cite{nashed2014parallel}. Additionally, a parallel version of relaxed averaged alternating reflections was adapted for GPU implementation \cite{enfedaque2019high}. However, these parallel algorithms often rely on GPUs, which may not be universally available. Alternatively, some algorithms designed for large-scale ptychography without high-performance computing resources have been proposed. For instance, a multigrid optimization framework accelerates gradient-based methods for phase retrieval \cite{fung2020multigrid}, while an overlapping domain decomposition method combined with ADMM enables a highly parallel algorithm with good load balance \cite{chang2021overlapping}. A batch-based ADMM algorithm {\cite{yang2022batch}} was proposed to iteratively process a subset of measurement data instead of all for image reconstruction. Overall, these algorithms require careful tailoring specific to ptychography to demonstrate their benefits. Notably, to the best of our knowledge, a batch-based, stochastic ADMM algorithm for large-scale ptychography has yet to be developed. By incorporating stochastic gradient descent (SGD) \cite{robbins1951stochastic}, stochastic ADMM may be able to find better solutions \cite{kleinberg2018alternative} than the deterministic ADMM algorithms \cite{chang2019blind,chang2018total, chang2016phase, yang2022batch} for the ptychography problem that is formulated as a nonconvex optimization problem. Moreover, being batch based, it would be computationally practical for practitioners without access to multiple cores for parallel computing.

To improve the image reconstruction quality in phase retrieval, total variation (TV) \cite{rudin1992nonlinear} has been incorporated  for the cases when the measurements are corrupted with Gaussian noise \cite{chang2016phase} or with Poisson noise \cite{chang2018total}. 
Both cases consider the isotropic TV approximation
\begin{align}
    \|\nabla z\|_{2,1} = \sum_{i=1}^{n^2} \sqrt{|(\nabla_x z)_i|^2 + |(\nabla_y z)_i|^2},
\end{align}
where $z \in \mathbb{C}^{n^2}$ is an $n \times n$ image, $\nabla_x$ and $\nabla_y$ are the horizontal and vertical difference operators, respectively, and $(\nabla_x z)_i$ and $(\nabla_y z)_i$ are the $i$th entries in lexicographic ordering of $\nabla_x z$ and $\nabla_y z$, respectively. 
However, it has been known that isotropic TV tends to blur oblique edges. An alternative approximation that preserves sharper edges is the anisotropic TV \cite{esedoglu2004decomposition}:
\begin{align}
    \|\nabla z\|_1 = \sum_{i=1}^{n^2} \left(|(\nabla_x z)_i| + |(\nabla_y z)_i|\right).
\end{align}
Overall, TV is meant to approximate the $\ell_0$ ``norm" of the image gradient, i.e., $\|\nabla z\|_0$, because TV is based on the $\ell_1$ norm, a convex relaxation of $\ell_0$. A nonconvex alternative to $\ell_1$ is $\ell_1 - \alpha \ell_2, 0 < \alpha \leq 1$, which performs well in recovering sparse solutions in various compressed sensing problems \cite{lou2018fast,lou2015computing, lou2016point,  yin2015minimization}. The superior performance of $\ell_1 - \alpha \ell_2$ in sparse recovery has motivated the development of the weighted difference of anisotropic and isotropic total variation (AITV) \cite{lou2015weighted}, which applies $\ell_1 -\alpha \ell_2$ on each gradient vector of an image. Mathematically, AITV is formulated by
\begin{align}
    \|\nabla z \|_1 - \alpha \|\nabla z\|_{2,1} = \sum_{i=1}^{n^2} \left[|(\nabla_x z)_i| + |(\nabla_y z)_i| - \alpha  \sqrt{|(\nabla_x z)_i|^2 + |(\nabla_y z)_i|^2} \right].
\end{align}
AITV has demonstrated better performance than TV in image denoising, image deconvolution, image segmentation, and MRI reconstruction \cite{bui2021weighted, lou2015weighted, park2016weighted}, especially in preserving sharper edges.

In this work, we focus on addressing the large-scale ptychography problem with measurements corrupted by Gaussian or Poisson noise. To enhance image reconstruction quality, we incorporate AITV regularization within a general variational framework. The problem is then solved using an ADMM algorithm, where some subproblems are approximated by SGD \cite{robbins1951stochastic}. However, direct application of SGD to the ptychography subproblem is not straightforward. Therefore, we demonstrate the appropriate adaptation of SGD to develop our specialized stochastic ADMM algorithm. Furthermore, we leverage the inherent structure embedded in the experiment such as per-pixel illumination strength and automate the selection of optimal step sizes within the minibatch regime. 
Unlike traditional methods using all measurements per iteration, our stochastic ADMM algorithm iteratively processes a random batch of measurements, resulting in accurate and efficient ptychographic reconstruction.

The paper is organized as follows. In Section \ref{sec:models}, we describe the AITV-regularized variational models to solve the image ptychography problem. Within this section, we design the stochastic ADMM algorithms to solve these models. Convergence analysis of the stochastic ADMM algorithm follows in Section \ref{sec:convergence}. In Section \ref{sec:experiment}, we illustrate the performance of our proposed stochastic ADMM algorithms and compare them with other competing algorithms. Lastly, in Section \ref{sec:conclude}, we conclude the paper with summary and future works. 

\section{Mathematical Model}\label{sec:models}
We first describe  basic notations used throughout the paper. Let $z=[z_i]_{i=1}^{n^2} \in \mathbb{C}^{n^2}$ be an $n^2 \times 1$ vector that represents an $n \times n$ image in lexicographical order, where the rows of an image matrix are transposed and stacked consecutively into one vector. The vector $\mathbf{1}$ is a vector whose entries are all ones. Superscripts $^{\top}$ and $^*$ denote the real and conjugate transpose of a matrix, respectively. The sign of a complex value $z' \in \mathbb{C}$ is given by
\begin{align*}
    \sign(z') = \begin{cases}
    \displaystyle \frac{z'}{|z'|} ,&\text{ if } z' \neq 0,\\
    c \in \{c' \in \mathbb{C}: |c'| \leq 1\}, &\text{ if } z' = 0.
    \end{cases}
\end{align*}
The sign of a vector $z \in \mathbb{C}^{n^2}$ is denoted by $\sign(z)$ and is defined elementwise by $
    \sign(z)_i = \sign(z_i), \; i=1, \ldots, n^2.$
The vector $\mathbf{e}_i$ has 1 as the $i$th component while all other entries are zeros. For $p = (p_x, p_y) \in \mathbb{C}^{n^2} \times \mathbb{C}^{n^2}$, its $i$th entry is $p_i = \begin{bmatrix}(p_x)_i \\ (p_y)_i\end{bmatrix} \in \mathbb{C}^2$. We define the following norms on $\mathbb{C}^{n^2} \times \mathbb{C}^{n^2}$:
\begin{align*}
    \|p\|_1 &= \sum_{i=1}^{n^2} |(p_x)_i| + |(p_y)_i|, \quad
    \|p\|_2 = \sqrt{ \sum_{i=1}^{n^2} |(p_x)_i|^2+ |(p_y)_i|^2 }, \\
    \|p\|_{2,1} &=  \sum_{i=1}^{n^2} \sqrt{|(p_x)_i|^2+|(p_y)_i|^2}
= \sum_{i=1}^{n^2} \|p_i\|_2.
\end{align*}
The discrete gradient operator $\nabla: \mathbb{C}^{n^2} \rightarrow \mathbb{C}^{n^2} \times \mathbb{C}^{n^2}$ when specifically applied to the image $z$ is given by $\nabla z = (\nabla_x z, \nabla_y z)$, where $\nabla_x$ and $\nabla_y$ are the forward horizontal and vertical difference operators. 
More specifically, the $i$th entry of $\nabla z$ is defined by
\begin{align}
    (\nabla z)_i = \begin{bmatrix}
    (\nabla_x z)_i\\
    (\nabla_y z)_i
    \end{bmatrix},
\end{align}
where 
\begin{align*}
    (\nabla_x z)_i = \begin{cases}
    z_i - z_{i-1}, & \text{ if } i \text{ mod }n \neq 1,\\
    z_i - z_{i+n-1}, & \text{ if } i \text{ mod } n = 1
    \end{cases}
\end{align*}
and
\begin{align*}
    (\nabla_y z)_i = \begin{cases}
    z_i - z_{i-n}, & \text{ if } i > n,\\
    z_{i} - z_{i+(n-1)n}, & \text{ if } 1 \leq i \leq n.
    \end{cases}
\end{align*}
We also define the proximal operator of a function $f: \mathbb{C}^{n^2} \rightarrow \mathbb{R} \cup \{+\infty\}$ as
\begin{align*}
    \text{prox}_{f(\cdot)}(\tilde{z}) = \argmin_z f(z) + \frac{1}{2}\|z-\tilde{z}\|_2^2, \, \forall \tilde{z} \in \mathbb{C}^{n^2}.
\end{align*}
Lastly, considering complex functions, we compute their gradient according to Wirtinger calculus~\cite{kreutz2009complex}.

We now describe the 2D ptychography in the discrete setting based on the notations in \cite{chang2019blind}. Let $z \in \mathbb{C}^{n^2}$ be the object of interest with $n \times n$ pixels and $\omega \in \mathbb{C}^{m^2}$ be the localized 2D probe with $m \times m$ pixels, where $m < n$. Both the object $z$ and the probe $\omega$ are expressed as vectors in lexicographical order. We denote the set of masks by $\{S_j\}_{j=1}^N$ for $N$ total scanning positions, where each $S_j \in \mathbb{R}^{m^2 \times n^2}$ is a binary matrix that corresponds to the $j$th scanning window of the probe $\omega$ on the object $z$. As a result, $S_jz$ represents a different $m \times m$ subwindow of the object $z$. The set of phaseless measurements $\{d_j\}_{j=1}^N$ is obtained by
    $d_j =  | \mathcal{F}(\omega \circ S_j z) |^2$,
where $\mathcal{F} \in \mathbb{C}^{m^2 \times m^2}$ is the normalized 2D discrete Fourier operator, $\circ$ is the elementwise multiplication, and $|\cdot |$ is the elementwise absolute value of a vector. We note that the Fourier operator $\mathcal{F}$ is unitary, i.e., $\mathcal{F}^* = \mathcal{F}^{-1}$ and $\|\mathcal{F}x\|_2 = \|x\|_2$ for any $x \in \mathbb{C}^{m^2}$.  For each $j=1,\ldots, N$, the measurement $d_j$ and the mask matrix $S_j$ form together the $j$th scan information $(d_j, S_j)$. Throughout the paper, we assume that for each $i=1, \ldots, n^2$, there exists $j \in \{1, \ldots, N\}$ such that $\|S_{j}\mathbf{e}_i\|_1=1$. This assumption ensures that each pixel of an image $z \in \mathbb{C}^{n^2}$ is scanned at least once.

The  blind ptychographic phase retrieval problem is expressed as follows:
\begin{align} \label{eq:blind_pr}
    \text{To find } \omega \in \mathbb{C}^{m^2} \text{ and } z \in \mathbb{C}^{n^2} \text{ such that } |\mathcal{F}(\omega \circ S_j z) |^2 = d_j, \; j=1, \ldots, N.
\end{align} When the probe $\omega$ is known, \eqref{eq:blind_pr} reduces to the non-blind case where we only find $z \in \mathbb{C}^{n^2}$. Suppose that the measurements $\{d_j\}_{j=1}^N$ are corrupted by independent and identically distributed (iid) noise. The blind variational model~\cite{chang2019blind} is formulated by
\begin{align} \label{eq:blind_model}
    \min_{\omega \in \mathbb{C}^{m^2}, z \in \mathbb{C}^{n^2}} \sum_{j=1}^N \mathcal{B}(|\mathcal{F}(\omega \circ S_j z)|^2, d_j),
\end{align} 
where
\begin{align}\label{eq:phase_fidelity}
    \mathcal{B}(g,f) = \displaystyle \begin{cases}
   \displaystyle \frac{1}{2} \|\sqrt{g} - \sqrt{f}\|_2^2, &\text{ amplitude Gaussian metric (AGM) \cite{wen2012alternating}},\\ \\
    \displaystyle \frac{1}{2} \langle g-f \circ \log(g), \mathbf{1} \rangle, & \text{ intensity Poisson metric (IPM) \cite{chang2018total}}.
    \end{cases}
\end{align}
Note that $\sqrt{\cdot}$ is elementwise square root. The AGM metric is suitable when the magnitude measurements $\{\sqrt{d_j}\}_{j=1}^N$ are corrupted by Gaussian noise, while the IPM metric is appropriate when the measurements $\{d_j\}_{j=1}^N$ are corrupted by Poisson noise. Nevertheless, AGM can also be used for the Poisson noise case because it is a variance-stablizing transform {\cite{marchesini2016sharp}} and an approximation for IPM {\cite{thibault2012maximum}}.

\subsection{AITV model}\label{subsec:nonblind_model}
To improve image recovery in blind ptychography, we propose a class of AITV-regularized variants of \eqref{eq:blind_model}:
\begin{align}\label{eq:AITV_blind_prob}
     \min_{\omega \in \mathbb{C}^{m^2}, z \in \mathbb{C}^{n^2}} \sum_{j=1}^N \mathcal{B}(|\mathcal{F}(\omega \circ S_j z)|^2, d_j) + \lambda \left( \|\nabla z\|_1 - \alpha \|\nabla z\|_{2,1} \right ), \; \lambda > 0,\; \alpha \in [0,1].
\end{align}
The non-smoothness of the objective function \eqref{eq:AITV_blind_prob} with respect to $\omega$ and $z$ motivates the development of an ADMM algorithm for solving it. Specifically, we introduce auxiliary variables $u = (u_1, \ldots, u_N) \in \mathbb{C}^{m^2 \times N}$ to represent $(\mathcal{F}(\omega \circ S_1 z), \ldots, \mathcal{F}(\omega \circ S_N z))$, and $v = (v_x, v_y) \in \mathbb{C}^{n^2} \times \mathbb{C}^{n^2}$ to represent $\nabla z = (\nabla_x z, \nabla_y z)$. These auxiliary variables play a crucial role in handling the non-smoothness of the objective function and enable the effective utilization of the ADMM framework for optimization. As a result, we obtain an  equivalent constrained optimization problem
\begin{gather} \label{eq:blind_constrained}
    \begin{aligned}
    \min_{u, \omega, v, z} \quad  \sum_{j=1}^N \mathcal{B}(|u_j|^2, d_j) + \lambda \left( \|v\|_1 - \alpha \|v\|_{2,1} \right ) \; \textrm{s.t. }  &u_j = \mathcal{F}(\omega \circ S_j z), \;j = 1, \ldots, N,\\ &v = \nabla z.
    \end{aligned}
\end{gather}
The augmented Lagrangian of \eqref{eq:blind_constrained} is
\begin{gather}\label{eq:blind_Lagrangian}
\begin{aligned}
   \mathcal{L}(u, \omega,  v, z, \Lambda, y)= &
    \sum_{j=1}^N \left[ \mathcal{B}(|u_j|^2, d_j)+ \mathbb{R}\left(\langle \Lambda_j, u_j- \mathcal{F}(\omega \circ S_j z) \rangle \right) + \frac{\beta_1}{2}  \|u_j - \mathcal{F}(\omega \circ S_j z) \|_2^2 \right]\\ &  + \lambda \left( \|v\|_1 - \alpha \|v\|_{2,1} \right )+ \mathbb{R} \left(\langle y, v - \nabla z \rangle \right) + \frac{\beta_2}{2} \|v - \nabla z \|_2^2,
\end{aligned}
\end{gather}
where $\mathbb{R}(\cdot)$ denotes the real component of a complex number; $\langle \cdot, \cdot \rangle$ denotes the complex inner product between two vectors; $\Lambda = (\Lambda_1, \ldots, \Lambda_N) \in \mathbb{C}^{m^2 \times N}$ and $y = (y_x, y_y) \in \mathbb{C}^{n^2} \times \mathbb{C}^{n^2}$ are Lagrange multipliers; and $\beta_1, \beta_2 >0$ are penalty parameters. The ADMM algorithm iterates as follows:
\begin{subequations}\label{eq:admm1}
\begin{align}
u^{t+1} &\in \argmin_{u} \mathcal{L}(u, \omega^t, v^{t}, z^{t}, \Lambda^t, y^t),\label{eq:u_sub}\\
\omega^{t+1}  &\in \argmin_{\omega} \mathcal{L}(u^{t+1}, \omega, v^{t}, z^{t}, \Lambda^t, y^t),\label{eq:omega_subprob}\\
v^{t+1} &\in \argmin_{v} \mathcal{L}(u^{t+1}, \omega^{t+1}, v, z^{t}, \Lambda^t, y^t), \label{eq:v_sub}\\
z^{t+1} &\in \argmin_{z} \mathcal{L}(u^{t+1}, \omega^{t+1}, v^{t+1}, z ,\Lambda^t, y^t), \label{eq:z_sub}\\
\Lambda_j^{t+1} &= \Lambda_j^t + \beta_1\left(u_j^{t+1} - \mathcal{F}(\omega^{t+1} \circ S_j z^{t+1}) \right), \quad j=1, \ldots, N, \label{eq:lambda_update}\\
y^{t+1} &= y^t + \beta_2\left(v^{t+1} - \nabla z^{t+1} \right). \label{eq:y_update}
\end{align}
\end{subequations}
We explain how to solve each subproblem and adapt it in the stochastic, batch-based setting, where at each iteration we only have a randomly sampled mini-batch $n^t  \subset \{1, \ldots, N\}$ of scan information available such that $|n^t| = b$. 
As a result, \eqref{eq:lambda_update} reduces to
\begin{align}
\Lambda_j^{t+1} = \begin{cases} \Lambda_j^t + \beta_1\left(u_j^{t+1} - \mathcal{F}(\omega^{t+1} \circ S_j z^{t+1}) \right), &\text{ if } j \in n^t,\\
 \Lambda_j^t,  &\text{ if } j \not \in n^t.
\end{cases} \label{eq:lambda_update2} \tag{\ref{eq:lambda_update}$'$}
\end{align}
\subsubsection{$u$-subproblem}Solving \eqref{eq:u_sub} simplifies to solving $u_j$ independently, so we only need to update $u_j$'s whose corresponding scan information is available. Let $D_{\omega} = (\omega \mathbf{1}^{\top}) \circ I_{m^2 \times m^2}$ be a diagonal matrix whose diagonal is $\omega \in \mathbb{C}^{m^2}$. Then for each $j \in n^t$,  we have 
\begin{gather}
\begin{aligned}\label{eq:u_min}
    &u_j^{t+1} \in \argmin_{u_j}\mathcal{B}(|u_j|^2, d_j) +  \mathbb{R}\left(\langle \Lambda_j^t, u_j- \mathcal{F}(P_j^t z^t) \rangle \right) + \frac{\beta_1}{2}  \|u_j - \mathcal{F}(P_j^t
    z^t) \|_2^2\\
    &=\argmin_{u_j} \frac{1}{\beta_1} \mathcal{B}(|u_j|^2, d_j) + \frac{1}{2}\left \|u_j - \mathcal{F}(P_j^t z^t) + \frac{1}{\beta_1} \Lambda_j^t \right \|_2^2\\
    &=\text{prox}_{\frac{1}{\beta_1}\mathcal{B}(|\cdot|^2, d_j)}\left( \mathcal{F}(P_j^t z^t) - \frac{1}{\beta_1} \Lambda_j^t  \right),
    \end{aligned}
\end{gather}
where $P_j^t \coloneqq D_{\omega^t}S_j$ since $\omega^t \circ S_j z^t = D_{\omega^t}S_j z^t$.
Because the proximal operator for each fidelity term in \eqref{eq:phase_fidelity} has a closed-form solution provided in \cite{chang2018total, chang2016phase}, the overall update is
{\scriptsize
\begin{gather}
\begin{aligned}\label{eq:u_closed_form}
     u_j^{t+1} = 
    \begin{cases}\begin{cases}
     \frac{\sqrt{d_j}+\beta_1\left| \mathcal{F}(P_j^t z^t) - \frac{1}{\beta_1} \Lambda_j^t \right|}{1+ \beta_1} \circ \sign\left( \mathcal{F}(P_j^t z^t) - \frac{1}{\beta_1} \Lambda_j^t \right), &\text{AGM,} \\
    \frac{\beta_1 |\mathcal{F}(P_j^t z^t)-\frac{1}{\beta_1} \Lambda_j^t| + \sqrt{\left(\beta_1|\mathcal{F}(P_j^t z^t) - \frac{1}{\beta_1} \Lambda_j^t|\right)^2 + 4(1+\beta_1)d_j}}{2(1+\beta_1)} \circ \sign\left(\mathcal{F}(P_j^t z^t) -\frac{1}{\beta_1}\Lambda_j^t \right), & \text{IPM},
    \end{cases} & \text{ if } j \in n^t,\\
    u_j^t, & \text{ if } j \not \in n^t.
    \end{cases}
\end{aligned}
\end{gather}}%
\subsubsection{$\omega$-subproblem}\label{sec:omega_subprob}
The $\omega$-subproblem \eqref{eq:omega_subprob} can be rewritten as
\begin{gather}\label{eq:full_omega_subprob}
\begin{aligned}
    \omega^{t+1} &\in \argmin_{\omega} \sum_{j=1}^N \left[  \frac{\beta_1}{2} \left\| \mathcal{F}^{-1} \left( u_{j}^{t+1}  + \frac{\Lambda_j^t}{\beta_1} \right) - \omega \circ S_j z^t \right\|_2^2 \right].
\end{aligned}
\end{gather}
Instead of using all $N$ scan information, we develop an alternative update scheme that uses only $b$ of them. Instead of solving \eqref{eq:omega_subprob} exactly, we linearize it as done in \cite{liu2019linearized, ouyang2015accelerated} to obtain
\begin{align} \label{eq:linearize_omega_subprob}
\omega^{t+1} \in \argmin_{\omega} \mathbb{R}(\langle \nabla_{\omega}\mathcal{L}(u^{t+1}, \omega^t, z^t, \Lambda^t, y^t), \omega - \omega^t \rangle) + \frac{1}{2 \delta_{\omega}^t} \|\omega - \omega^t\|_2^2
\end{align}
for some constant $\delta_{\omega}^t > 0$ at iteration $t$. Let $\mathcal{G}^t_j(\omega) = \frac{\beta_1}{2} \left\| \mathcal{F}^{-1} \left( u_{j}^{t+1}  + \frac{\Lambda_j^t}{\beta_1} \right) - \omega \circ S_j z^t \right\|_2^2$ so that we have
\begin{align*}
\nabla \mathcal{G}_j^t(\omega) = -\beta_1 (S_j z^t)^* \circ \left[ \mathcal{F}^{-1}\left( u_{j}^{t+1}  + \frac{\Lambda_j^t}{\beta_1} \right) - \omega \circ S_j z^t \right].
\end{align*} 
Then \eqref{eq:linearize_omega_subprob} is equivalent to performing gradient descent with step size $\delta_{\omega}^t$:
\begin{gather}
\begin{aligned}\label{eq:omega_gd}
\omega^{t+1} &= \omega^t - \delta_{\omega}^t \nabla_{\omega} \mathcal{L}(u^{t+1}, \omega^{t}, v^t, z^t, \Lambda^t, y^t)\\
&= \omega^t - \delta_{\omega}^t \sum_{j=1}^N \nabla \mathcal{G}_j^t(\omega^t) \\
& =  \omega^t - N\delta_{\omega}^t\left(\frac{1}{N} \sum_{j=1}^N \nabla \mathcal{G}_j^t(\omega^t) \right).
\end{aligned}
\end{gather}
Because we only have $b$ scan information available, we replace the update term in \eqref{eq:omega_gd} with its stochastic estimate $\tilde{\nabla}_{\omega} \mathcal{L}$ as:
\begin{align}\label{eq:sgd_omega}
    \omega^{t+1} = \omega^t - \delta_{\omega}^t \tilde{\nabla}_{\omega} \mathcal{L}(u^{t+1}, \omega^{t}, v^t, z^t, \Lambda^t, y^t). 
\end{align}
For simplicity, we choose the SGD estimator \cite{bottou2018optimization, robbins1951stochastic} as $\tilde{\nabla}_{\omega} \mathcal{L}$: 
\begin{gather}
\begin{aligned}\label{eq:omega_sgd_estimator}
\tilde{\nabla}_{\omega}^{SGD}  \mathcal{L}(u^{t+1}, \omega^{t}, v^t, z^t, \Lambda^t, y^t) &=  \frac{1}{b} \sum_{j \in n^t} \nabla \mathcal{G}_j^t(\omega^t)\\ &=-\frac{\beta_1}{b}  \sum_{j \in n^t} (S_j z^t)^* \circ \left[ \mathcal{F}^{-1}\left( u_{j}^{t+1}  + \frac{\Lambda_j^t}{\beta_1} \right) - \omega^t \circ S_j z^t \right].
\end{aligned}
\end{gather}
Similar to the PIE family algorithm \cite{ maiden2017further, maiden2009improved, rodenburg2004phase}, we further extend \eqref{eq:omega_sgd_estimator} by incorporating spatially varying illumination strength modeled as
\begin{equation}
\label{eq:omega_PIE}
 \Phi_j^t =  \displaystyle \frac{\mathbf{1}}{(1-\gamma_{\omega}) |S_j z^t|^2 + \gamma_{\omega} \|S_j z^t\|_{\infty}^2\mathbf{1} },
\end{equation}
where $\gamma_{\omega} \in [0,1]$ and division is elementwise. Incorporating $\Phi_j^t$ into \eqref{eq:omega_sgd_estimator}, we have another class of stochastic estimators
\begin{gather}\label{eq:omega_stochastic2}
\begin{aligned}
    \tilde{\nabla}_{\omega}^{PIE}  \mathcal{L}(u^{t+1}, \omega^{t}, v^t, z^t, \Lambda^t, y^t) =- \frac{\beta_1}{b} \sum_{j \in n^t}  \Phi_j^t \circ (S_j z^t)^* \circ \left[ \mathcal{F}^{-1}\left( u_{j}^{t+1}  + \frac{\Lambda_j^t}{\beta_1} \right) - \omega^t \circ S_j z^t \right]. 
\end{aligned}
\end{gather}

\subsubsection{$v$-subproblem}
Expanding \eqref{eq:v_sub} gives
\begin{gather} \label{eq:v_sub2}
    \begin{aligned}
    v^{t+1} & \in  \argmin_v \frac{\lambda}{\beta_2}\left(\|v\|_1 - \alpha \|v\|_{2,1} \right) + \frac{1}{2} \left\|v - \nabla z^t + \frac{y^t}{\beta_2}\right\|_2^2\\
    &= \argmin_v \sum_{i=1}^{n^2} \frac{\lambda}{\beta_2} \left( \|v_i\|_1 - \alpha \|v_i\|_2  \right)+ \frac{1}{2} \left\|v_i - (\nabla z^t)_i + \frac{(y^t)_i}{\beta_2}  \right\|_2^2,
    \end{aligned}
\end{gather}
which means that the solution $v^{t+1}$ can be solved elementwise. As a result, the subproblem simplifies to
\begin{align}\label{eq:v_sub3}
    (v^{t+1})_i = \text{prox}_{\frac{\lambda}{\beta_2} ( \|\cdot \|_1 - \alpha \|\cdot\|_2)} \left((\nabla z^t)_i - \frac{(y^t)_i}{\beta_2} \right).
\end{align}
 A closed-form solution for the proximal operator of $\ell_1 - \alpha \ell_2$ is provided in \cite{lou2018fast}  but only for real-valued vectors. We generalize it to the complex case in Lemma \ref{lemma:prox_l1l2}, whose proof is delayed to Appendix \ref{sec:prox}.
\begin{lemma} \label{lemma:prox_l1l2}
Given $x' \in \mathbb{C}^n$, $\lambda >0$, and $\alpha \geq 0$, we have the following cases:
\begin{enumerate}
    \item When $\|x'\|_{\infty} > \lambda$, we have
    \begin{align*}
        x^* = (\|\xi\|_2 + \alpha \lambda) \frac{\xi}{\|\xi\|_2}, \text{ where } \xi = \sign(x') \circ \max(|x'|-\lambda,0).
    \end{align*}
    
    \item When $(1-\alpha) \lambda < \|x'\|_{\infty} \leq \lambda$, we have $x^*$ as a 1-sparse vector such that one chooses an index $i \in \argmax_j (|(x')_j|)$ and have
    \begin{align*}
        (x^*)_j = \begin{cases}
        (|(x')_j| + (\alpha-1)\lambda) \sign((x')_j), &\text{ if } j = i, \\
        0, &\text{ if } j \neq i.
        \end{cases}
    \end{align*}
\item When $\|x'\|_{\infty} \leq (1- \alpha)\lambda$, we have $x^* = 0$. 
\end{enumerate}
Then $x^*$ is an optimal solution to
\begin{align}\label{eq:prox_l1_l2}
    \mathrm{prox}_{\lambda(\|\cdot\|_1 - \alpha \|\cdot\|_2)}(x') = \argmin_x \lambda \left(\|x\|_1 -\alpha \|x\|_2 \right)+ \frac{1}{2} \|x-x'\|_2^2.
\end{align}
\end{lemma}

\subsubsection{$z$-subproblem} \label{subsubsec:z_subprob}
\eqref{eq:z_sub} can be rewritten as
\begin{gather} \label{eq:z_sub2}
\begin{aligned}
    z^{t+1}&\in\argmin_z \sum_{j=1}^N \left[\frac{\beta_1}{2} \left \|u_j^{t+1} - \mathcal{F}(P_j^{t+1}z) + \frac{\Lambda_j^t}{\beta_1} \right\|_2^2 \right] + \frac{\beta_2}{2} \left\| v^{t+1} - \nabla z + \frac{y^t}{\beta_2} \right\|_2^2,
\end{aligned}
\end{gather}
which implies that $z^{t+1}$ must satisfy the first-order optimality condition
\begin{gather}
\begin{aligned} \label{eq:z_opt}
\left( \beta_1 \sum_{j=1}^N (P_j^{t+1})^{*} P_j^{t+1} - \beta_2 \Delta \right) z^{t+1} = \sum_{j=1}^N \beta_1 (P_j^{t+1})^{*} \mathcal{F}^{-1} \left( u_j^{t+1} + \frac{\Lambda_j^t}{\beta_1} \right) + \beta_2 \nabla^{\top}\left( v^{t+1} + \frac{y^t}{\beta_2} \right),
\end{aligned}
\end{gather}
where the Laplacian $\Delta =-\nabla^{\top} \nabla$.
Since the coefficient matrix of $z^{t+1}$ is invertible, solving \eqref{eq:z_opt} can be performed exactly, but it could be computationally expensive if the matrix system is extremely large because of the size of $z$. Since the coefficient matrix tends to be sparse, conjugate gradient \cite{hestenes1952methods} can be used to solve \eqref{eq:z_opt} like in \cite{chang2018total, chang2016phase}, but it needs access to all $N$ scan information and requires at most $n^2$ iterations to attain an exact solution, assuming exact arithmetic. Moreover, it is sensitive to roundoff error \cite{greenbaum1989behavior}.

Alternatively, we linearize \eqref{eq:z_sub2} to obtain the gradient descent step with step size $\delta_z^t > 0$:
\begin{align}\label{eq:z_gd}
    z^{t+1} = z^t - \delta_z^t \nabla_z \mathcal{L}(u^{t+1}, \omega^{t+1}, v^{t+1}, z^t, \Lambda^t, y^t).
\end{align}
Approximating $\nabla_z \mathcal{L}$ by its stochastic estimator $\tilde{\nabla}_z \mathcal{L}$ that only has access to $b \leq N$ scans, we have the SGD step:
\begin{align}\label{eq:z_sgd}
    z^{t+1} = z^t - \delta_z^t \tilde{\nabla}_z \mathcal{L}(u^{t+1}, \omega^{t+1}, v^{t+1}, z^t, \Lambda^t, y^t).
\end{align}
To design candidates for $\tilde{\nabla}_z \mathcal{L}$, we will use the following lemma:
\begin{lemma}\label{lemma:mask}
Let $S \in \mathbb{R}^{m^2 \times n^2}$. If $\mathbf{e}_i \in \text{ker}(S)$ for some index $i$, then for any $x \in \mathbb{C}^{m^2}$, we have $(S^{\top} x)_i = 0$.
\end{lemma}
\begin{proof}
We have
$(S^{\top}x)_i = \langle S^{\top}x, \mathbf{e}_i \rangle = \langle x, S\mathbf{e}_i \rangle = \langle x ,0 \rangle =0.$
\end{proof}
For brevity, we denote the vectors
\begin{gather}
\begin{aligned}\label{eq:matrix_def}
    A_j^t &= -\beta_1(P_j^{t+1})^* \left[\mathcal{F}^{-1}\left(u_j^{t+1} + \frac{\Lambda_j^t}{\beta_1} \right) - P_j^{t+1} z^t \right],\\
    B^t &= - \beta_2 \left[ \nabla^{\top} \left(v^{t+1} + \frac{y^t}{\beta_2} \right)+ \Delta z^t \right].
\end{aligned}
\end{gather}
At each pixel $i=1, \ldots, n^2$, \eqref{eq:z_gd} becomes
\begin{gather}\label{eq:pixel_gd}
\begin{aligned}
    (z^{t+1})_i &= (z^t)_i - \delta_z^t(\nabla_z \mathcal{L}(u^{t+1}, \omega^{t+1}, v^{t+1}, z^t, \Lambda^t, y^t))_i=(z^t)_i-\delta_z^t \left[\sum_{j=1}^N (A_j^t)_i + (B^t)_i \right].
\end{aligned}
\end{gather}
 By Lemma \ref{lemma:mask}, since $(P_j^{t+1})^* = S_j^{\top}(D_{\omega^{t+1}})^*$, we have $(A_j^t)_i = 0$ if $\mathbf{e}_i \in \text{ker}(S_j) \subset \text{ker}(D_{\omega^{t+1}}S_j)$, which means that pixel $i$ is not scanned by the mask matrix $S_j$. For each $i=1, \ldots, n^2$, we define $N_i = \{j : \mathbf{e}_i \not \in \text{ker}(S_j)\} \subset \{1, \ldots, N\}$ to be the set of indices corresponding to the mask matrices that scan pixel $i$.
As a result, \eqref{eq:pixel_gd} reduces to and can be rewritten as
\begin{gather} \label{eq:pixel_gd2}
    \begin{aligned}
    (z^{t+1})_i &= (z^t)_i - \delta_z^t \left[\sum_{j \in N_i} (A_j^t)_i + (B^t)_i \right] = (z^t)_i - |N_i|\delta_z^t \left[ \frac{1}{|N_i|}  \sum_{j \in N_i}\left( (A_j^t)_i + \frac{1}{|N_i|}(B^t)_i\right) \right].
    \end{aligned}
\end{gather}

Let $n_i^t \coloneqq \{j \in n^t: \mathbf{e}_i \not \in \text{ker}(S_j)\} \subset N_i$ be the set of indices corresponding to the mask matrices available at iteration $t$ that scan pixel $i$. Since the SGD estimator \cite{bottou2018optimization, robbins1951stochastic} of $\frac{1}{|N_i|}  \sum_{j \in N_i}\left( (A_j^t)_i + \frac{1}{|N_i|}(B^t)_i\right)$ is $\frac{1}{|n_i^t|} \sum_{j \in n_i^t} \left( (A_j^t)_i + \frac{1}{|N_i|}(B^t)_i\right)$, one candidate stochastic estimator (up to a constant multiple at each pixel) for $\nabla_{z} \mathcal{L}$ is
 $\tilde{\nabla}_z^{SGD} \mathcal{L}$, where
\begin{align}\label{eq:z_sgd_estimator}
    \left(\tilde{\nabla}_z^{SGD} \mathcal{L}(u^{t+1}, \omega^{t+1}, v^{t+1}, z^t, \Lambda^t, y^t) \right)_i=
    \begin{cases}
    \frac{1}{|n_i^t|} \displaystyle\sum_{j \in n_i^t} \left( (A_j^t)_i + \frac{1}{|N_i|}(B^t)_i\right), & \text{ if } n_i^t \neq \emptyset, \\
    0, &\text{ if } n_i^t = \emptyset
    \end{cases}
\end{align}
for $i = 1, \ldots, n^2$. 
Similar to the $\omega$-subproblem, we further extend \eqref{eq:z_sgd_estimator} to incorporate spatially varying illumination strength. 
Let
\begin{equation}
\label{eq:PIE_step}
\Psi_{i,j} = \displaystyle
    \frac{1}{(1-\gamma_z)\|P_j^{t+1}\mathbf{e}_i\|_1^2 + \gamma_z \|\omega^{t+1}\|_{\infty}^2},
\end{equation}
where $\gamma_z \in [0,1]$. Then we have the following stochastic estimator $\tilde{\nabla}_z^{PIE} \mathcal{L}(u^{t+1},\omega^{t+1}, v^{t+1}, z^t, \Lambda^t, y^t)$ as
\begin{align}\label{eq:sgd_z_PIE}
 \left(\tilde{\nabla}_z^{PIE} \mathcal{L}(u^{t+1},\omega^{t+1}, v^{t+1}, z^t, \Lambda^t, y^t) \right)_i=
 \begin{cases}\frac{1}{|n_i^t|} \displaystyle \sum_{j \in n_i^t} \Psi_{i,j}\left( (A_j^t)_i + \frac{1}{|N_i|}(B^t)_i \right), &\text{ if } n_i^t \neq \emptyset, \\
 0, &\text{ if } n_i^t = \emptyset
 \end{cases}
\end{align}
for $i=1,\ldots, n^2$.

The overall stochastic ADMM algorithm that solves \eqref{eq:AITV_blind_prob} is provided by Algorithm \ref{alg:stochastic_ADMM_blind}.
\begin{algorithm}
	\caption{Stochastic ADMM to solve \eqref{eq:AITV_blind_prob}}
	\label{alg:stochastic_ADMM_blind}\smallskip
	\textbf{Input: } set of scan information $\{(S_j, d_j)\}_{j=1}^N$; model parameters $\lambda > 0$, $\alpha \in [0,1]$; penalty parameters $\beta_1, \beta_2> 0$; sequence of step sizes $\{(\delta_{\omega}^t, \delta_{z}^t) \}_{t=0}^{\infty}$; batch size $b \leq N$; PIE factors $\gamma_{\omega}, \gamma_{z} \in [0,1]$.
	\begin{algorithmic}[1]
		\STATE Initialize $\omega^0$, $z^0$, $\{u_j^0\}_{j=1}^N = \{\Lambda_j^0\}_{j=1}^N$, $y^0 = \nabla z^0$.
		\FOR{$t=0$ to $T-1$}
		\STATE Uniformly sample without replacement the mini-batch $n^t \subset \{1, \ldots, N\}$ of batch size $b$.  
		\STATE Update $u_j^{t+1}$
  according to \eqref{eq:u_closed_form}.
  \IF{$\omega$ is unknown}
		\STATE Update $\omega^{t+1}= \omega^t - \delta_{\omega}^t \tilde{\nabla}_{\omega} \mathcal{L}(u^{t+1}, \omega^{t}, v^t, z^t, \Lambda^t, y^t)$.
		See \eqref{eq:omega_sgd_estimator} and \eqref{eq:omega_stochastic2} for a candidate $\tilde{\nabla}_{\omega} \mathcal{L}$.
  \ELSE
    \STATE $\omega^{t+1} = \omega^t$.
    \ENDIF
		\STATE Compute{
		\begin{align*}
		    (v^{t+1})_i = \text{prox}_{\frac{\lambda}{\beta_2} ( \|\cdot \|_1 - \alpha \|\cdot\|_2)} \left((\nabla z^t)_i - \frac{(y^t)_i}{\beta_2} \right), \; \forall i=1, \ldots,n^2.
		\end{align*}}
		See Lemma \ref{lemma:prox_l1l2}.
		\STATE Update $z^{t+1} = z^t - \delta_z^t \tilde{\nabla}_z \mathcal{L}(u^{t+1}, \omega^{t+1}, v^{t+1}, z^t, \Lambda^t, y^t)$.
		See \eqref{eq:z_sgd_estimator} and \eqref{eq:sgd_z_PIE} for a candidate $\tilde{\nabla}_z \mathcal{L}$.
		\STATE Compute {
		\begin{align*}\Lambda_j^{t+1} &= \begin{cases}\Lambda_j^t + \beta_1\left(u_j^{t+1} - \mathcal{F}(\omega^{t+1} \circ S_j z^{t+1}) \right),\; &\text{ if } j \in n^t, \\
  \Lambda_j^t, \; &\text{ if } j \not \in n^t,
  \end{cases}\; \forall j =1, \ldots, N,\\
		y^{t+1} &= y^t + \beta_2\left(v^{t+1} - \nabla z^{t+1} \right).
		\end{align*}}
		\ENDFOR
	\end{algorithmic}
	\textbf{Output: $\omega^* = \omega^{T}$,  $z^* = z^{T}$} 
\end{algorithm}

\section{Convergence Analysis} \label{sec:convergence}
We discuss the convergence of Algorithm \ref{alg:stochastic_ADMM_blind}. Although global convergence for ADMM can be established using Kurdyka-Łojasiewicz assumptions \cite{wang2019global}, the result does not apply for our models because they contain the gradient operator, which does not satisfy the necessary surjectivity assumption. Hence, we will prove up to subsequential convergence.  The convergence analysis is based on the analyses done in \cite{chang2018total,chang2016phase, wen2012alternating}, where under certain assumptions, they showed that the iterate subsequences of the ADMM algorithms converge to Karush-Kuhn-Tucker (KKT) points. To simplify notation, let $Z= (u, \omega, v, z)$ and $\Omega = (\Lambda, y)$. A KKT point $(Z^{\star}, \Omega^{\star})$ of the Lagrangian \eqref{eq:blind_Lagrangian} satisfies the KKT conditions given by
\begin{subequations}
    \begin{align}
        u_j^{\star} &= \mathcal{F}(\omega^{\star} \circ S_j z^{\star})\quad \text{for } j=1, \ldots, N,  \label{eq:lagrange_kkt1}\\
    v^{\star} &= \nabla z^{\star}, \label{eq:v_equal_nablaz}\\
    0 &\in \begin{cases}
    \partial|u_j^{\star}| \circ ( |u_j^{\star}| - \sqrt{d_j}) + \Lambda_j^{\star}, &\text{if AGM}, \\
    \partial|u_j^{\star}| \circ\left( |u_j^{\star}| - \displaystyle \frac{d_j}{|u_j^{\star}|} \right) + \Lambda_j^{\star}, &\text{if IPM,}
    \end{cases}\quad \text{for } j=1, \ldots, N, \label{eq:u_kkt}\\
    -\frac{y^{\star}}{\lambda} &\in \partial(\|v^{\star}\|_1 - \alpha \|v^{\star}\|_{2,1}), \label{eq:v_kkt}\\
      \nabla_{\omega}\mathcal{L}(Z^{\star}, \Omega^{\star})&=0, \label{eq:omega_kkt}\\ 
      \nabla_z \mathcal{L}(Z^{\star}, \Omega^{\star})&=0. \label{eq:z_kkt}
\end{align}
\end{subequations}

Let $\mathbb{E}_t$ denote the expectation conditioned on the past sequence of iterates $\{(Z^k, \Omega^k)\}_{k=0}^t$. More specifically, if $\mathcal{F}_t$ is the sigma algebra generated by the mini-batches $\{n^k\}_{k=0}^{t-1}$, then $\mathbb{E}_t[\;\cdot \;] \coloneqq \mathbb{E}[\;\cdot \; | \; \mathcal{F}_t]$. We impose the following assumption:
\begin{assumption} \label{assume:gradient_bound} Let $\{(Z^t, \Omega^t)\}_{t=1}^{\infty}$ be a sequence of iterates generated by Algorithm \ref{alg:stochastic_ADMM_blind}. For brevity, define $\mathcal{L}(\omega) \coloneqq \mathcal{L}(u^{t+1}, \omega, v^t, z^t, \Lambda^t, y^t)$ and $\mathcal{L}(z) \coloneqq \mathcal{L}(u^{t+1}, \omega^{t+1}, v^{t+1}, z, \Lambda^t, y^t)$ at each iteration $t$.
Suppose that at each iteration $t$, the stochastic gradient estimators \\$\tilde{\nabla}_\omega \mathcal{L} (\omega^t) \coloneqq \tilde{\nabla}_\omega \mathcal{L}(u^{t+1}, \omega^t, v^{t}, z^{t}, \Lambda^t, y^t)$ and $\tilde{\nabla}_z \mathcal{L} (z^t) \coloneqq \tilde{\nabla}_z \mathcal{L}(u^{t+1}, \omega^{t+1}, v^{t+1}, z^{t}, \Lambda^t, y^t)$ satisfy  the following:
\begin{enumerate}[label=(\alph*)]
\item Unbiased estimation:
\begin{align}\mathbb{E}_t\left[\tilde{\nabla}_\omega \mathcal{L} (\omega^t) \;\middle |\; u^{t+1}\right] &= \nabla_{\omega} \mathcal{L}(\omega^t),\\
\mathbb{E}_t\left[\tilde{\nabla}_z\mathcal{L} (z^t)\;\middle | \;u^{t+1}, \omega^{t+1}, v^{t+1}\right] &= \nabla_{z} \mathcal{L}(z^t).
\end{align}
\item Expected smoothness: there exist constants $A_1, A_2, A_3 > 0$ such that 
\begin{align}
    &\mathbb{E}_t \left[ \|\tilde{\nabla}_{\omega} \mathcal{L}(\omega^t)\|_2^2 \;\middle |\; u^{t+1}\right]  \leq A_1 \mathcal{L}(\omega^t) + (A_2+1)  \| \nabla_{\omega} \mathcal{L}(\omega^t)\|_2^2 +A_3, \label{eq:expect_smooth_omega}\\
    &\mathbb{E}_t \left[ \|\tilde{\nabla}_{z} \mathcal{L}(z^t)\|_2^2\;\middle |\; u^{t+1}, \omega^{t+1}, v^{t+1}\right]  \leq A_1 \mathcal{L}(z^t) + (A_2+1) \| \nabla_{z} \mathcal{L}(z^t)\|_2^2 + A_3.\label{eq:expect_smooth_z}
\end{align}
\end{enumerate}
\end{assumption}
\begin{remark}
    Assumption {\ref{assume:gradient_bound}}(a) is standard in the analysis of many stochastic optimization algorithms \mbox{\cite{bottou2018optimization, lan2020first}}. Assumption {\ref{assume:gradient_bound}}(b) was first proposed in {\cite{khaled2020better}} as a more general, weaker assumption on the second moment of the stochastic gradient than the other standard assumptions, such as bounded variance {\cite{lan2020first}} and relaxed growth condition {\cite{bottou2018optimization}}.\end{remark}
To prove the convergence of Algorithm \ref{alg:stochastic_ADMM_blind}, we require the following preliminary results. Under some conditions, Lemma {\ref{lemma:bounded_sequence}} bounds the iterates $\{(Z^t, \Omega^t)\}_{t=1}^{\infty}$ generated by Algorithm {\ref{alg:stochastic_ADMM_blind}} and the gradients $\{(\nabla_{\omega} \mathcal{L}(\omega^t), \nabla_z \mathcal{L}(z^t))\}_{t=1}^{\infty}$. Moreover, it shows that the stochastic gradients have bounded variance. Lemma \ref{lemma:sufficient_decrease} provides useful inequalities while Proposition {\ref{prop:finite_length}} establishes that the gradients $\{(\nabla_{\omega} \mathcal{L}(\omega^t), \nabla_z \mathcal{L}(z^t))\}_{t=1}^{\infty}$ subsequentially converge to zero. The convergence of Algorithm \ref{alg:stochastic_ADMM_blind} is finally established in Theorem \ref{thm:kkt}. All proofs are delayed to Appendix~\ref{sec:appendix}.
\begin{lemma} \label{lemma:bounded_sequence}
Let $\{(Z^t, \Omega^t)\}_{t=1}^{\infty}$ be a sequence of iterates generated by Algorithm \ref{alg:stochastic_ADMM_blind} that satisfies Assumption \ref{assume:gradient_bound}. Suppose that $\{(\omega^t, z^t)\}_{t=1}^{\infty}$ is bounded and  $\sum_{t=1}^{\infty} \|\Omega^{t+1} - \Omega^t\|_2^2 < \infty$. Then $\{(Z^{t}, \Omega^t)\}_{t=1}^{\infty}$ and $\{(\nabla_{\omega} \mathcal{L}(\omega^t), \nabla_z \mathcal{L}(z^t))\}_{t=1}^{\infty}$  are bounded and  there exists a constant $\sigma > 0$ such that
\begin{align}
\mathbb{E}_t \left[ \|\tilde{\nabla}_{\omega} \mathcal{L}(\omega^t)\|_2^2 \;\middle |\; u^{t+1}\right] - \| \nabla_{\omega} \mathcal{L}(\omega^t)\|_2^2  \leq \sigma^2, \label{eq:omega_var_bound} \\
\mathbb{E}_t \left[ \|\tilde{\nabla}_{z} \mathcal{L}(z^t)\|_2^2\;\middle |\; u^{t+1}, \omega^{t+1}, v^{t+1}\right] - \| \nabla_{z} \mathcal{L}(z^t)\|_2^2\leq \sigma^2. \label{eq:z_var_bound}
\end{align}
\end{lemma}
\begin{lemma}\label{lemma:sufficient_decrease}
Let $\{(Z^t, \Omega^t)\}_{t=1}^{\infty}$ be a sequence of iterates generated by Algorithm \ref{alg:stochastic_ADMM_blind} that satisfies Assumption \ref{assume:gradient_bound}. Suppose that $\{(\omega^t, z^t)\}_{t=1}^{\infty}$ is bounded and  $\sum_{t=1}^{\infty} \|\Omega^{t+1} - \Omega^t\|_2^2 < \infty$. For each iteration $t$, we have
\begin{align}
    \mathbb{E}[\mathcal{L}(\omega^{t+1})] -\mathbb{E}[\mathcal{L}(\omega^t)] &\leq \frac{-2\delta_{\omega}^t +L_{\omega} (\delta_{\omega}^t)^2}{2} \mathbb{E}\left[ \|\nabla_{\omega} \mathcal{L}(\omega^t)\|_2^2 \right]+ \frac{L_{\omega}(\delta_{\omega}^t)^2 \sigma^2}{2},\label{eq:lemma4_5_omega_ineq}\\
     \mathbb{E}[\mathcal{L}(z^{t+1})] - \mathbb{E}[\mathcal{L}(z^{t})]&\leq \frac{-2\delta_{z}^t +L_{z} (\delta_{z}^t)^2}{2} \mathbb{E}\left[\|\nabla_{z} \mathcal{L}(z^t)\|_2^2 \right]+  \frac{L_{z}(\delta_{z}^t)^2 \sigma^2}{2}\label{eq:lemma4_5_z_ineq}
\end{align}
for some constants $L_{\omega}, L_z > 0$.
\end{lemma}

\begin{proposition}\label{prop:finite_length}
Let $\{(Z^t, \Omega^t)\}_{t=1}^{\infty}$ be a sequence of iterates generated by Algorithm \ref{alg:stochastic_ADMM_blind} that satisfies Assumption \ref{assume:gradient_bound}. Suppose $\{(\omega^t, z^t)\}_{t=1}^{\infty}$ is bounded, $\sum_{t=1}^{\infty} \|\Omega^{t+1} - \Omega^t\|_2^2 < \infty$, and 
\begin{gather}
\begin{aligned} \label{eq:step_size_condition}
    &\delta_{\omega}^t < \frac{2}{L_{\omega}}\;, \;\sum_{t=1}^{\infty} \delta_{\omega}^t = \infty, \quad  \sum_{t=1}^{\infty} (\delta_{\omega}^t)^2 < \infty,\\
   &\delta_{z}^t < \frac{2}{L_z},\; \sum_{t=1}^{\infty} \delta_{z}^t = \infty, \quad  \sum_{t=1}^{\infty} (\delta_{z}^t)^2 < \infty.
\end{aligned}
\end{gather}
Then
\begin{align}
    \liminf_{t \rightarrow \infty} \mathbb{E} \left[ \left\| \nabla_{\omega} \mathcal{L}(\omega^t) \right\|_2^2 \right] = 0, \label{eq:liminf1}\\
    \liminf_{t \rightarrow \infty}\mathbb{E} \left[ \left\| \nabla_{z} \mathcal{L}(z^t) \right\|_2^2 \right] = 0. \label{eq:liminf2}
\end{align}
\end{proposition}
\begin{theorem}\label{thm:kkt} Let $\{(Z^t, \Omega^t)\}_{t=1}^{\infty}$ be generated by Algorithm \ref{alg:stochastic_ADMM_blind}.
Under the same assumption as Proposition \ref{prop:finite_length}, there exists a subsequence of $\{(Z^t, \Omega^t)\}_{t=1}^{\infty}$  whose accumulation point $(Z^{\star}, \Omega^{\star})$ is a.s.\ a KKT point of \eqref{eq:blind_Lagrangian} that satisfies \eqref{eq:lagrange_kkt1}-\eqref{eq:z_kkt}.
\end{theorem}

\begin{remark} Similar to \cite{chang2018total,wen2012alternating}, we impose the assumption that $\{(\omega^t, z^t)\}_{t=1}^{\infty}$ is bounded to assist with the convergence analysis. This assumption can be removed by imposing box constraints onto the primal variables $(\omega, z)$ as done in \cite{chang2019blind,chang2016phase, hesse2015proximal}, but we leave it as a future direction to develop a globally convergent algorithm to solve \eqref{eq:AITV_blind_prob} with these box constraints. 
We note that the requirement $\sum_{t=1}^{\infty} \|\Omega^{t+1} - \Omega^t\|_2^2 < \infty$ is rather strong, but similar assumption was made in other nonconvex ADMM algorithms \cite{ke2017alternating,liu2022entropy, tu2023new, xu2012alternating} that do not satisfy the necessary assumptions for global convergence \cite{wang2019global}. Lastly, the step size condition \eqref{eq:step_size_condition} is a standard assumption in proving convergence of stochastic algorithms \cite{bottou2018optimization, robbins1951stochastic}.
\end{remark}

\section{Numerical Results}\label{sec:experiment}
In this section, we evaluate the performance of Algorithm \ref{alg:stochastic_ADMM_blind} on two complex images presented in Figure \ref{fig:test_image}. The chip image (Figures \ref{fig:chip_mag}-\ref{fig:chip_phase}) has image size $348 \times 348$, while the cameraman/baboon image (Figures \ref{fig:cameraman_mag}-\ref{fig:cameraman_phase}) has image size $350 \times 350$. The probe size used for both images is $256 \times 256$, and the scanning patterns of the probes are shown in Figures \ref{fig:chip_scan_pattern},\ref{fig:cameraman_scan_pattern}. In total, we have $N=100$ measurements per image. The measurements $\{d_j\}_{j=1}^N$ are either corrupted by Gaussian noise or Poisson noise as follows:
\begin{align}\label{eq:noise_d_j}
d_j = \begin{cases}
    (|\mathcal{F}(P_j z)|+\mathfrak{N}(0,s^2 I_{m^2 \times m^2}))^2, &\text{ if noise is Gaussian},\\
    \text{Poisson}(|\mathcal{F}(P_jz_{\zeta})|^{2}), & \text{ if noise is Poisson,}
\end{cases}
\end{align}
where $\mathfrak{N}(0,s^2 I_{m^2 \times m^2})$ is a multivariate Gaussian distribution with zero mean and covariance matrix $s^2 I_{m^2 \times m^2}$  and $z_{\zeta} = \zeta z$ for some constant $\zeta > 0$. Note that Poisson noise is stronger when $\zeta$ is smaller.
\afterpage{
\begin{figure}[t!!!]
     \centering
     \begin{subfigure}[b]{0.225\textwidth}
         \centering
         \includegraphics[width=\textwidth]{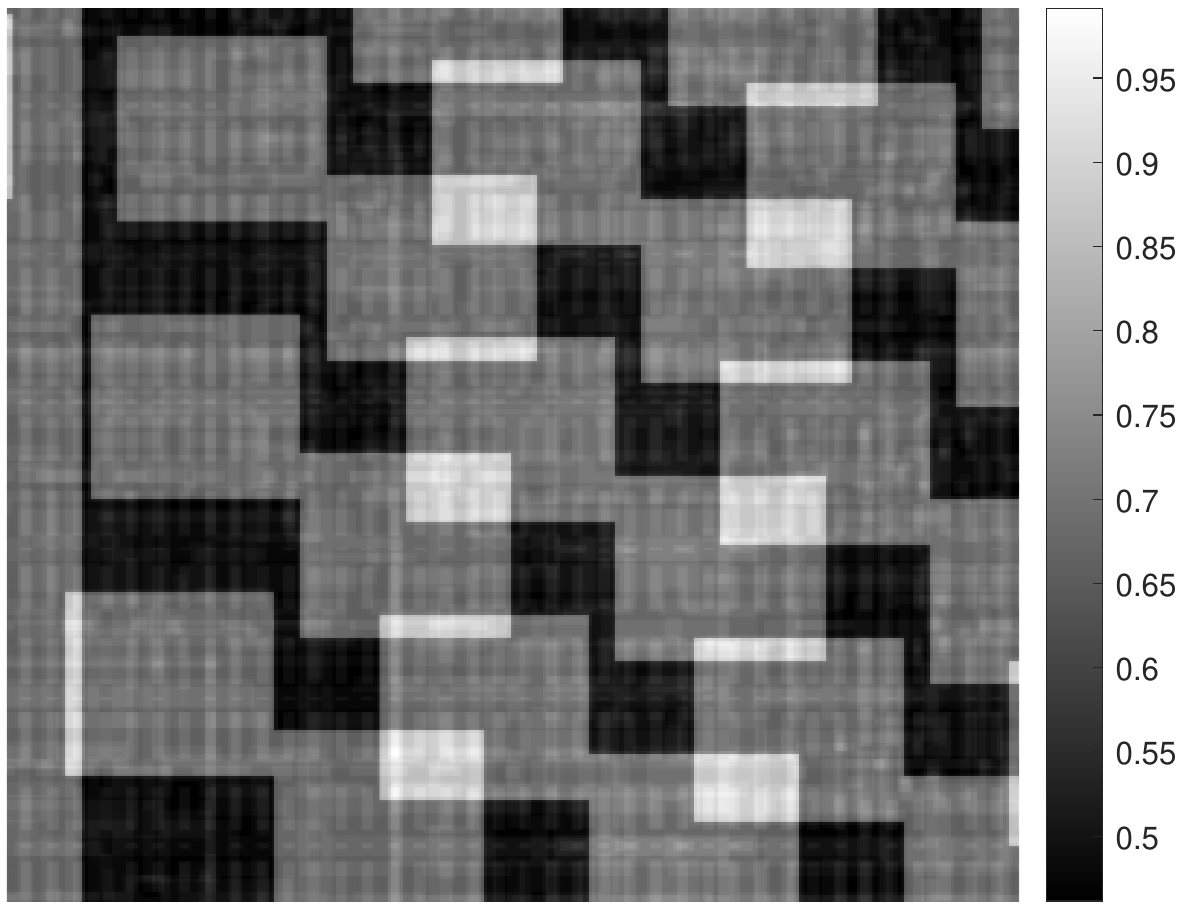}
         \caption{}
         \label{fig:chip_mag}
     \end{subfigure}
     \begin{subfigure}[b]{0.225\textwidth}
         \centering
         \includegraphics[width=\textwidth]{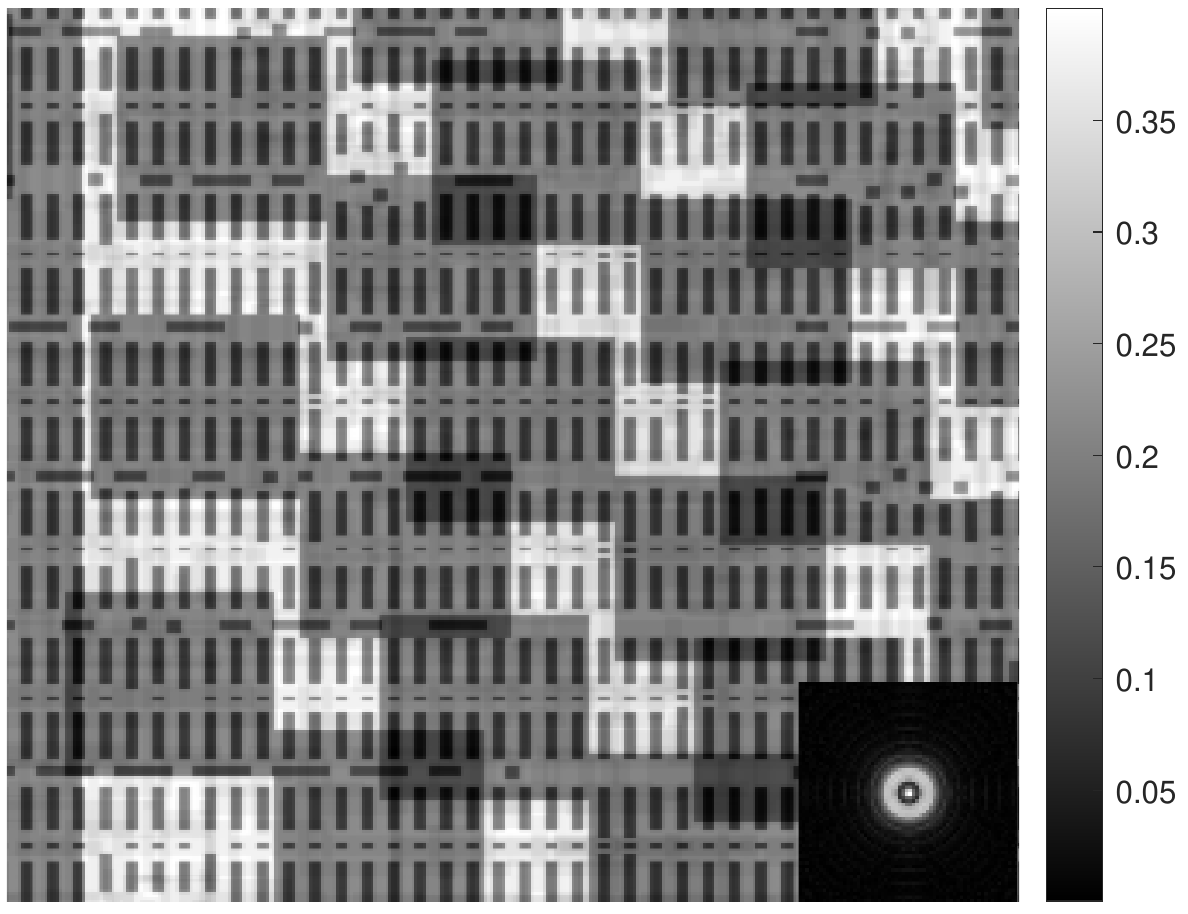}
         \caption{}
         \label{fig:chip_phase}
     \end{subfigure}
\begin{subfigure}[b]{0.225\textwidth}
         \centering
         \includegraphics[width=\textwidth]{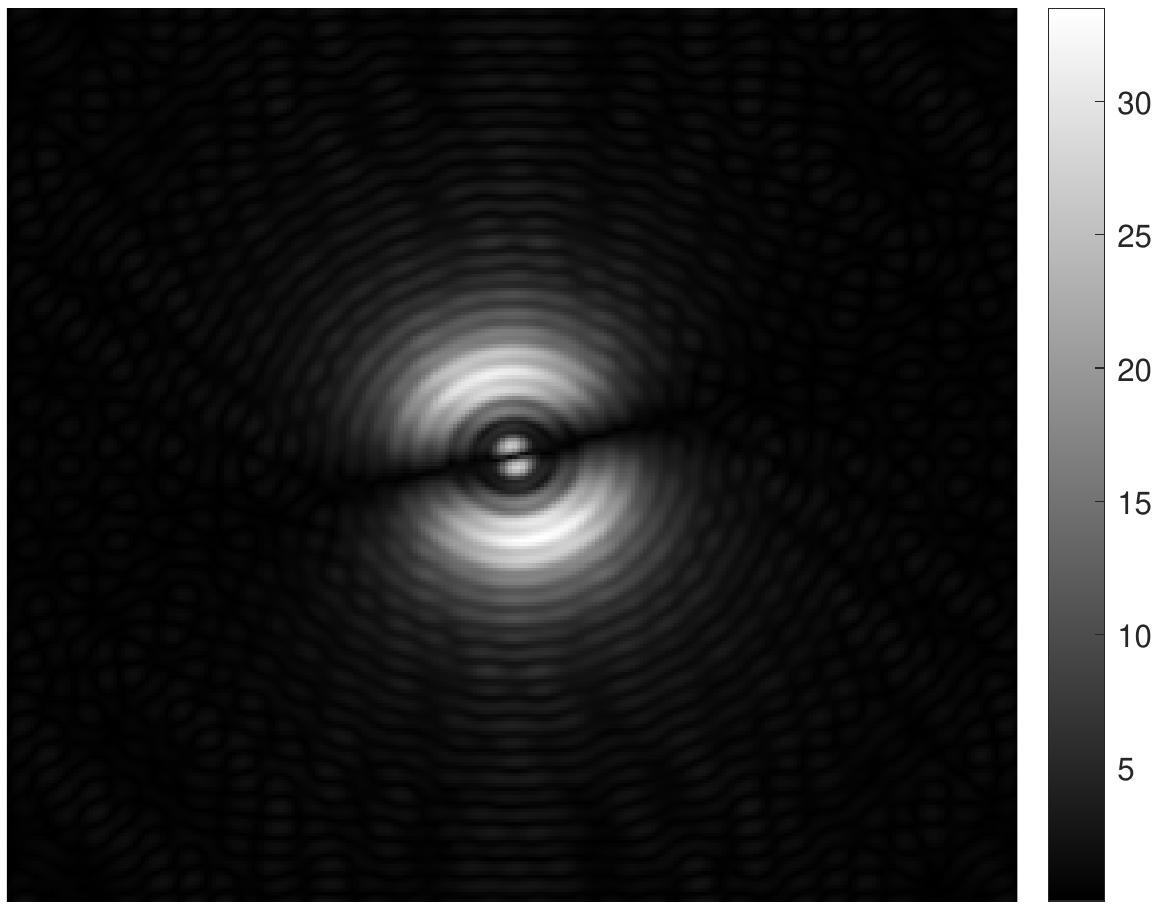}
         \caption{}
         \label{fig:chip_diff_probe}
     \end{subfigure}
     \begin{subfigure}[b]{0.225\textwidth}
         \centering
         \includegraphics[width=\textwidth]{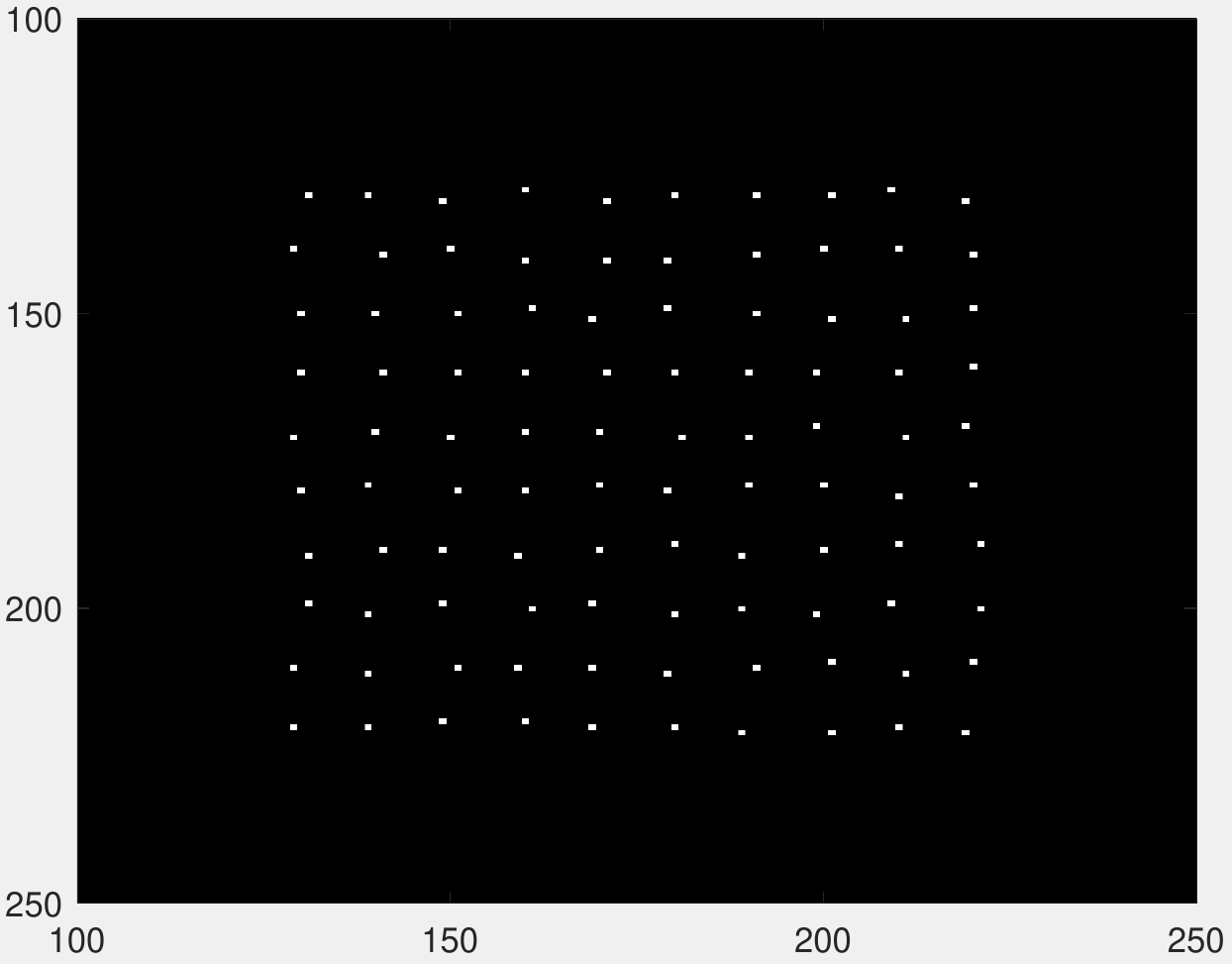}
         \caption{}
         \label{fig:chip_scan_pattern}
     \end{subfigure}\\
     \begin{subfigure}[b]{0.225\textwidth}
         \centering
         \includegraphics[width=\textwidth]{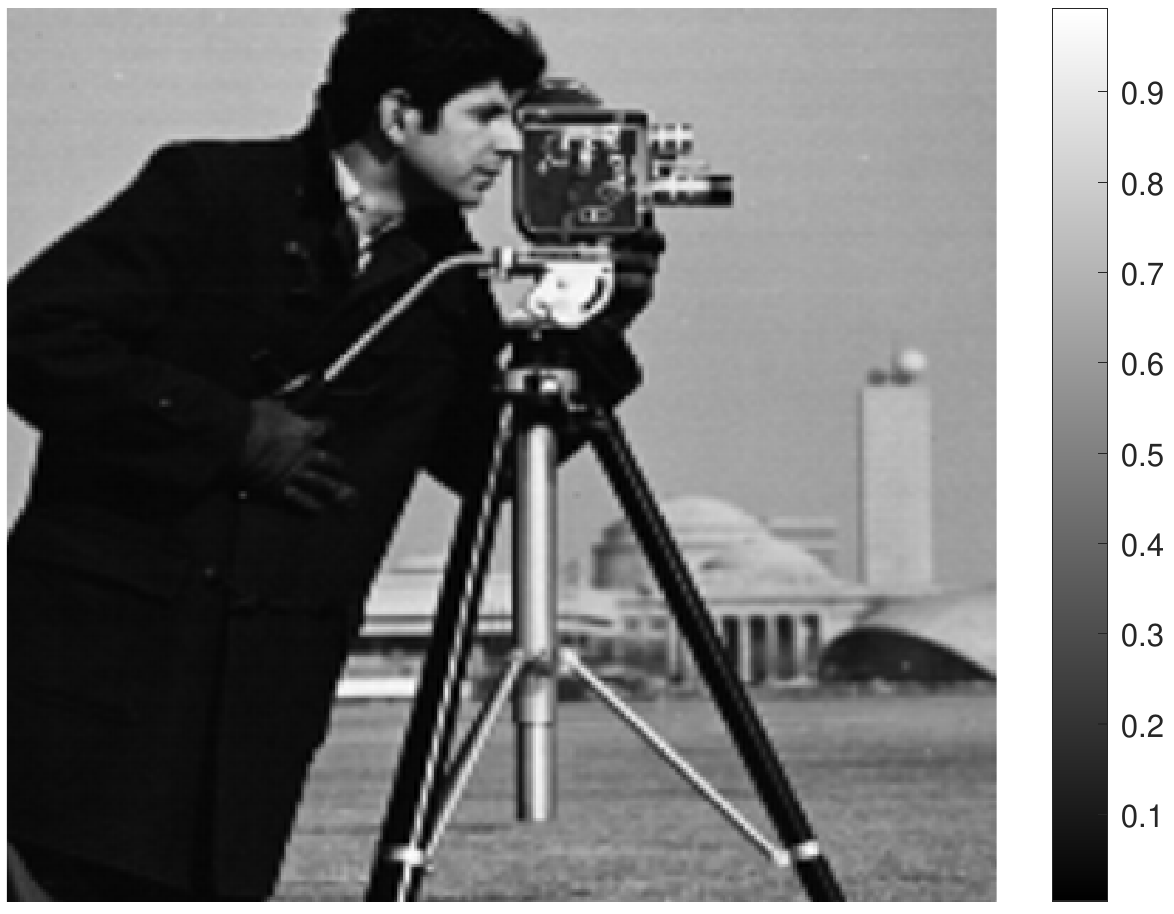}
         \caption{}
         \label{fig:cameraman_mag}
     \end{subfigure}
     \begin{subfigure}[b]{0.225\textwidth}
         \centering
         \includegraphics[width=\textwidth]{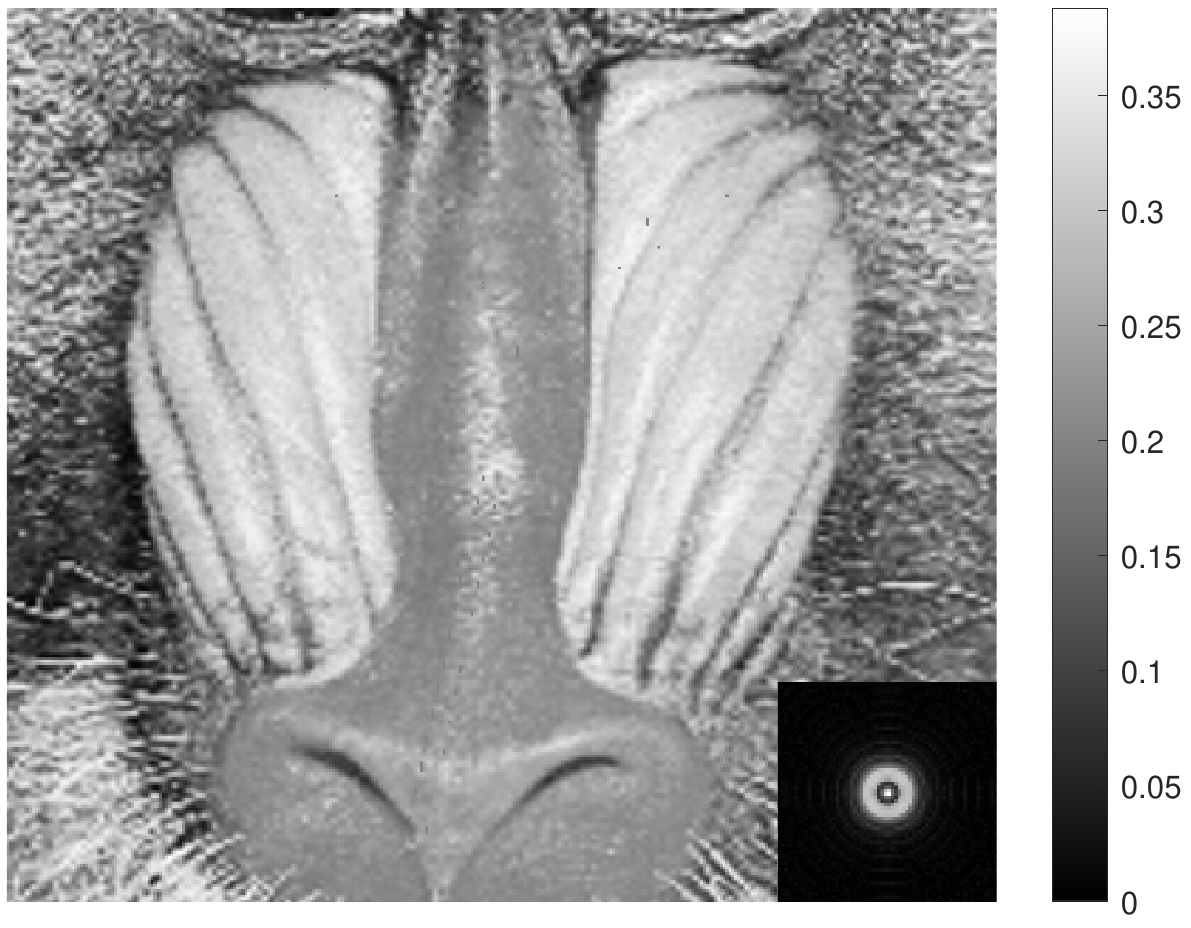}
         \caption{}
         \label{fig:cameraman_phase}
     \end{subfigure}
                    \begin{subfigure}[b]{0.225\textwidth}
         \centering
         \includegraphics[width=\textwidth]{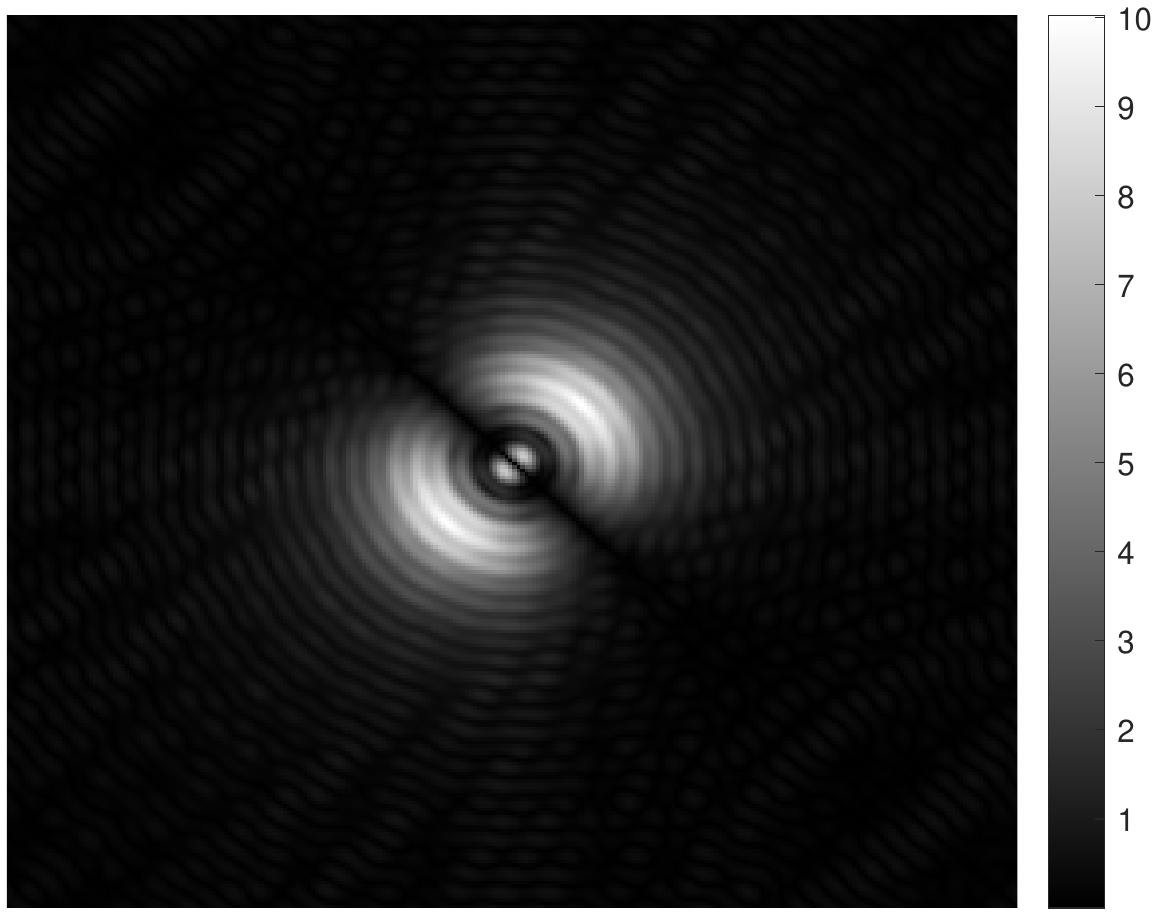}
         \caption{}
         \label{fig:cameraman_diff_probe}
     \end{subfigure}
        \begin{subfigure}[b]{0.225\textwidth}
         \centering
         \includegraphics[width=\textwidth]{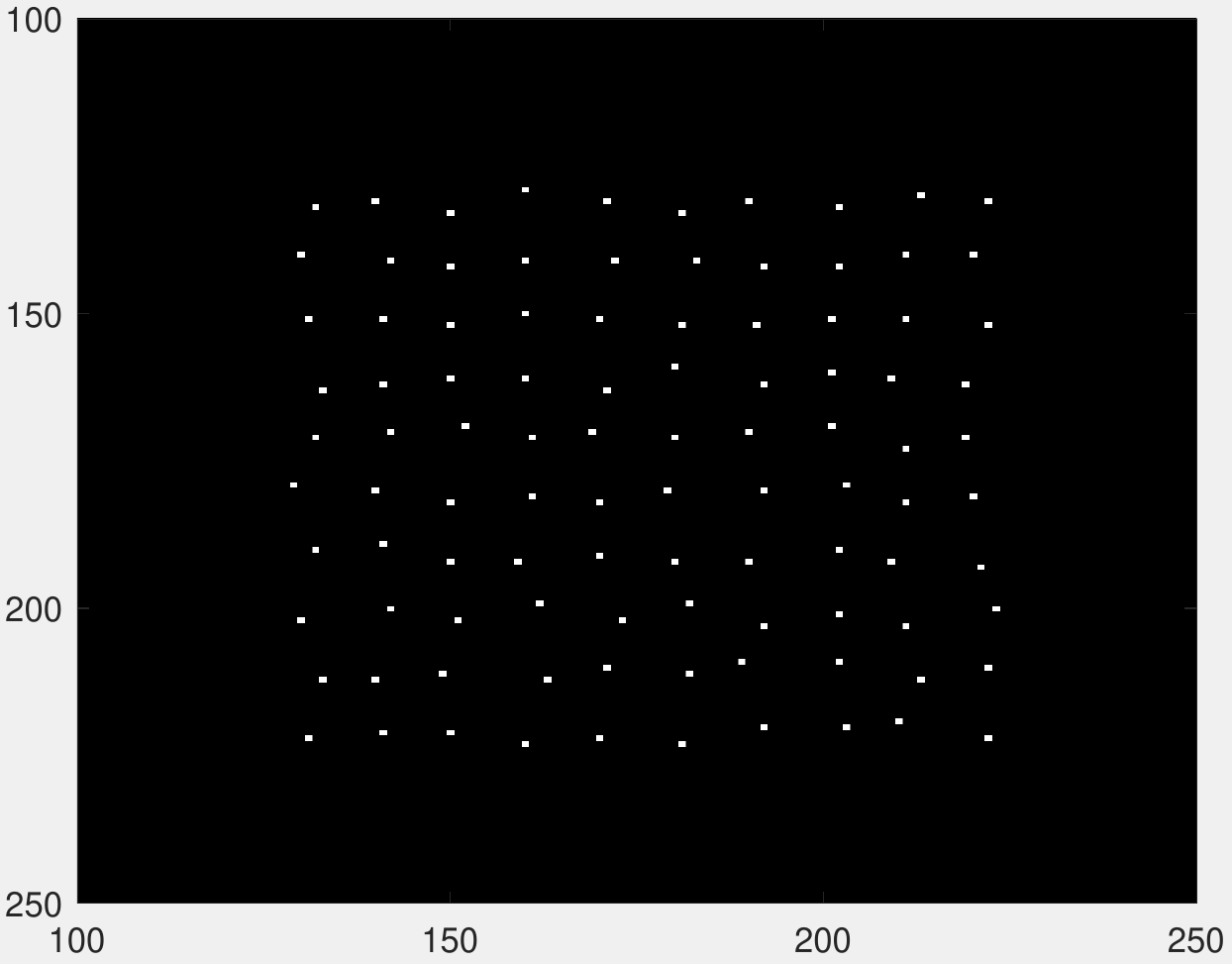}
         \caption{}
         \label{fig:cameraman_scan_pattern}
     \end{subfigure}
        \caption{Two complex sample images and their respective probes and scanning patterns examined in the experiments. First column: sample magnitude; second column: sample phase with inserted proportionally sized probe magnitude; third column: the magnitude differences between the ground-truth probe and the initial probe $\omega^0$; fourth column: scanning pattern of the probe zoomed in at $[100,250] \times [100,250]$, where the white dots represent the scanning lattice points  The image sizes are (a)-(b): $348 \times 348$; (d)-(e): $350 \times 350$. The probe sizes are $256 \times 256$.}
        \label{fig:test_image}
\end{figure}
\begin{table}[h!]
    \centering
    \caption{Parameter settings for each method. Note that $b$ refers to the batch size.}
    \label{tab:alg_param}
    \tiny
\begin{tabular}{|c|c|c|c|c|c|c|}
    \hline
        $\mathcal{B}(g,f)$& \makecell{Total \\ Epochs} &  $\beta_1 = \beta_2$ & $\delta_z^t$ & $\Psi_{j}^t$ & $\delta_{\omega}^t$ & $\Phi_{i,j}^t$ \\ [5ex]\hline
        \makecell{AGM} & 600 & \multirow{2}{*}{0.25} &  \multirow{2}{*}{\makecell{$\begin{cases}2 \sqrt{b} & \text{ if } 1 \leq t \leq 300, \\
        \frac{1}{5}\sqrt{b} &\text{ if } 300 < t \leq 450, \\
        \frac{1}{50} \sqrt{b} &\text{ if } 450 < t \leq 600
        \end{cases}$}}  & \multirow{2}{*}{$\gamma_z = 0.1$}&  \multirow{2}{*}{\makecell{$\begin{cases} \sqrt{b}\times 10^{-3} &\text{ if } 1 \leq t \leq 300, \\
        \sqrt{b} \times 10^{-4} &\text{ if } 300 < t \leq 450, \\
        \sqrt{b} \times 10^{-5} & \text{ if }450 \leq t \leq 600
        \end{cases}$}}& \multirow{2}{*}{$\gamma_{\omega} = 0.025$} \\ [2.5ex]\cline{1-2}
        \makecell{IPM} & 300 &  &   &  &  & \\[2.5ex]
        \hline
    \end{tabular}
\end{table}}%

For numerical evaluation, we compute the Structure Similarity Index Measure (SSIM) \cite{wang2004image} between the reconstructed image $z^{**}$ and the ground-truth image $z^g$ for the magnitude and phase, separately, where $z^{**}_{i} = \zeta^* z^*_{i+t^*}$ is adjusted for scaling by $\zeta^*$ and translation by $t^*$ and $(\zeta^*, t^*) = \displaystyle \argmin_{\zeta \in \mathbb{C}, t \in \mathbb{Z}} \sum_{i=1}^{n^2} |\zeta z^*_{i+t} - z^g_i|^2.$ 
We compare the proposed stochastic ADMM algorithms with its deterministic, full-batch counterparts (i.e., \eqref{eq:full_omega_subprob} and \eqref{eq:z_opt} are solved exactly) and its isotropic TV (isoTV) counterparts based on \cite{chang2018total, chang2016phase}. 
 The results are also compared with Douglas-Rachford splitting \cite{chang2019blind, thibault2009probe}, rPIE \cite{maiden2009improved}, and PHeBIE \cite{hesse2015proximal}. 

We initialize $z^0 = \frac{1}{\sqrt{2}} (\mathbf{1}+i\mathbf{1})$ when using AGM for Gaussian-corrupted measurements and $z^0 = \frac{\zeta}{\sqrt{2}} (\mathbf{1}+i\mathbf{1})$ when using IPM for Poisson-corrupted measurements. When performing the blind experiments using Algorithm \ref{alg:stochastic_ADMM_blind}, $\omega^0$ is initialized as the perturbation of the ground-truth probe. The magnitude differences between the initial and ground-truth probes are shown in Figures \ref{fig:chip_diff_probe},\ref{fig:cameraman_diff_probe}. 
The selected parameters, except for $\lambda$, are summarized in Table \ref{tab:alg_param}. The initial step sizes for $\delta_z^t$ and $\delta_{\omega}^t$ are determined empirically, and motivated by \eqref{eq:step_size_condition}, we decrease them by a factor of 10 at the $1/2$ and $3/4$ of the total epochs. Decreasing the step size in this way is a  popular technique, especially in the deep learning community \cite{goyal2017accurate, he2016deep}. Inspired from \cite{goyal2017accurate}, the step sizes are multiplied by a factor of $\sqrt{b}$ so that they scale by the batch size. For AITV regularization, we examine $\alpha \in \{0.2, 0.4, 0.6, 0.8\}$ and determine that $\alpha = 0.8$ yields the best results across all of our numerical examples. The batch sizes we examine are $b \in \{5, 10, 20, 50\}$ for Gaussian noise and $b \in \{5, 10, 20, 25\}$ for Poisson noise. For each parameter setting and image, we run three trials to obtain the mean SSIM values. 

All experiments are performed in MATLAB R2022b on a Dell laptop with a 1.80 GHz Intel Core i7-8565U
processor and 16.0 GB RAM. The code for the experiments is available at \url{https://github.com/kbui1993/Stochastic_ADMM_Ptycho}.
\subsection{Gaussian noise}
The SNR of the noisy measurements \cite{chang2018total} is given by
\begin{align*}
    \text{SNR}\left(\{\sqrt{d_j}\}_{j=1}^N, \{|\mathcal{F}(P_j z)|\}_{j=1}^N \right) = - 10 \log_{10} \left(  \frac{\displaystyle \sum_{j=1}^N \|\sqrt{d_j} -  |\mathcal{F}(P_jz)|\|_2^2}{\displaystyle \sum_{j=1}^N \|\mathcal{F}(P_j z)\|_2^2}  \right),
\end{align*}
so determined by the SNR, the noise level  $s$ in \eqref{eq:noise_d_j} is calculated by
\begin{align*}
    s = \sqrt{\frac{10^{-\text{SNR}/10} \displaystyle \sum_{j=1}^N \|\mathcal{F}(P_j z) \|_2^2}{Nm^2}}.
\end{align*}
 \afterpage{
\begin{table*}[t!!]
    \centering
    \caption{SSIM results of the algorithms applied to the Gaussian corrupted measurements with SNR = 40. The stochastic algorithms (e.g., AITV and isoTV, $b \in \{5, 10, 20, 50\}$) are ran three times to obtain the average SSIM values. \textbf{Bold} indicates best value; \underline{underline} indicates second best value.} 
    \label{tab:gaussian_result}
    \scriptsize
    \begin{tabular*}{\textwidth}{l|cc||cc|cc||cc}
    \hline
    & \multicolumn{4}{c|}{Non-blind} &  \multicolumn{4}{c}{Blind}\\ \hline
    &  \multicolumn{2}{c||}{Chip} &  \multicolumn{2}{c|}{Cameraman/Baboon}  &  \multicolumn{2}{c||}{Chip} &  \multicolumn{2}{c}{Cameraman/Baboon}  \\
        ~ & \makecell{mag.\\ SSIM} & \makecell{phase\\ SSIM} & \makecell{mag.\\ SSIM} & \makecell{phase\\ SSIM} &  \makecell{mag.\\ SSIM} & \makecell{phase\\ SSIM}& \makecell{mag.\\ SSIM} & \makecell{phase\\ SSIM} \\ \hline
        DR & 0.8130 & 0.8089 & 0.8701 & 0.5191 & 0.8008 & 0.7642 & 0.8009 & 0.3207 \\\hline
        rPIE & 0.8886 & 0.9073 & 0.8930 & 0.6055 & 0.9070 & 0.9120 & 0.8890 & 0.6145 \\ \hline
        PHeBIE & 0.8004 & 0.8019 & 0.8725  & 0.5718 & 0.8612 & 0.8438 &  0.8846  & 0.5756 \\\hline
        isoTV $(b=5)$ & 0.9501 & 0.9027 & 0.9393 & 0.7578 &	0.9426	& 0.8919 &0.9324 &0.7547
 \\ \hline
        isoTV $(b=10)$ & 0.9498 & 0.9004 & 0.9387 & 0.7475 &	0.9429 & 0.8891 & 0.9326 & 0.7477 \\ \hline
        isoTV $(b=20)$ & 0.9514	& 0.8981 &0.9385 &0.7302&0.9447	&0.8850&	0.9298&	0.7289
 \\ \hline
        isoTV $(b=50)$ & 0.9355 & 0.9193&	0.9294&	0.7050&	0.9322	&0.9047&	0.9153&	0.7025
 \\ \hline
        isoTV (full batch) & 0.9578 & 0.9145 & \underline{0.9769} & 0.7338 & 0.9527 & 0.8698 & 0.9589 & 0.5774\\ \hline
        
        AITV $(b=5)$ &0.9585 & 0.9556	& 0.9438 & \underline{0.7720}	&0.9490	&\underline{0.9477}	&0.9373	&\textbf{0.7775}
 \\ \hline
        AITV $(b=10)$ & 0.9620	&\underline{0.9579}	&0.9515	&\textbf{0.7747}&	0.9534	&\textbf{0.9481}&	0.9450	&\underline{0.7772}
 \\ \hline
        AITV $(b=20)$ & \underline{0.9629}&	\textbf{0.9583}&	0.9538&	0.7707&	\underline{0.9547}	&0.9470&	\underline{0.9468}&	0.7690
 \\ \hline
        AITV $(b=50)$ & 0.9585&	0.9550&	0.9490&	0.7358&	0.9514	&0.9432&	0.9391 &	0.7342
 \\ \hline
        AITV (full batch) & \textbf{0.9674} & 0.9513 & \textbf{0.9814} & 0.7463 & \textbf{0.9676} & 0.9296 & \textbf{0.9725} & 0.5956 \\  \hline
    \end{tabular*}
\end{table*}
\begin{figure}[t]
	\begin{minipage}{\linewidth}
		\centering
		\resizebox{\textwidth}{!}{%
			\begin{tabular}{c@{}c@{}c@{}c@{}c@{}c}
				\subcaptionbox{isoTV ($b=50$)}{\includegraphics[width = 1.35in]{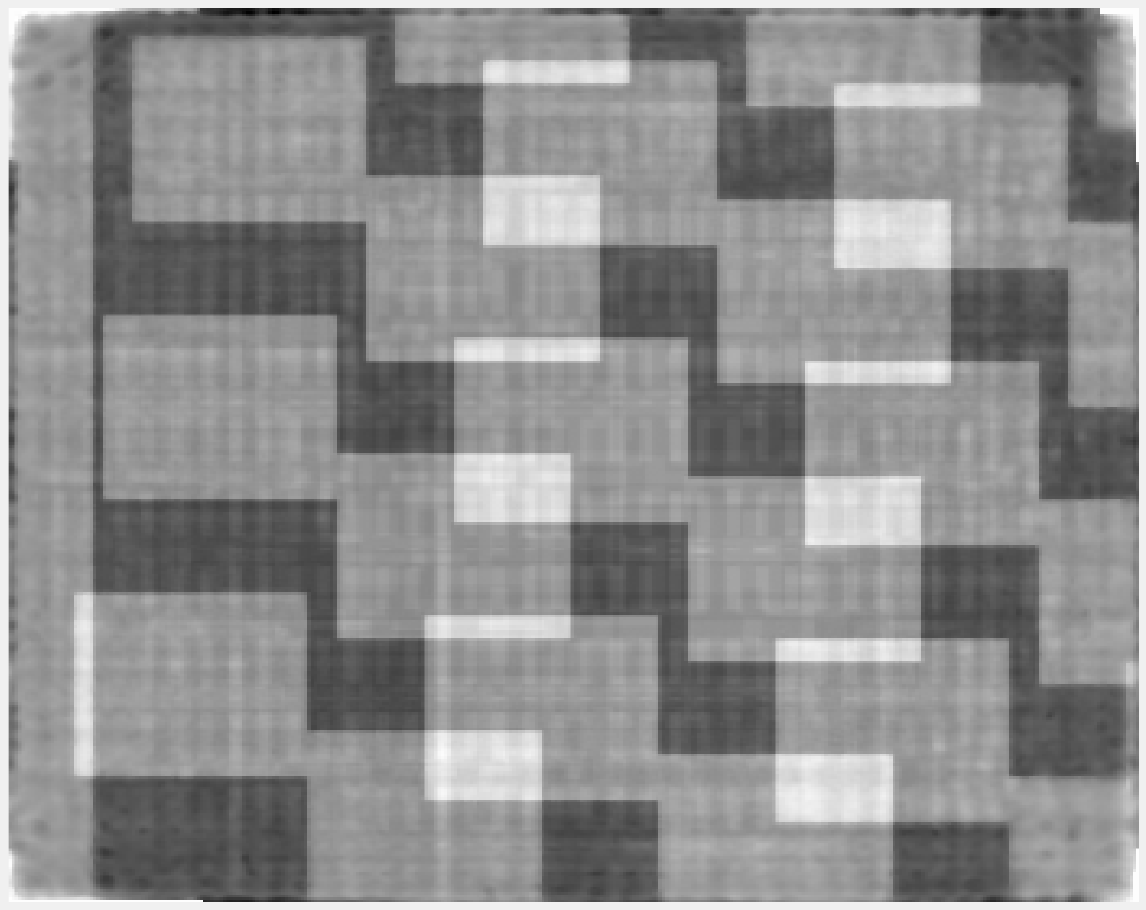}} & \subcaptionbox{AITV ($b=20$)}{\includegraphics[width = 1.35in]{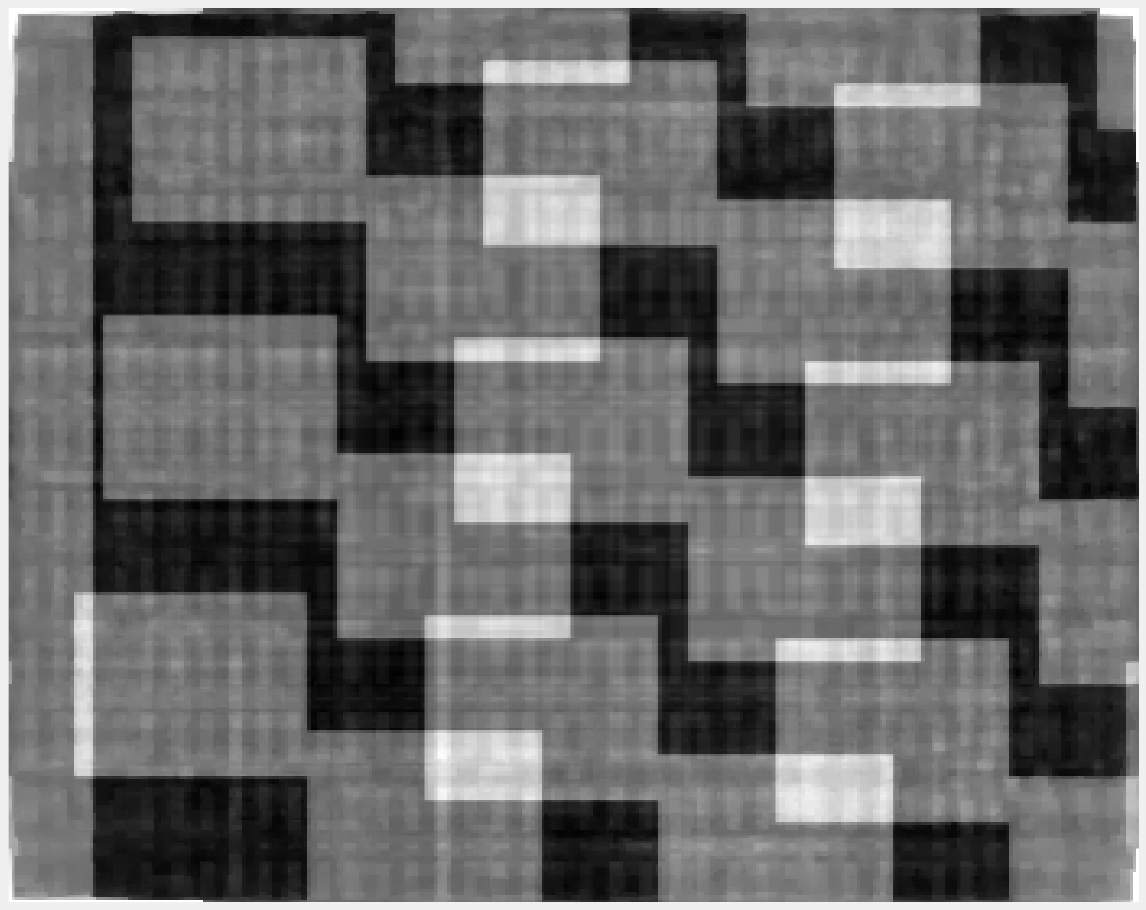}}   & \subcaptionbox{AITV (full)}{\includegraphics[width = 1.35in]{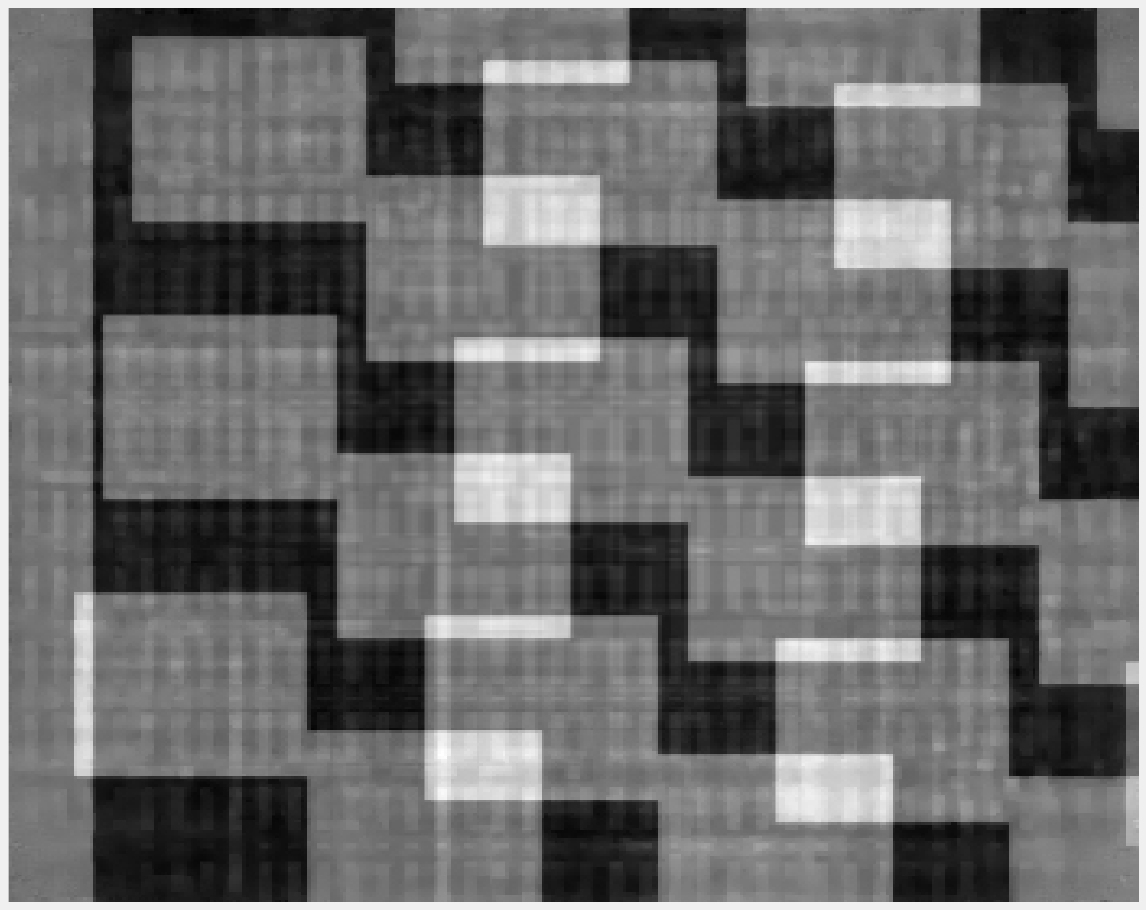}}  & \subcaptionbox{DR}{\includegraphics[width = 1.35in]{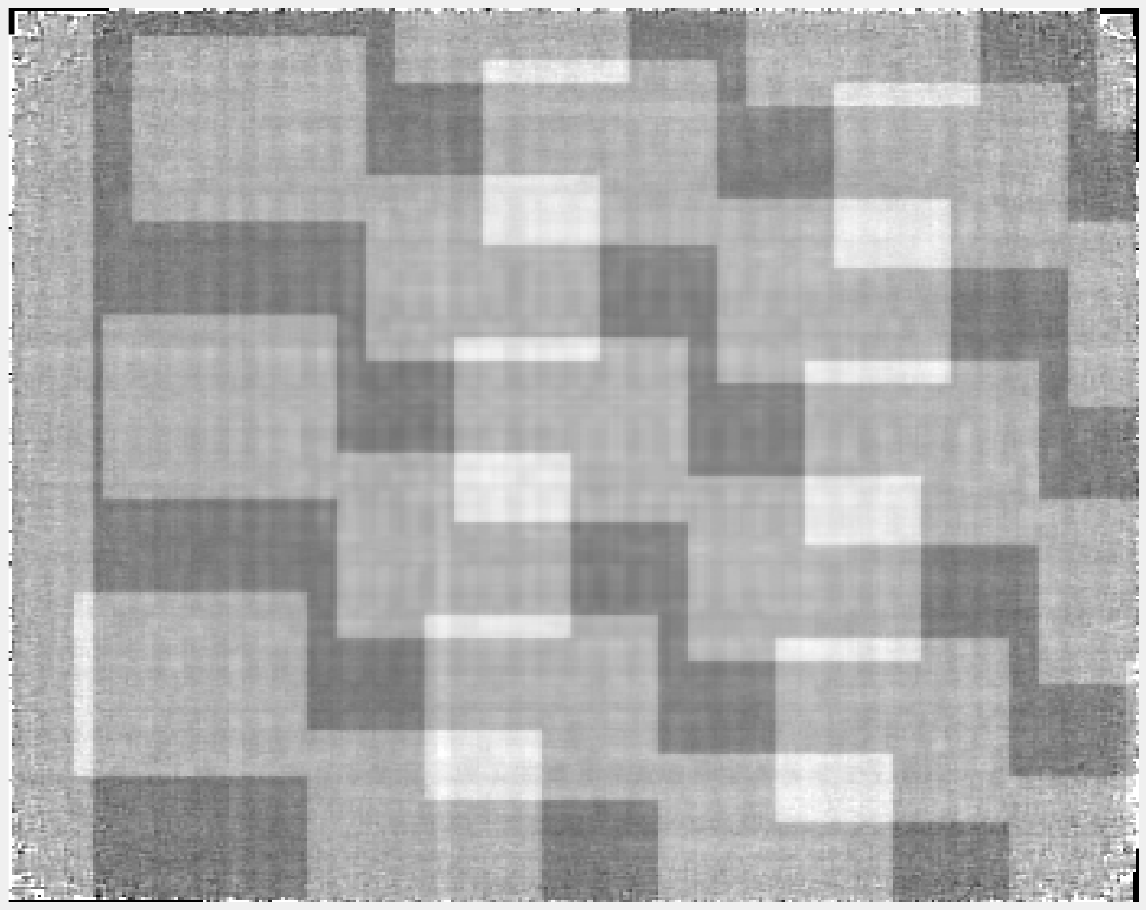}} & \subcaptionbox{rPIE}{\includegraphics[width = 1.35in]{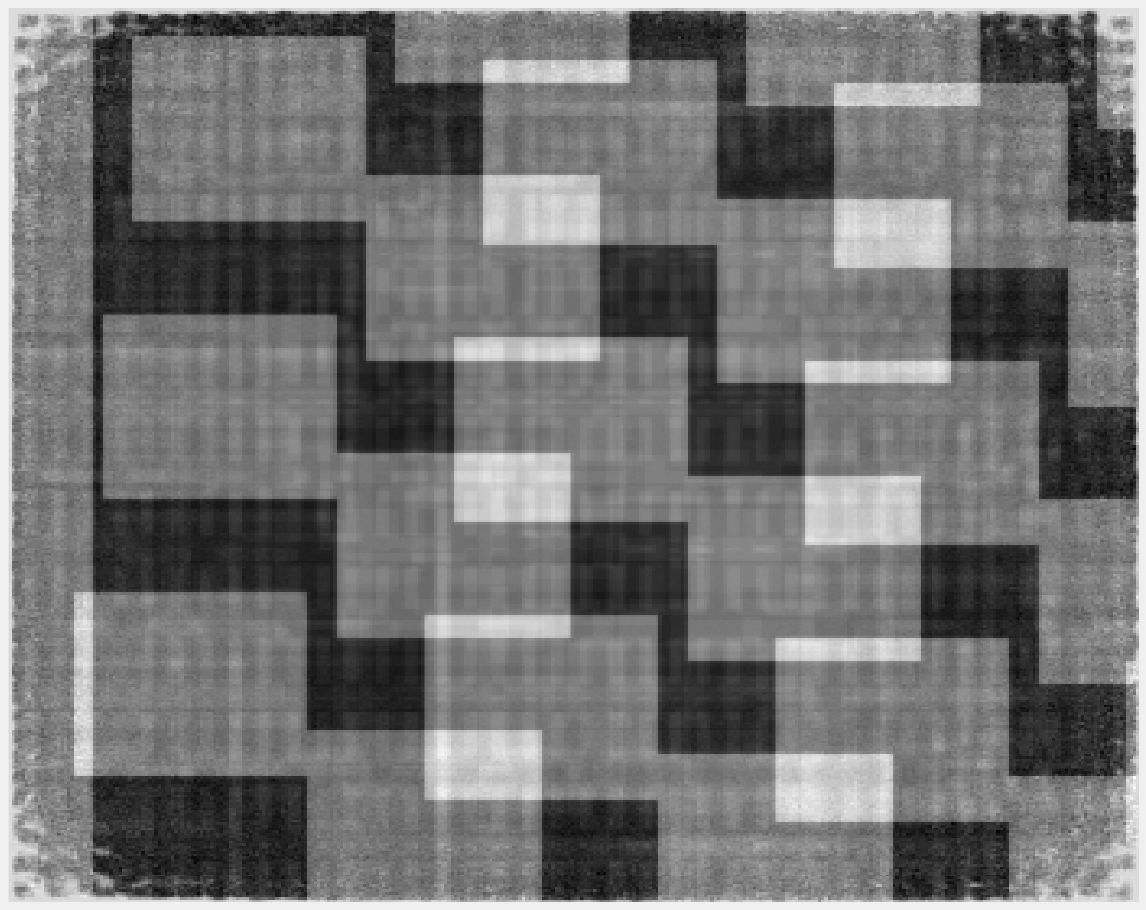}}& \subcaptionbox{PHeBIE}{\includegraphics[width = 1.35in]{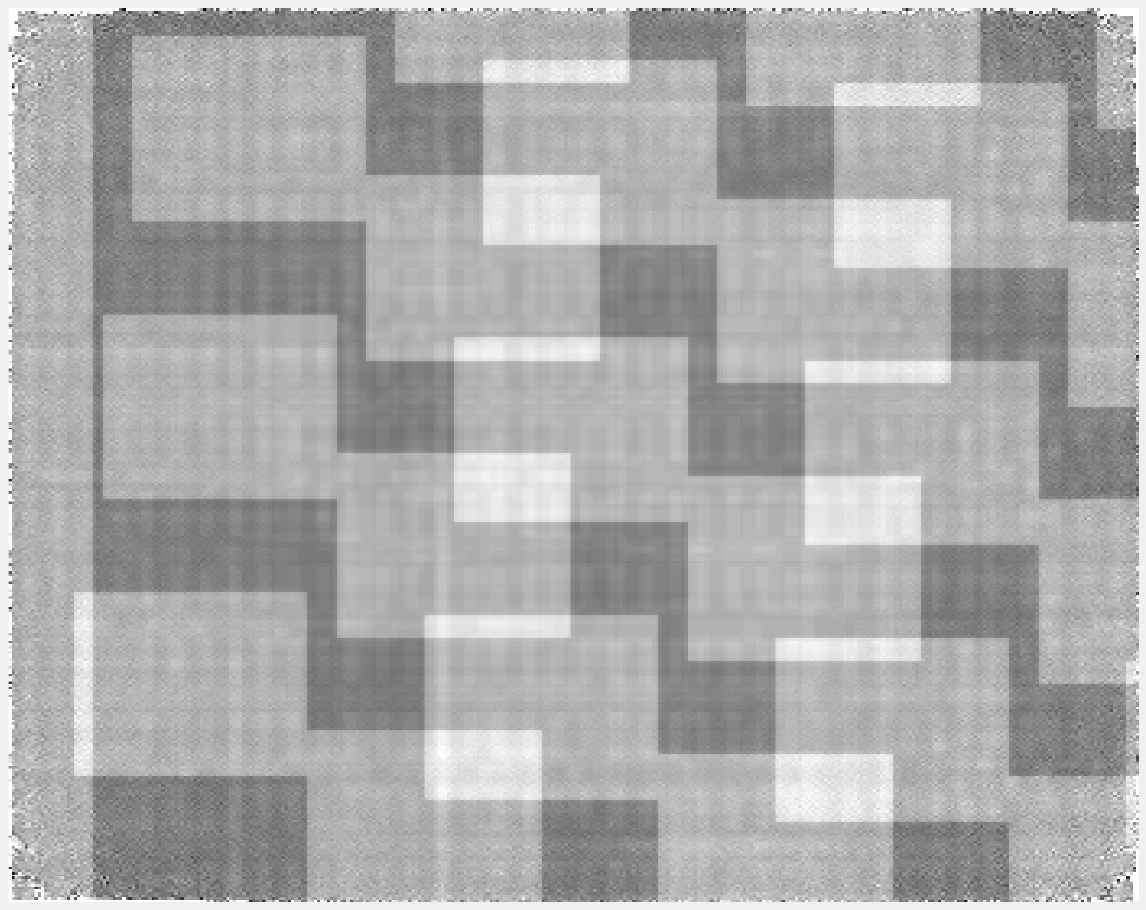}} \\	\subcaptionbox{isoTV $(b=50)$}{\includegraphics[width = 1.35in]{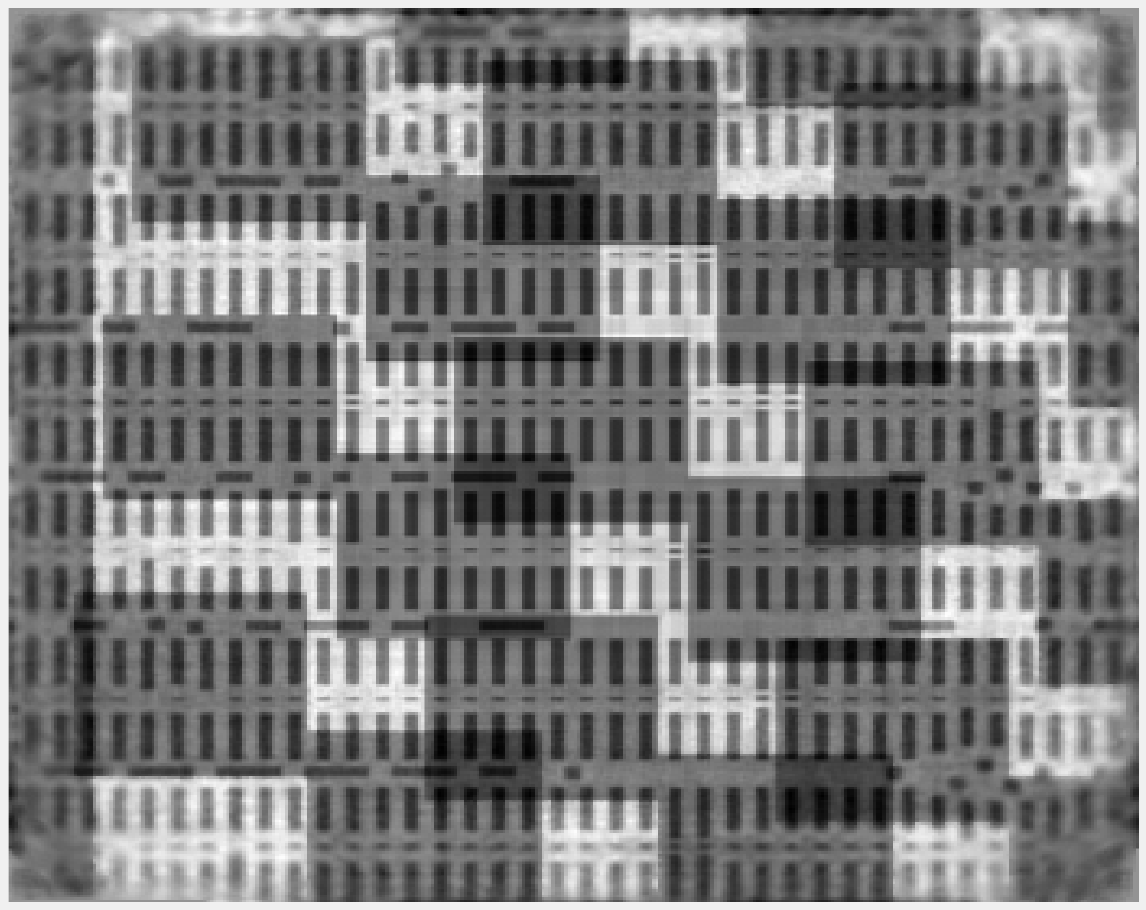}} & \subcaptionbox{AITV ($b=20$)}{\includegraphics[width = 1.35in]{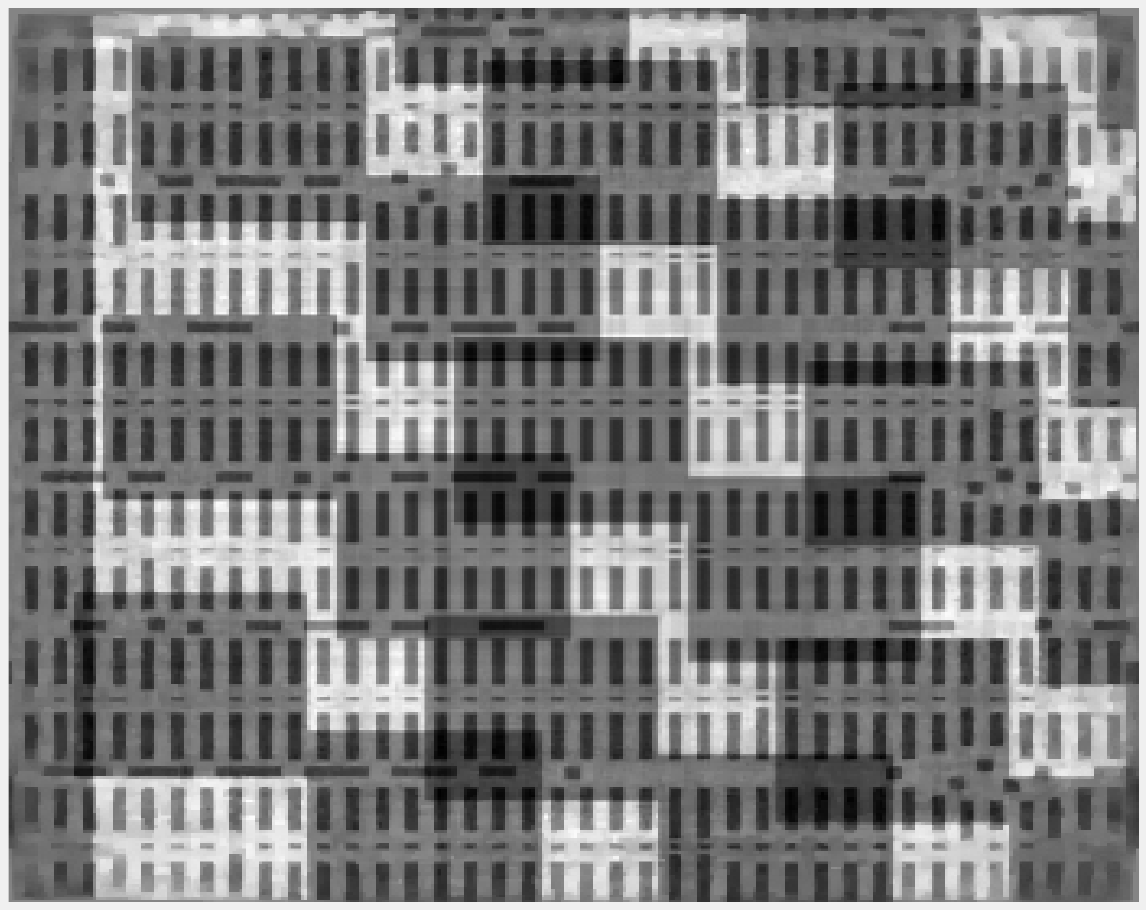}}   & \subcaptionbox{AITV (full)}{\includegraphics[width = 1.35in]{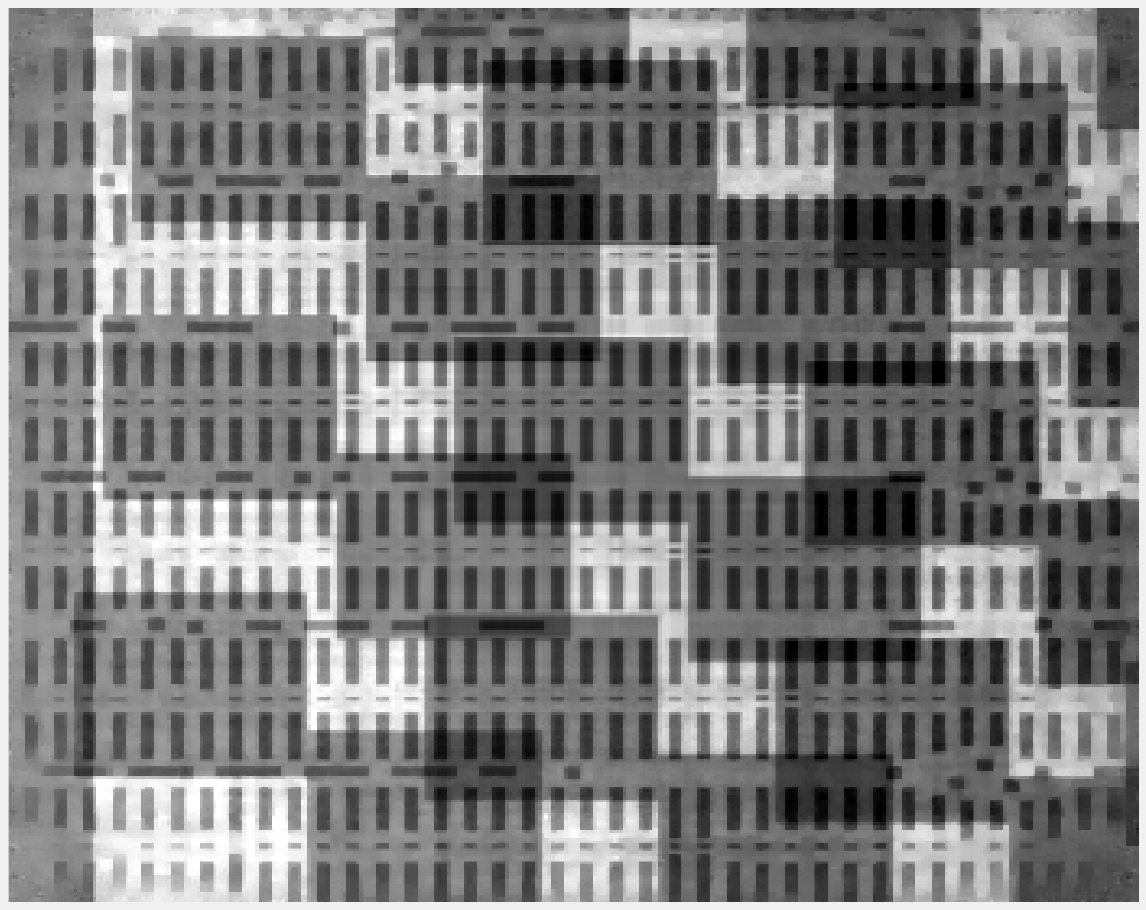}} &  \subcaptionbox{DR}{\includegraphics[width = 1.35in]{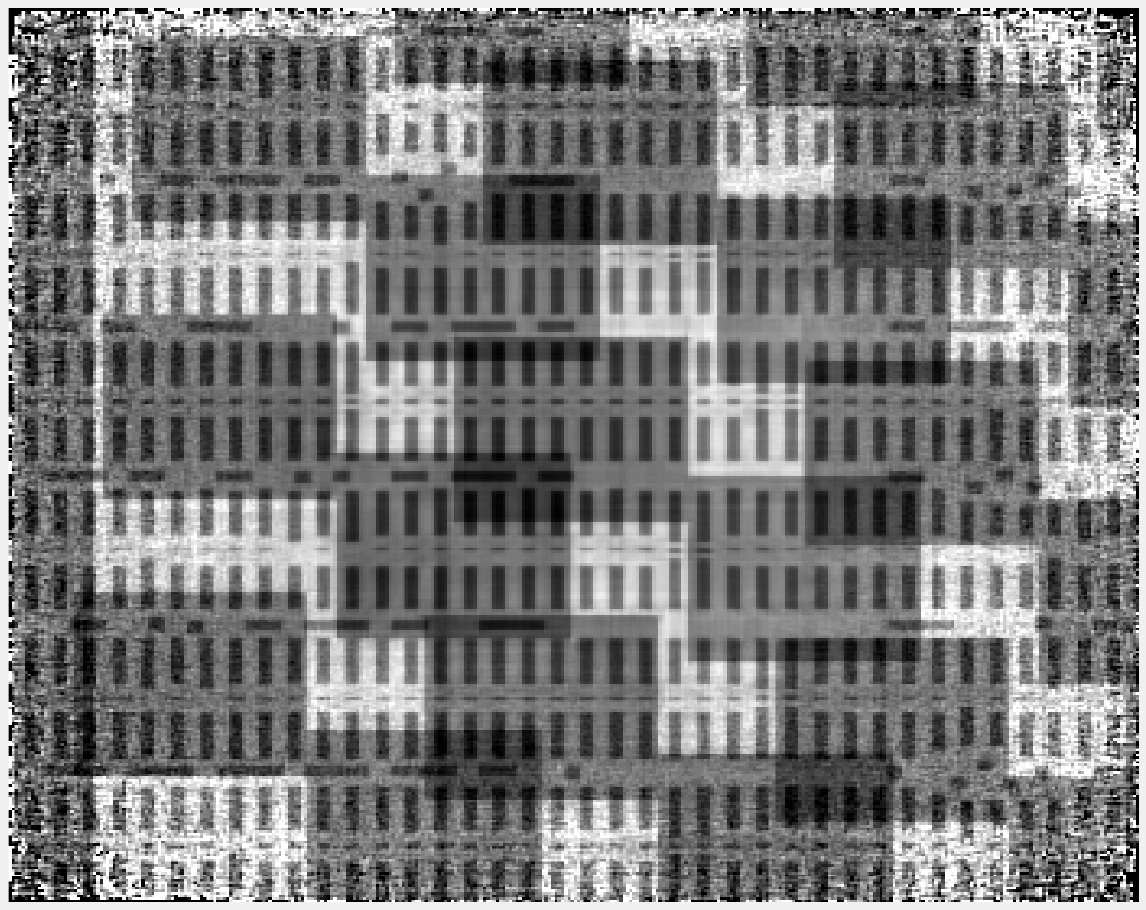}} & \subcaptionbox{rPIE}{\includegraphics[width = 1.35in]{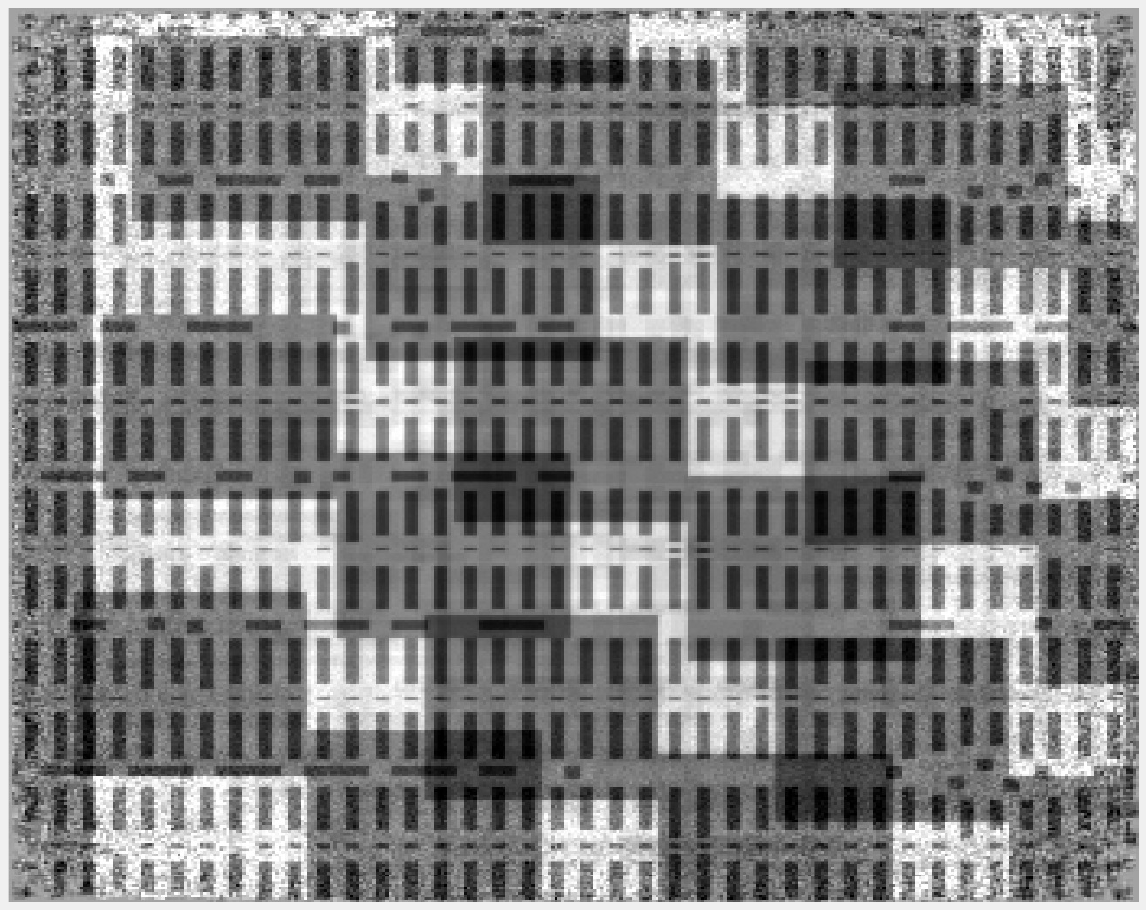}}& \subcaptionbox{PHeBIE}{\includegraphics[width = 1.35in]{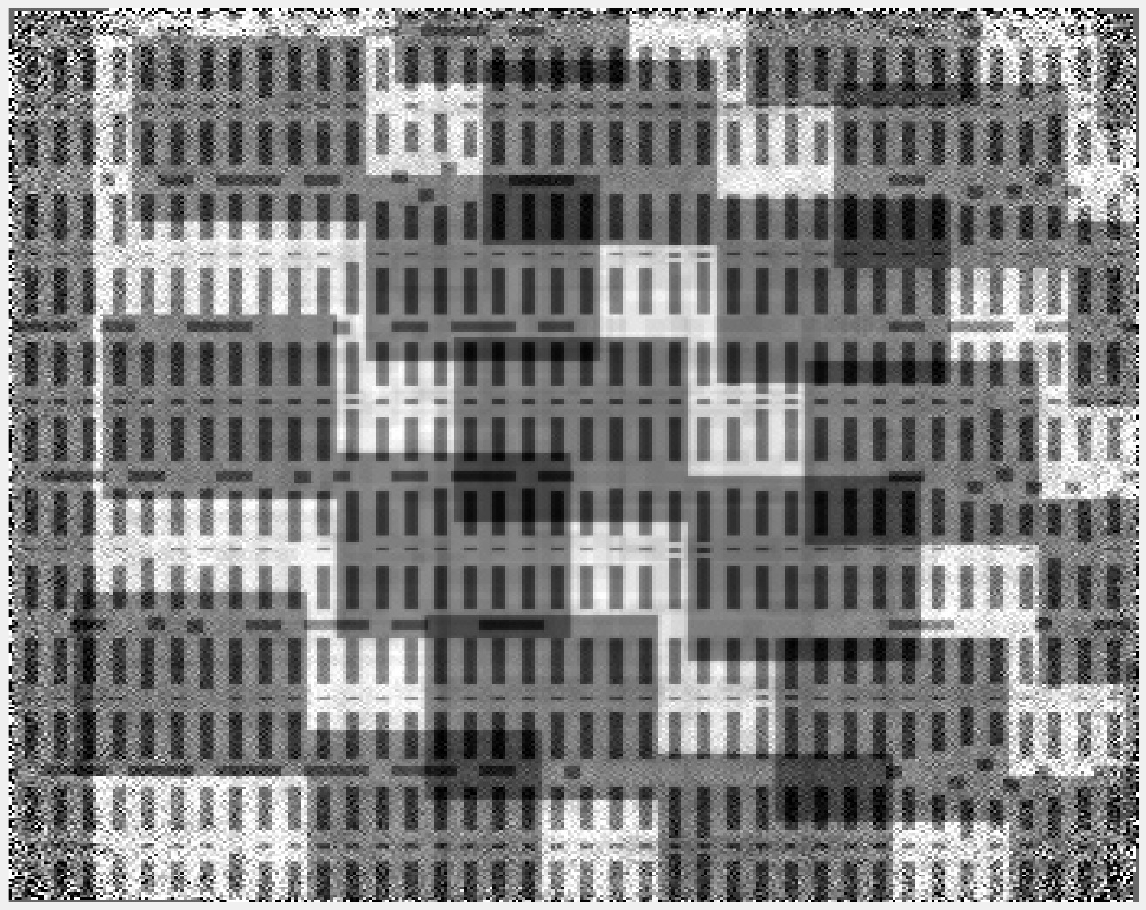}}
				\\
	\subcaptionbox{isoTV ($b=5$)}{\includegraphics[width = 1.35in]{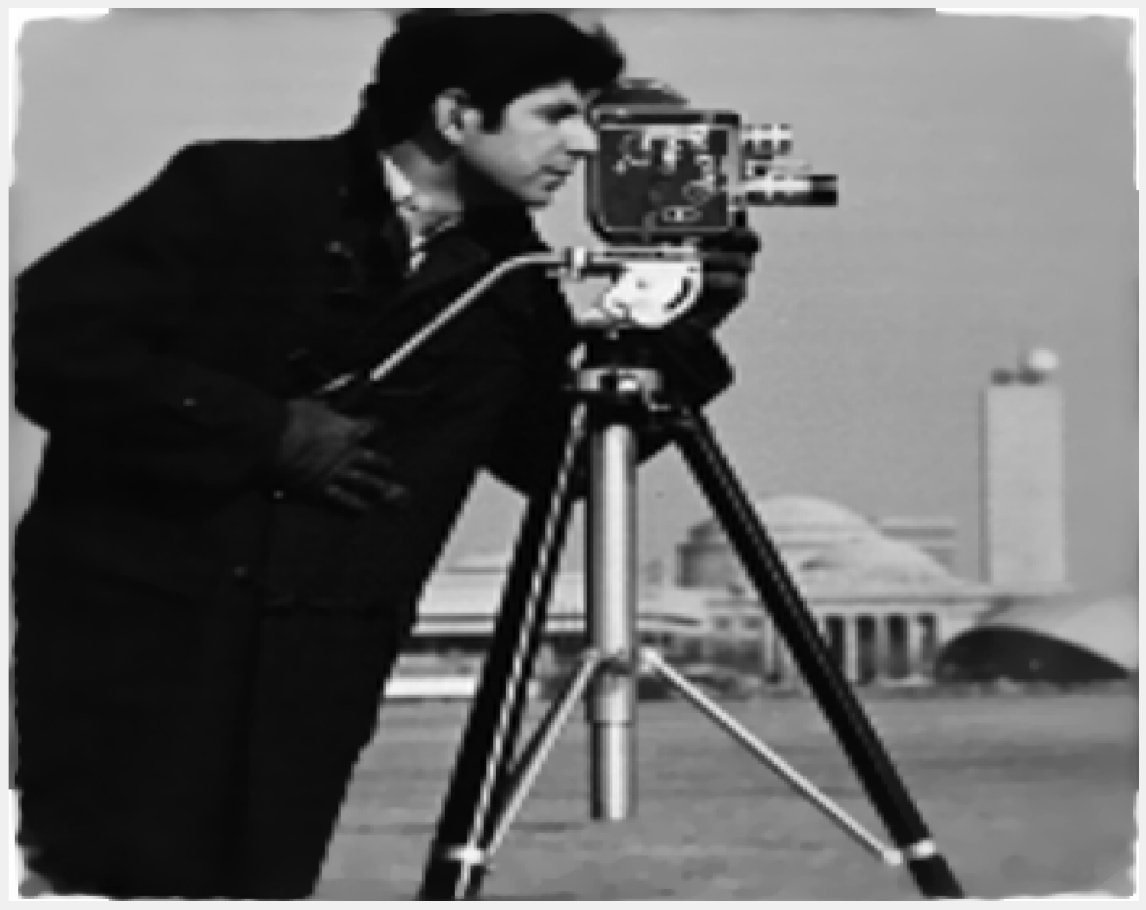}} & \subcaptionbox{AITV ($b=10$)}{\includegraphics[width = 1.35in]{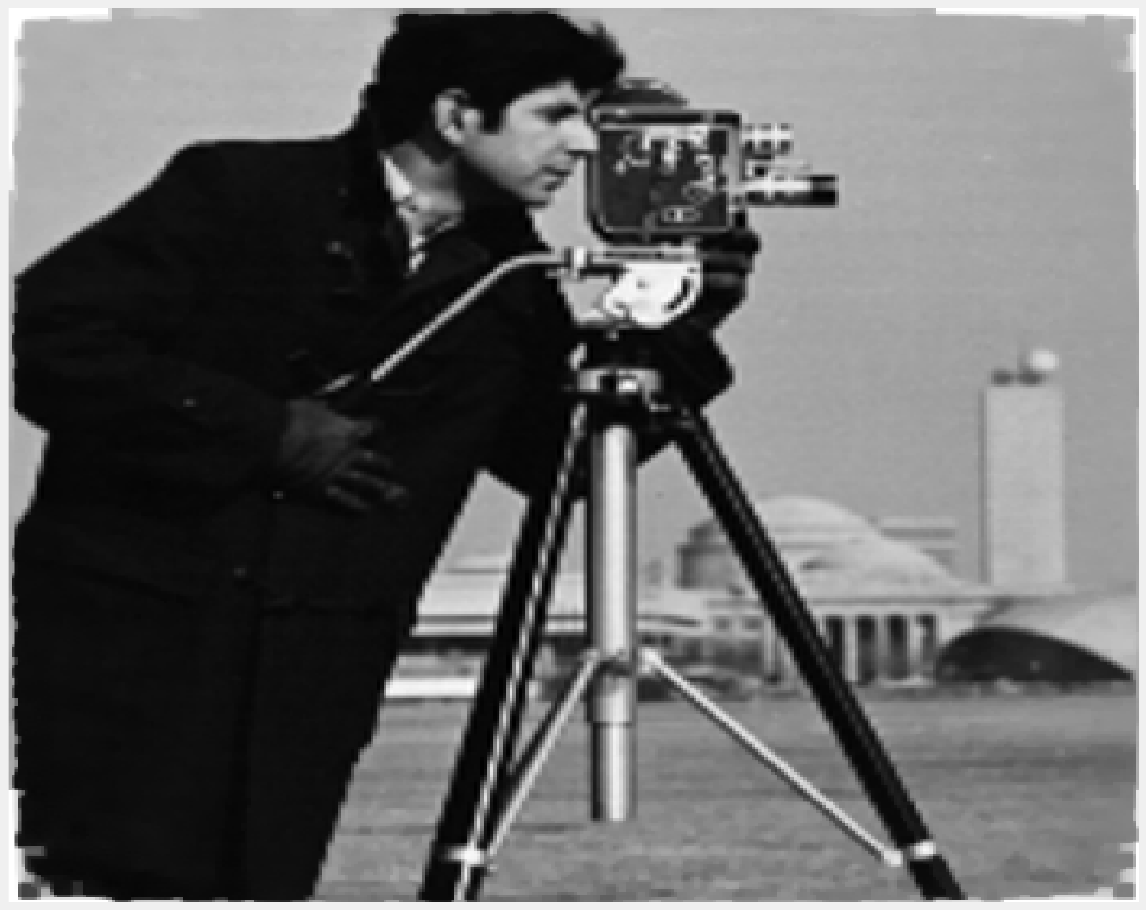}}  & \subcaptionbox{AITV (full)}{\includegraphics[width = 1.35in]{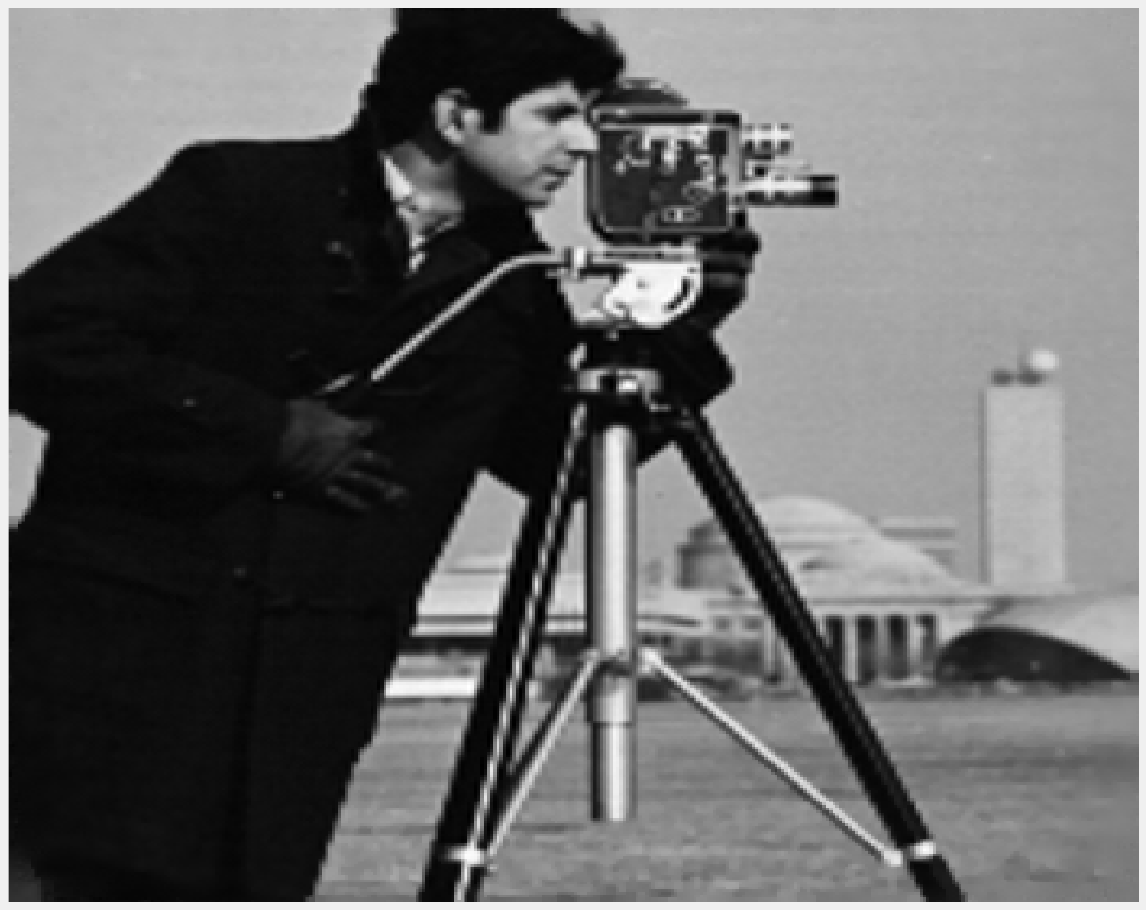}} & \subcaptionbox{DR}{\includegraphics[width = 1.35in]{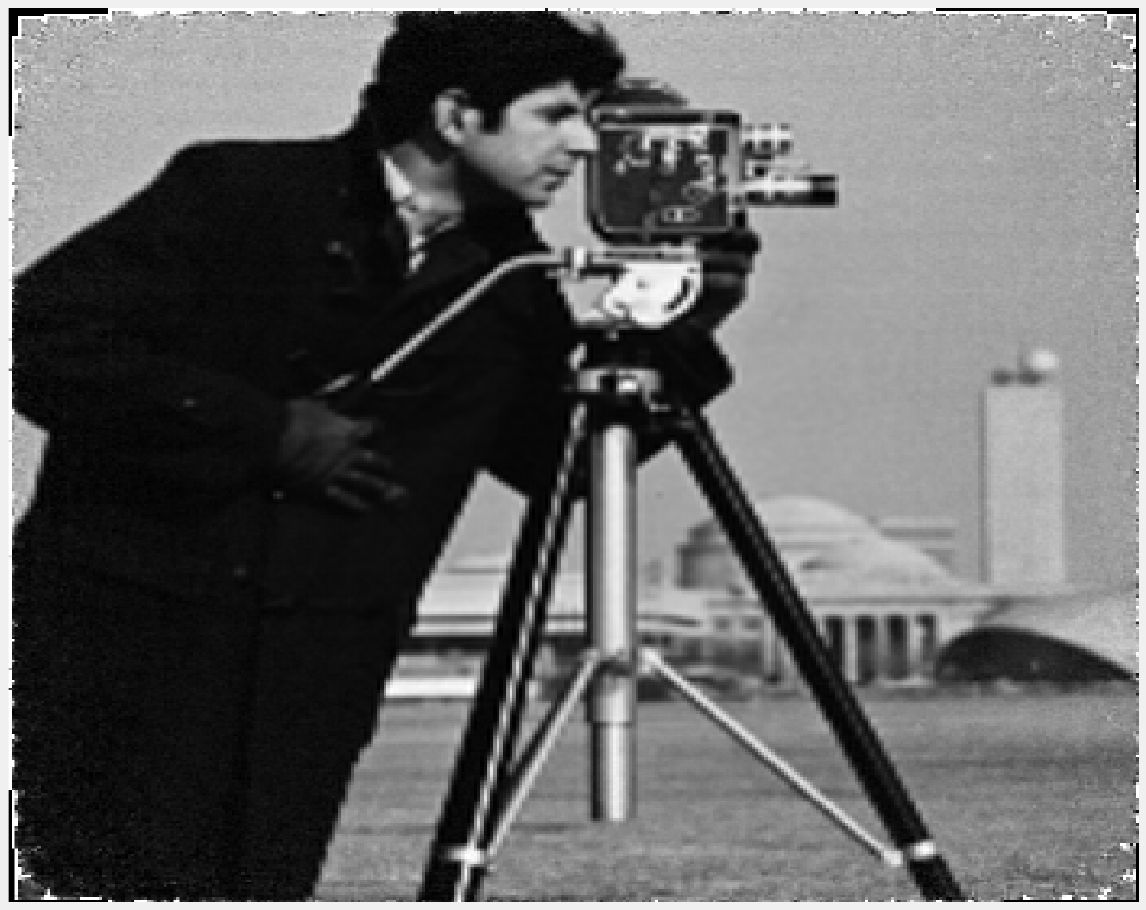}}  & \subcaptionbox{rPIE}{\includegraphics[width = 1.35in]{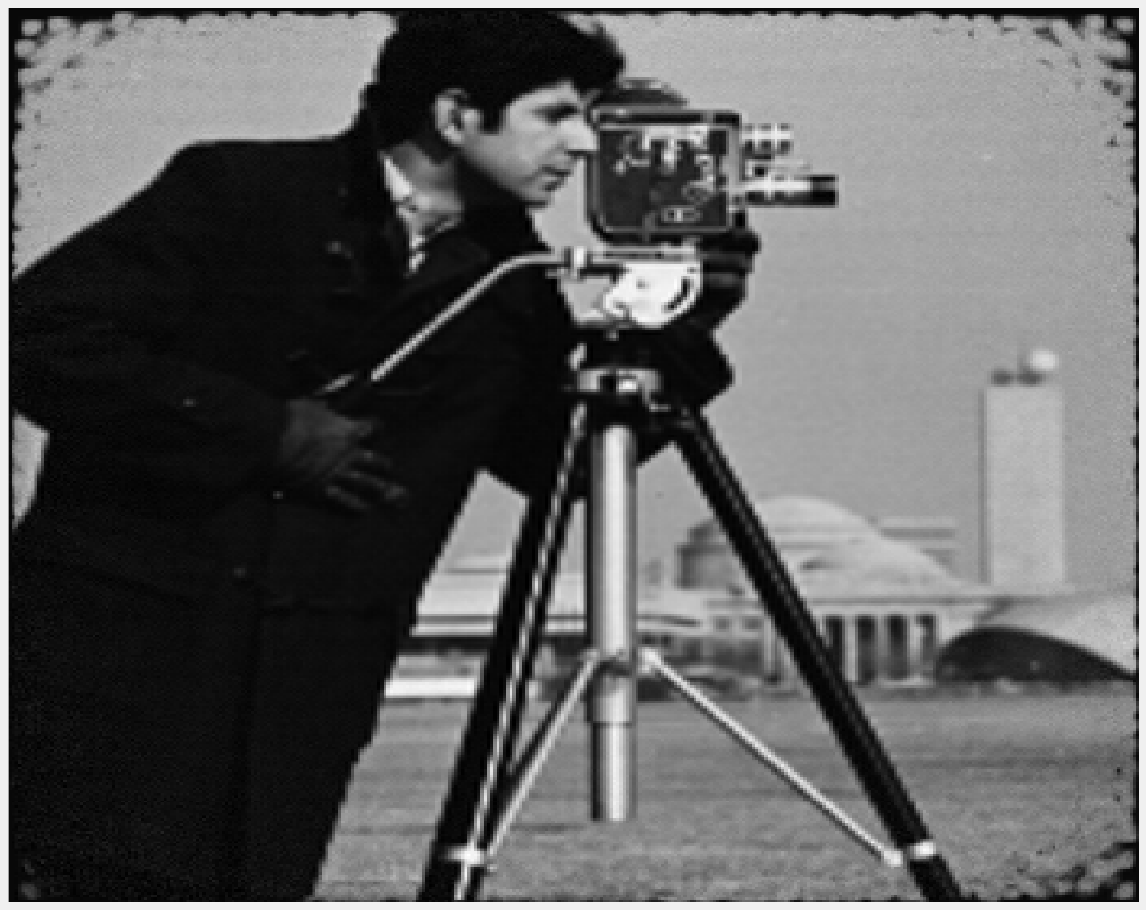}}  & \subcaptionbox{PHeBIE}{\includegraphics[width = 1.35in]{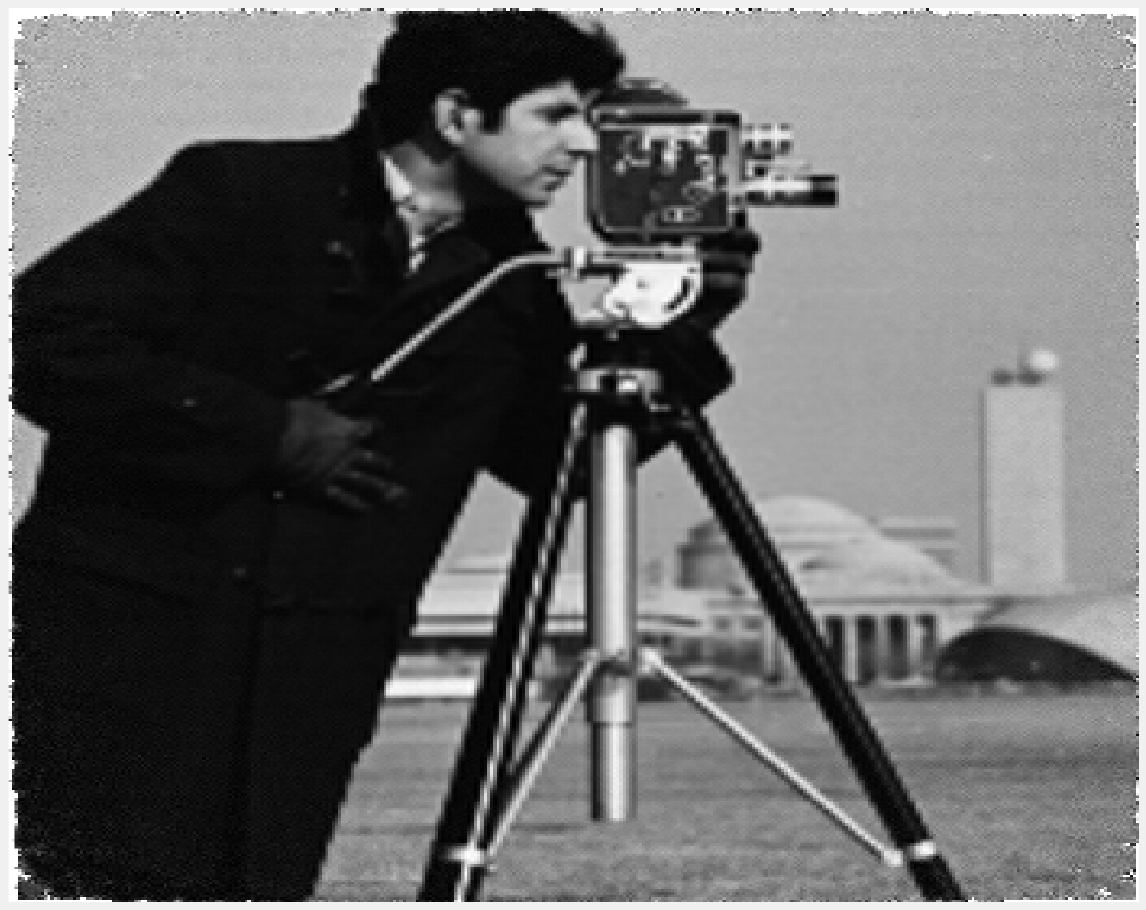}}
				\\
				\subcaptionbox{isoTV ($b=5$)}{\includegraphics[width = 1.35in]{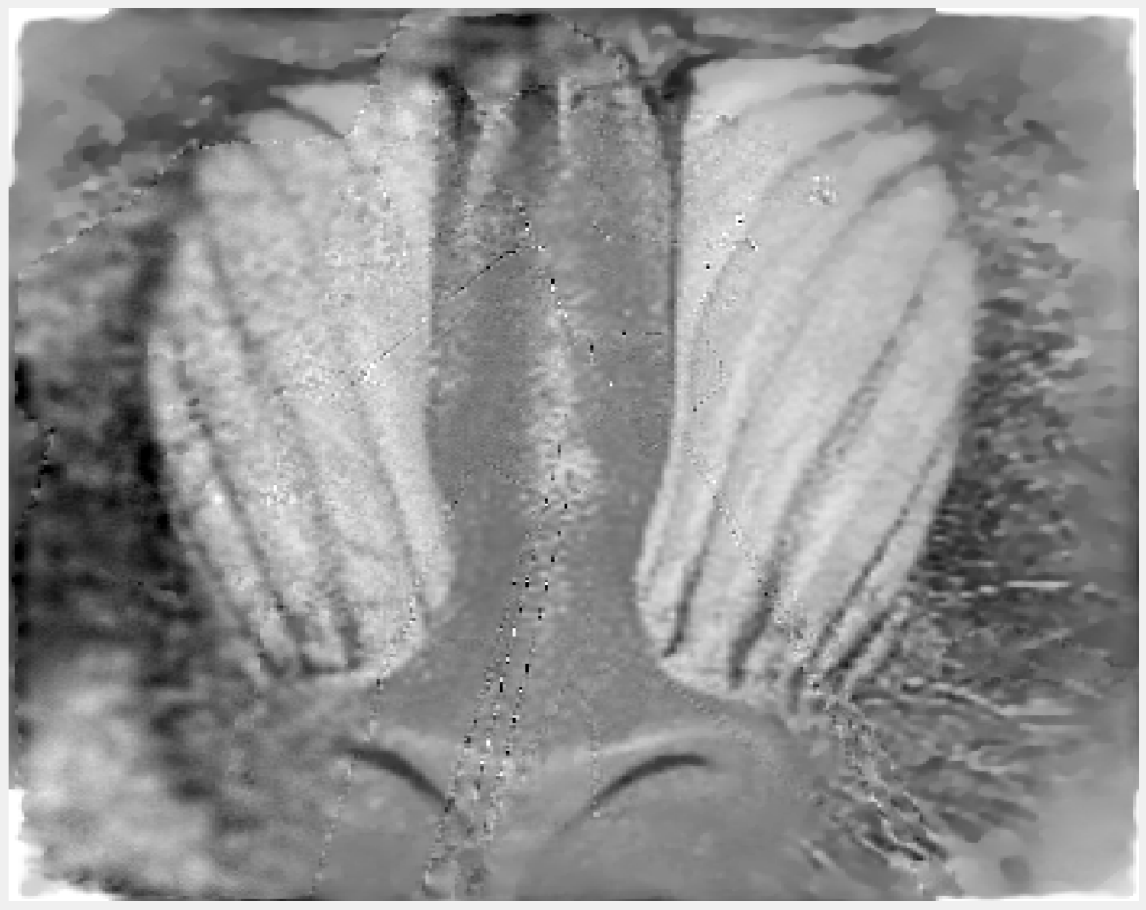}} & \subcaptionbox{AITV ($b=10$)}{\includegraphics[width = 1.35in]{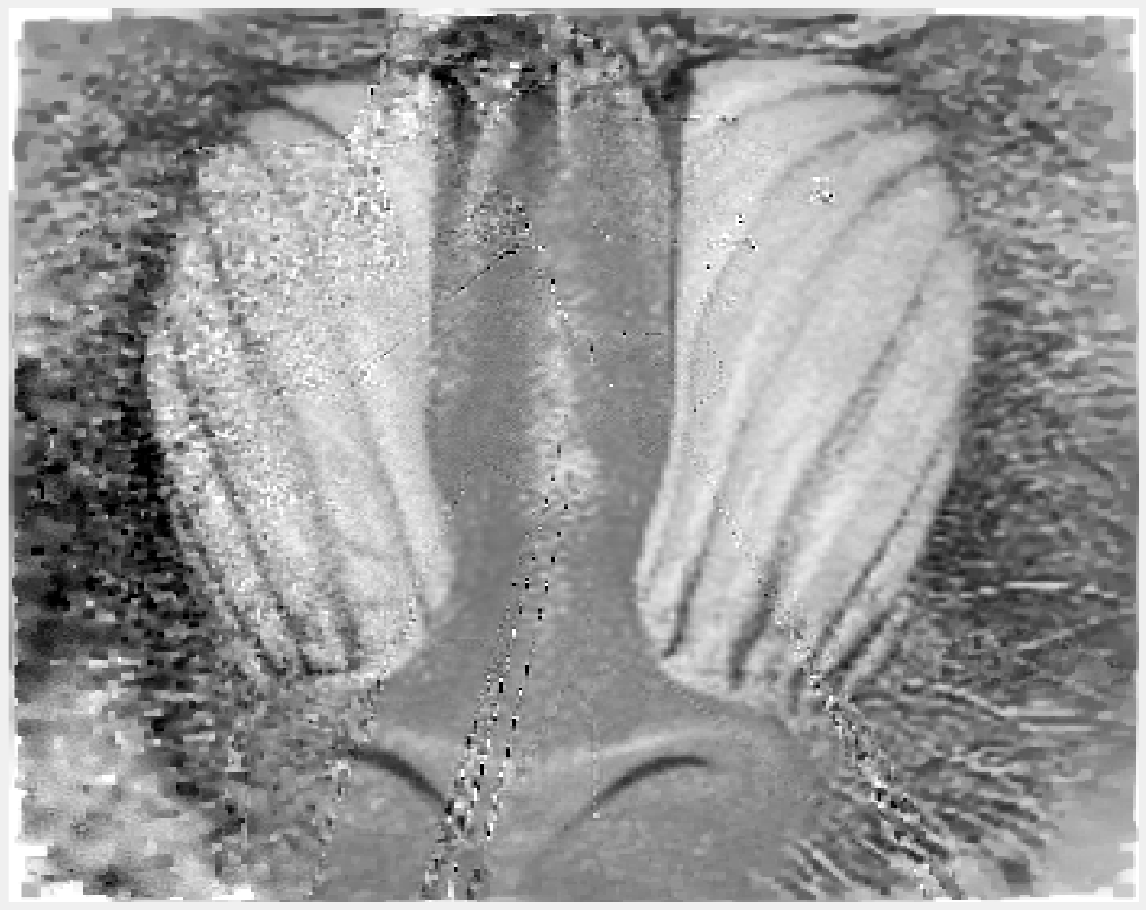}}   & \subcaptionbox{AITV (full)}{\includegraphics[width = 1.35in]{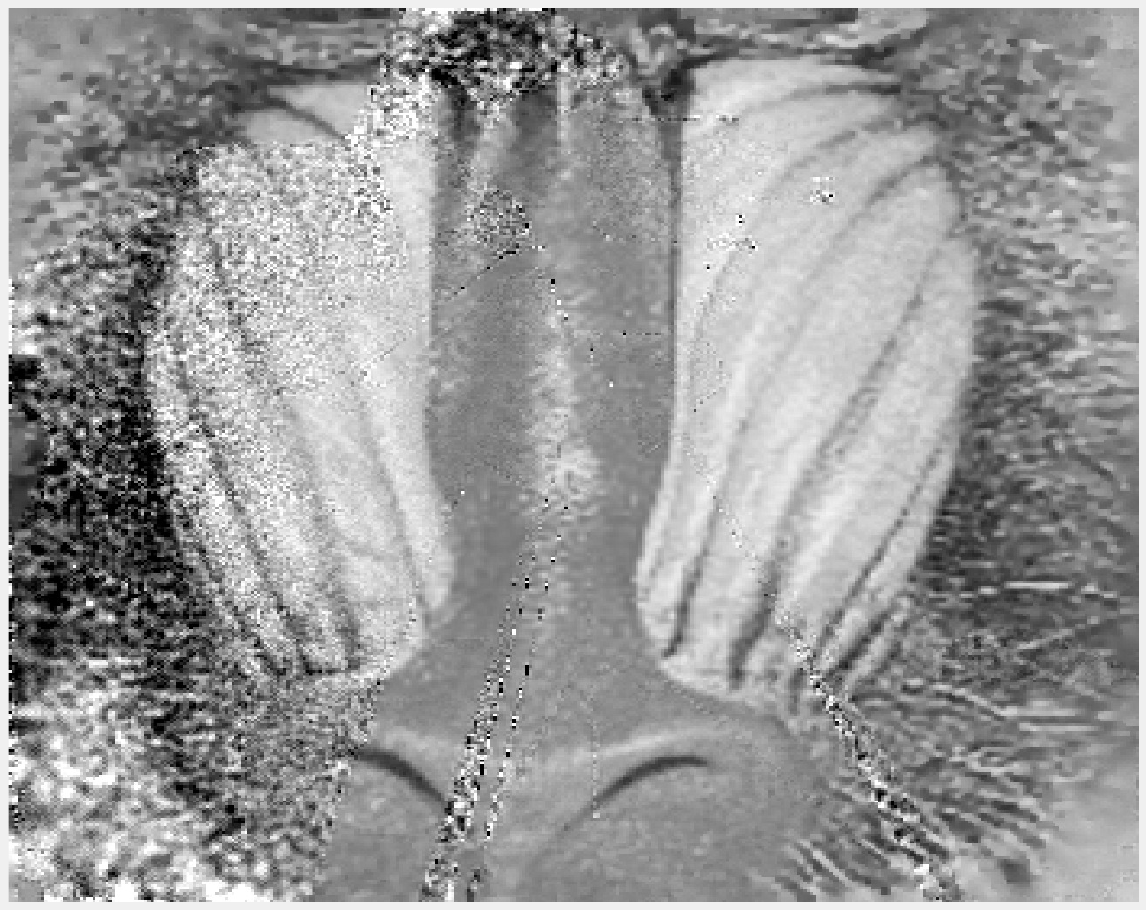}}  & \subcaptionbox{DR}{\includegraphics[width = 1.35in]{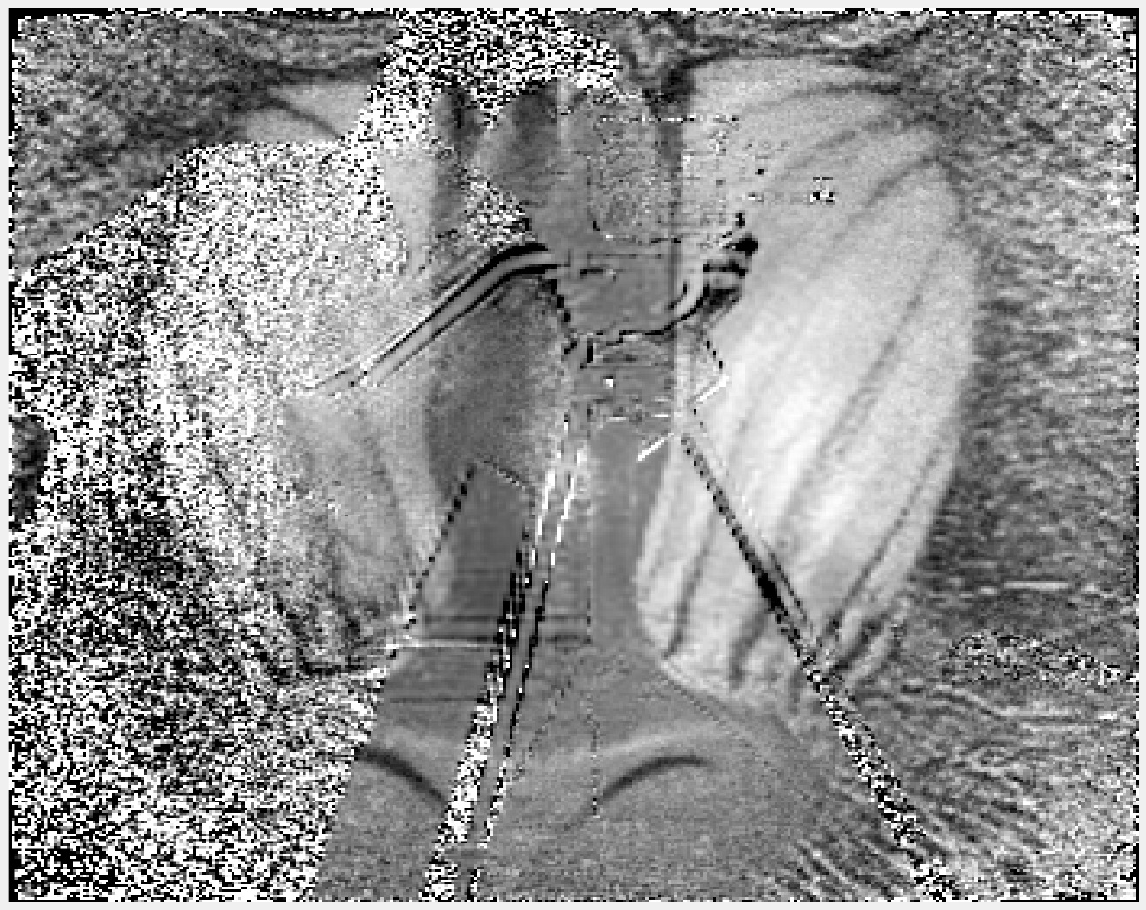}} & \subcaptionbox{rPIE}{\includegraphics[width = 1.35in]{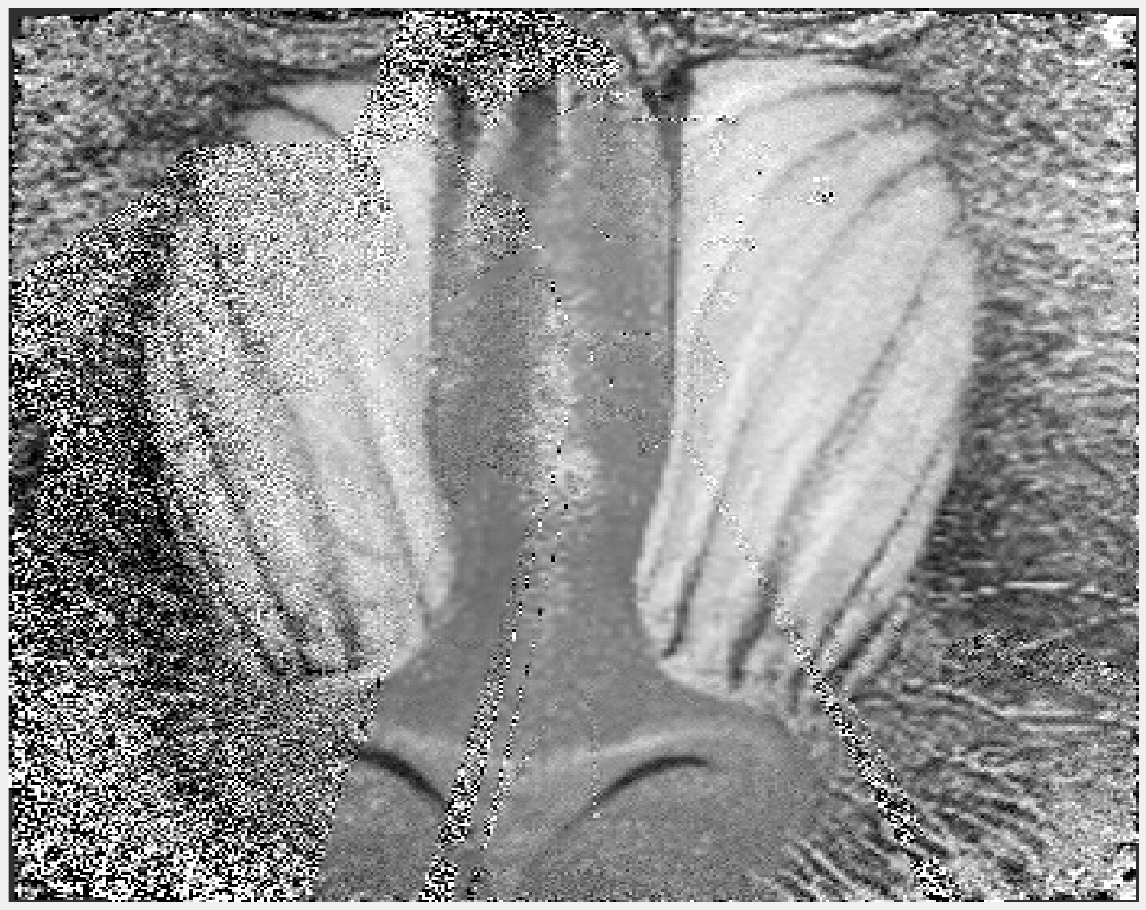}}& \subcaptionbox{PHeBIE}{\includegraphics[width = 1.35in]{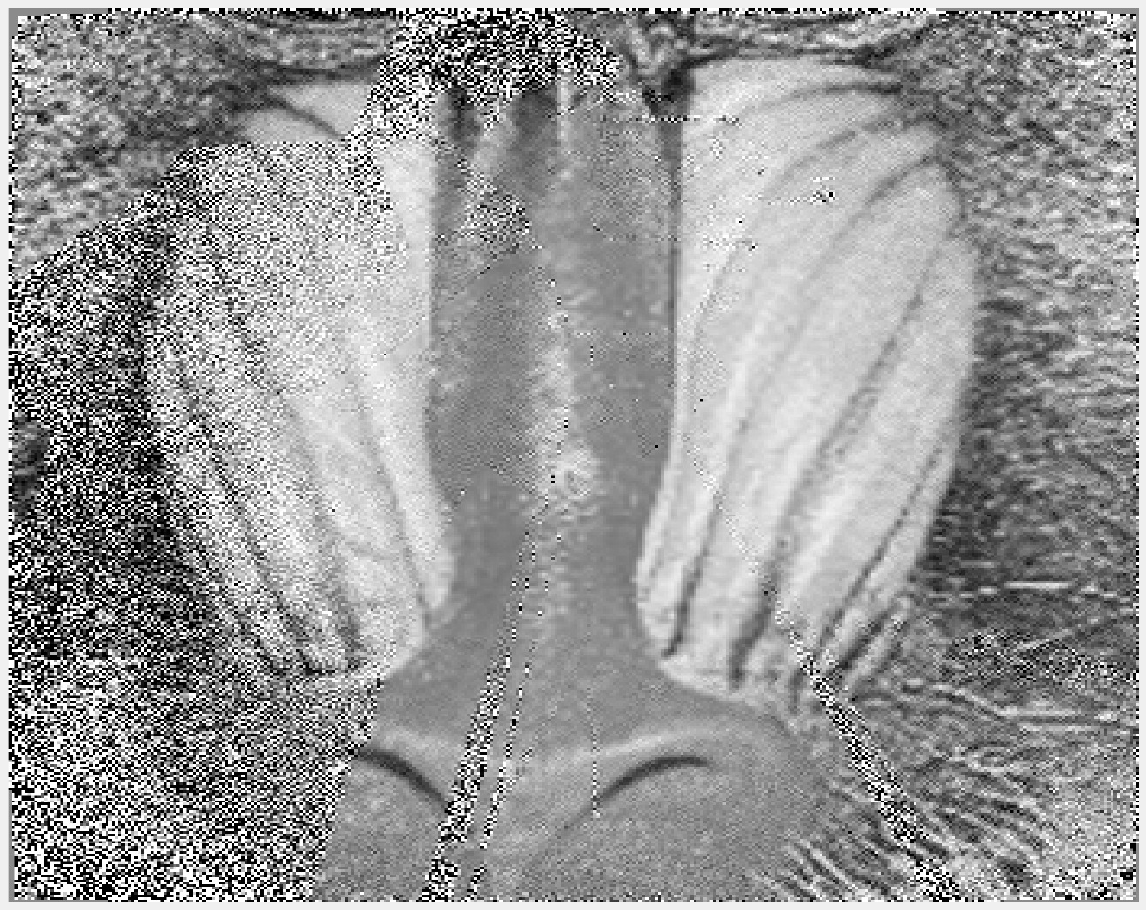}}
		\end{tabular}}
		\caption{Reconstructions of the non-blind case for the Gaussian noise with SNR = 40. Top two rows: reconstructions of Figures \ref{fig:chip_mag}-\ref{fig:chip_phase}; bottom two rows: reconstructions of Figs. \ref{fig:cameraman_mag}-\ref{fig:cameraman_phase}.}
		\label{fig:agm_nonblind}
\end{minipage}
	
\end{figure}\clearpage}%

\afterpage{
\begin{figure}[t]
	\begin{minipage}{\linewidth}
		\centering
		\resizebox{\textwidth}{!}{%
			\begin{tabular}{c@{}c@{}c@{}c@{}c@{}c}
				\subcaptionbox{isoTV ($b=5$)}{\includegraphics[width = 1.35in]{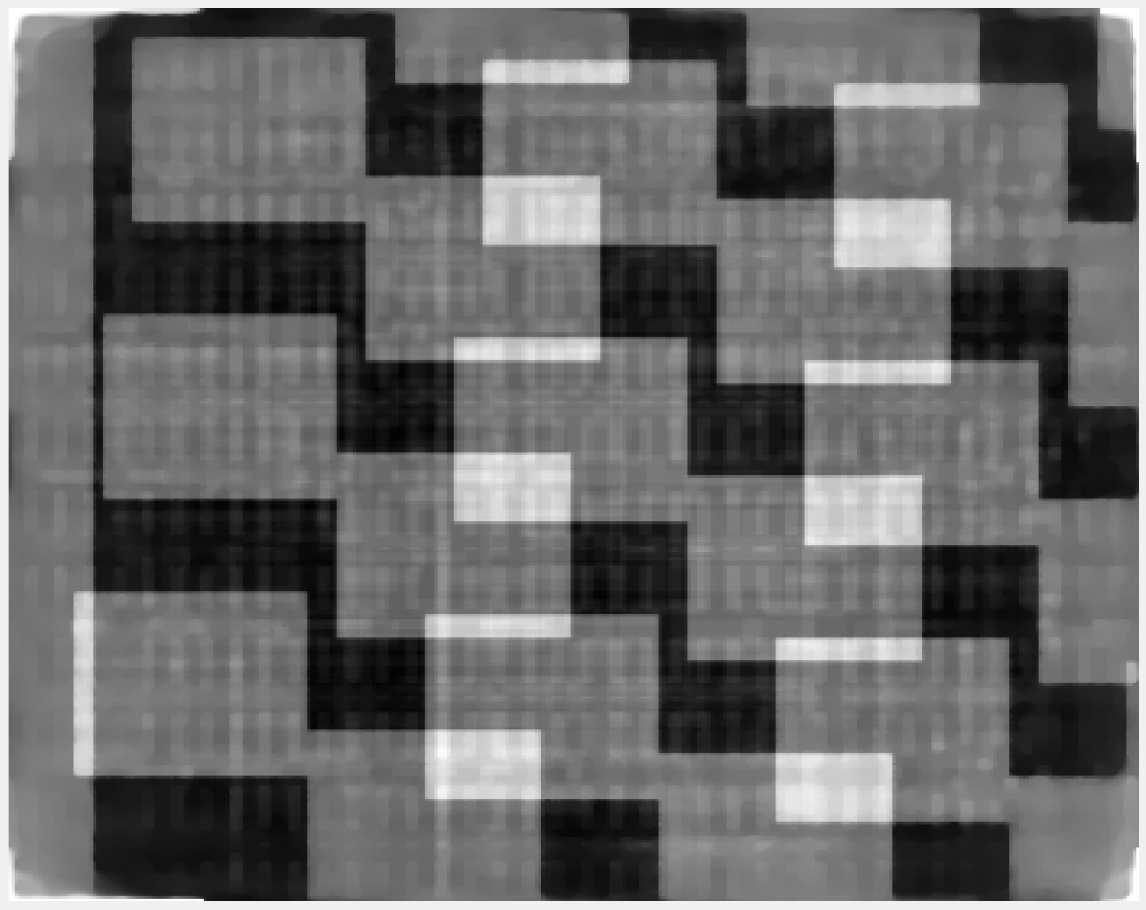}} & \subcaptionbox{AITV ($b=20$)}{\includegraphics[width = 1.35in]{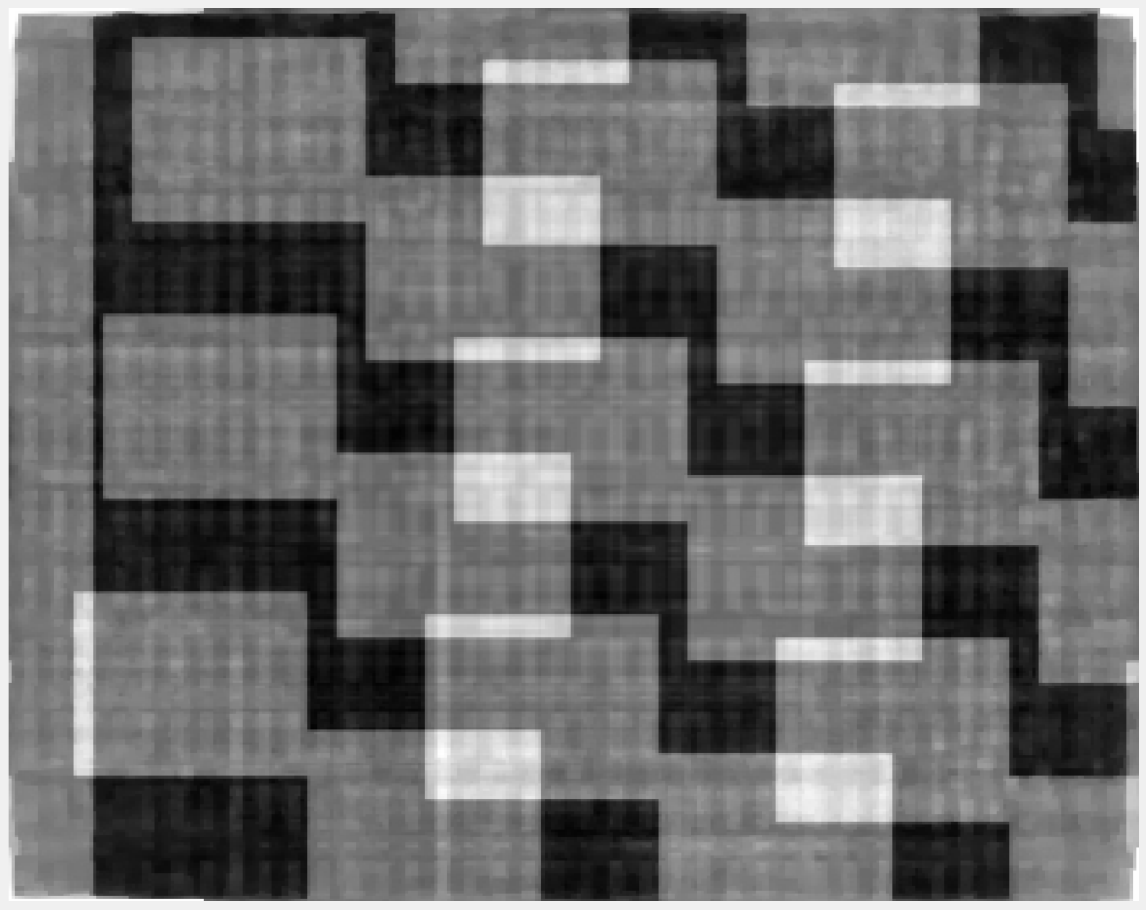}}  & \subcaptionbox{AITV (full)}{\includegraphics[width = 1.35in]{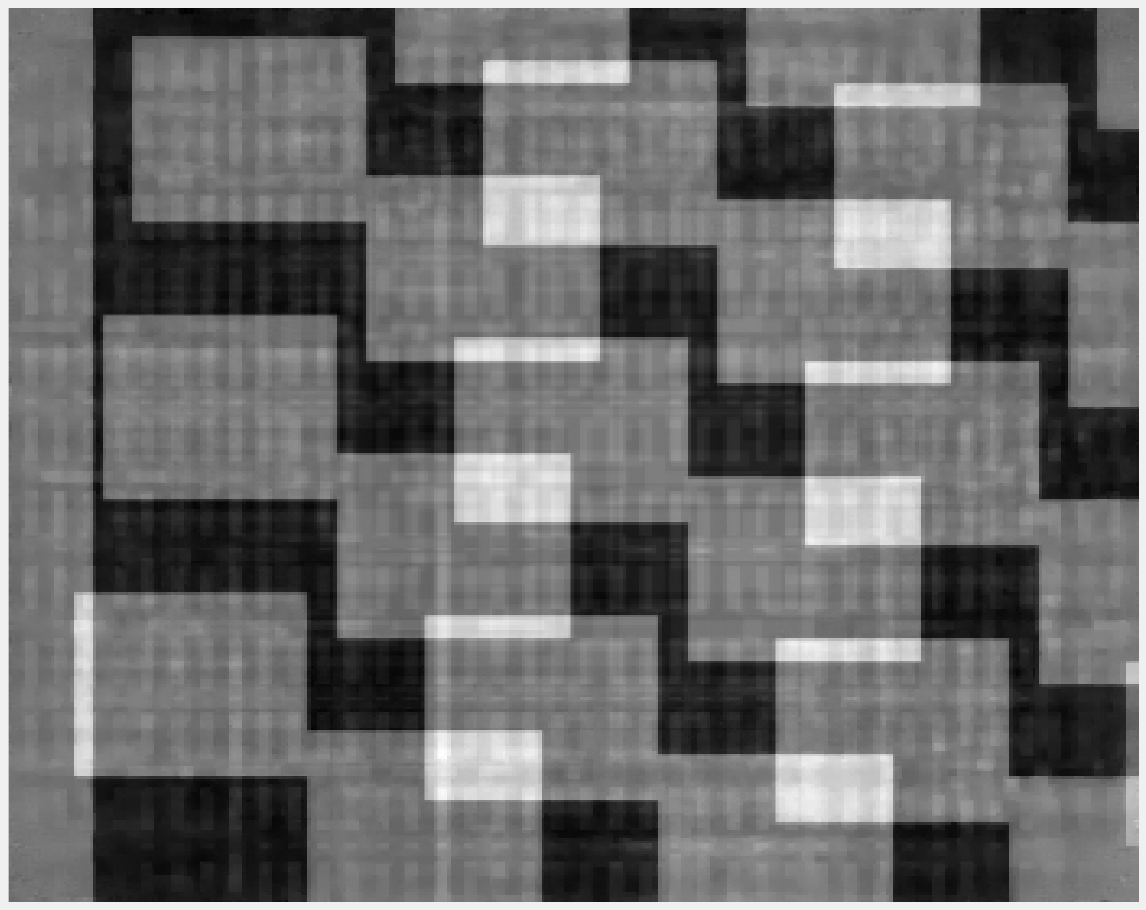}}  & \subcaptionbox{DR}{\includegraphics[width = 1.35in]{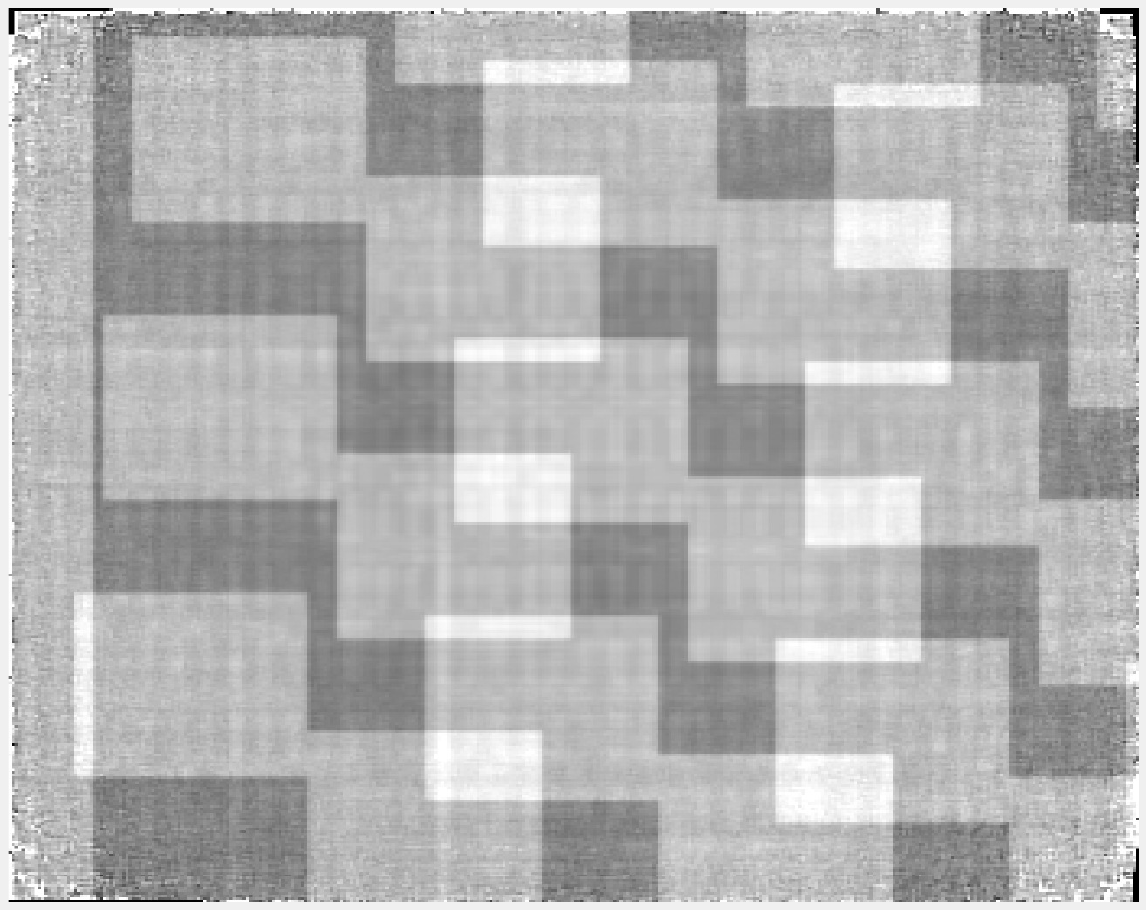}} & \subcaptionbox{rPIE}{\includegraphics[width = 1.35in]{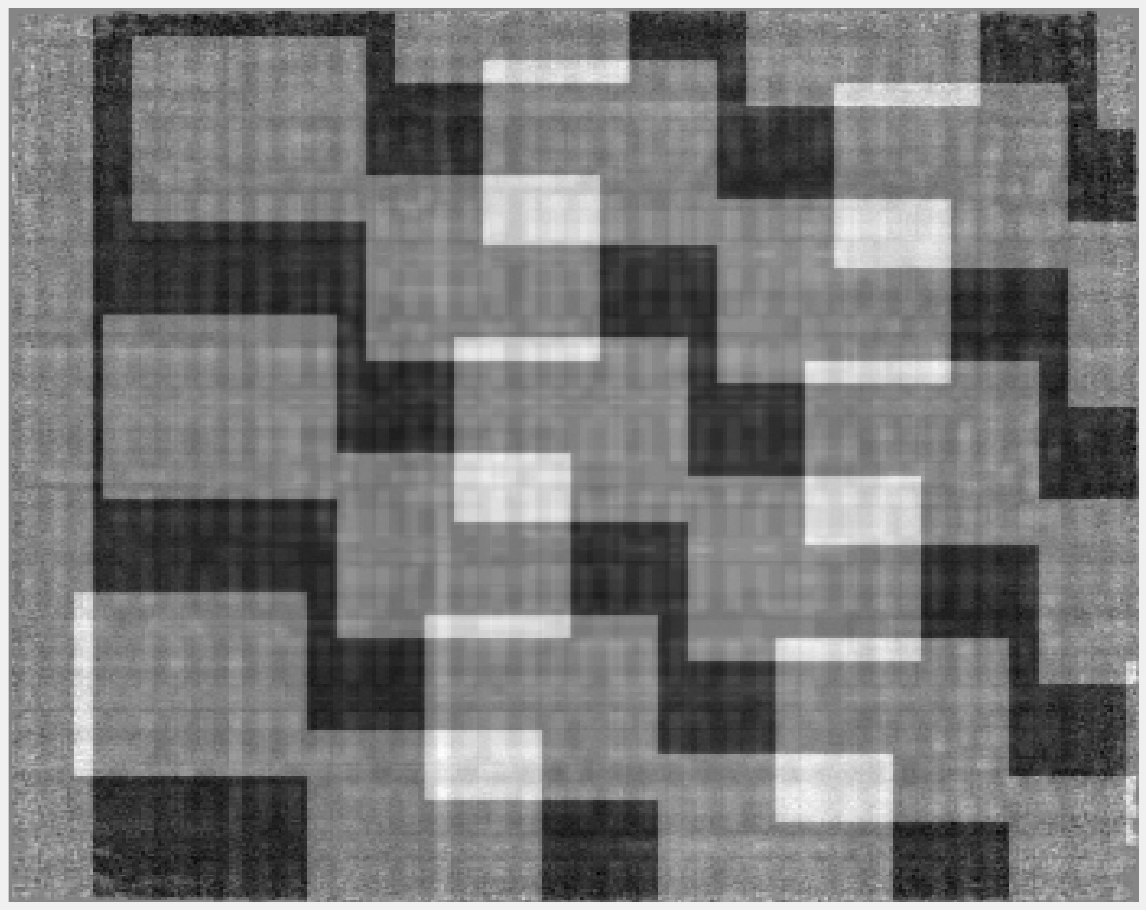}}  & \subcaptionbox{PHeBIE}{\includegraphics[width = 1.35in]{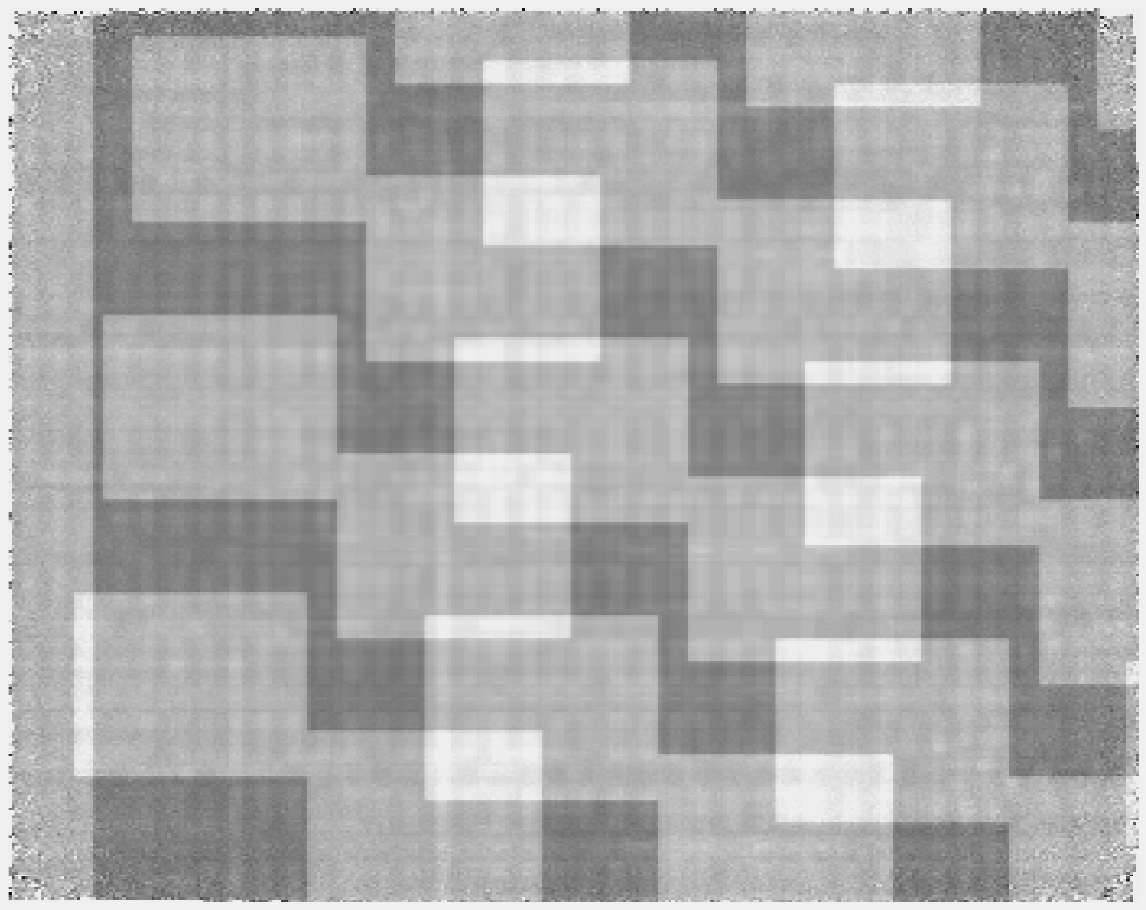}}\\	\subcaptionbox{isoTV $(b=5)$}{\includegraphics[width = 1.35in]{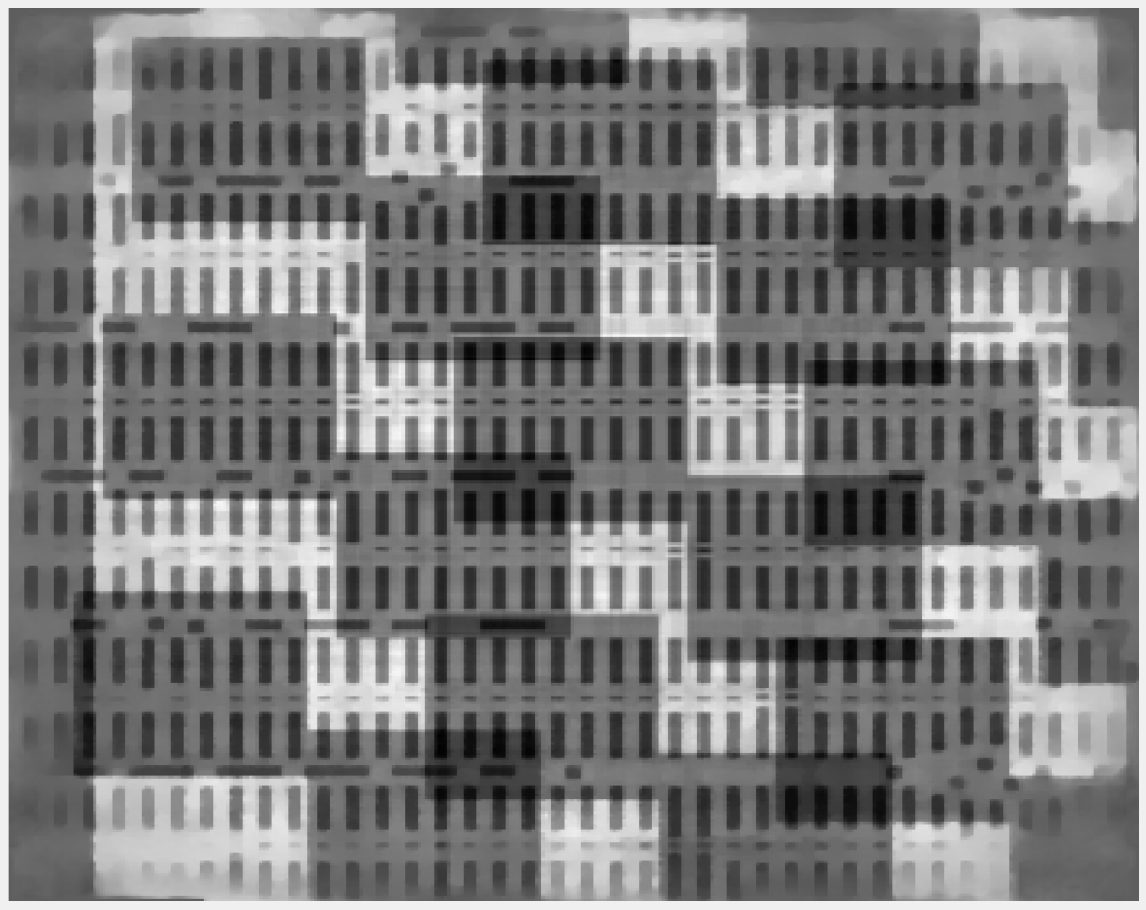}} & \subcaptionbox{AITV $(b=20)$}{\includegraphics[width = 1.35in]{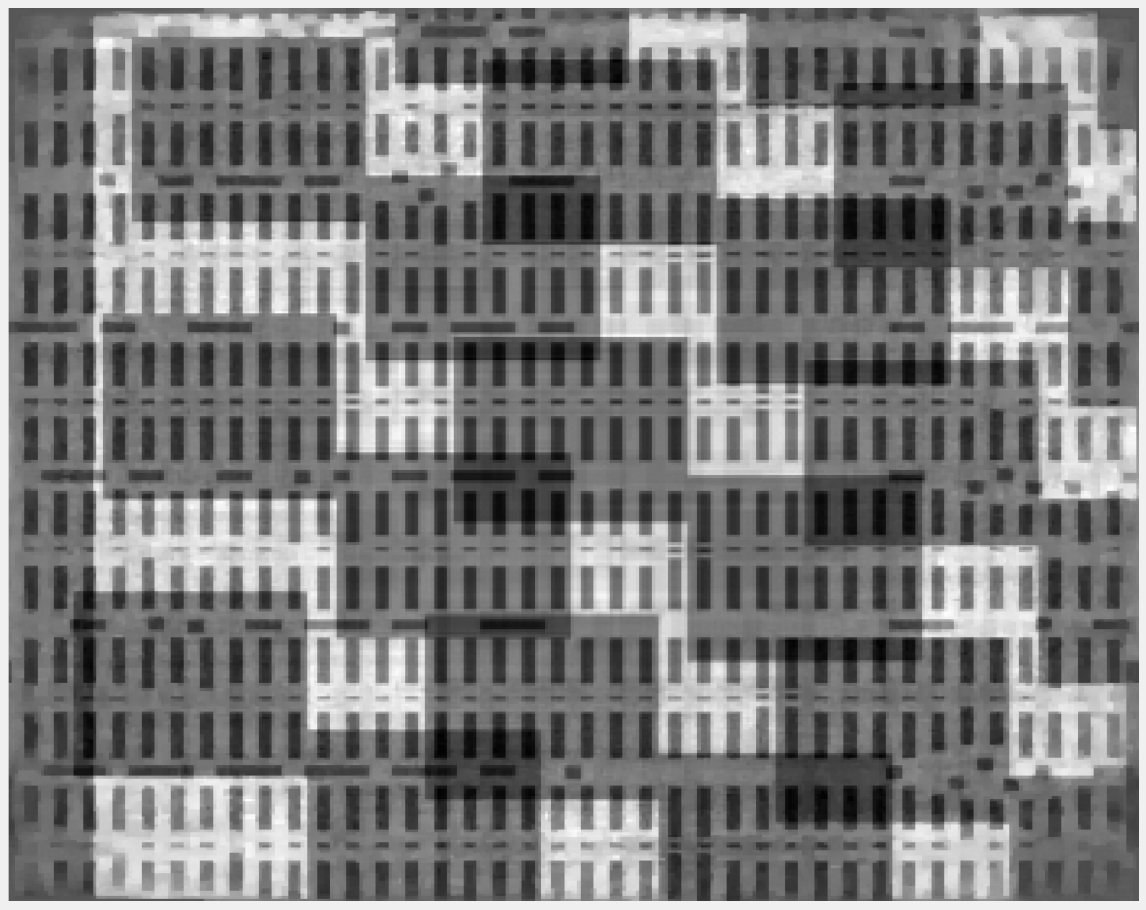}}   & \subcaptionbox{AITV (full)}{\includegraphics[width = 1.35in]{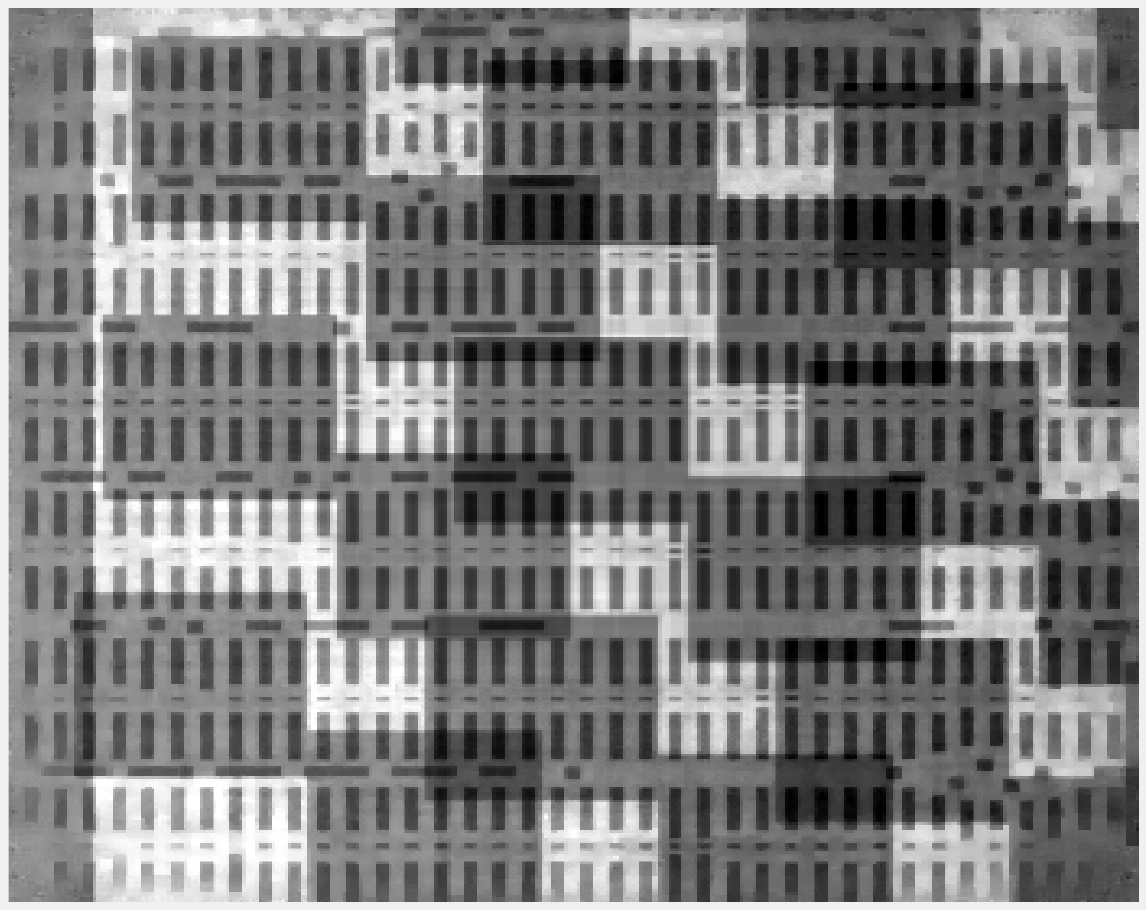}} &  \subcaptionbox{DR}{\includegraphics[width = 1.35in]{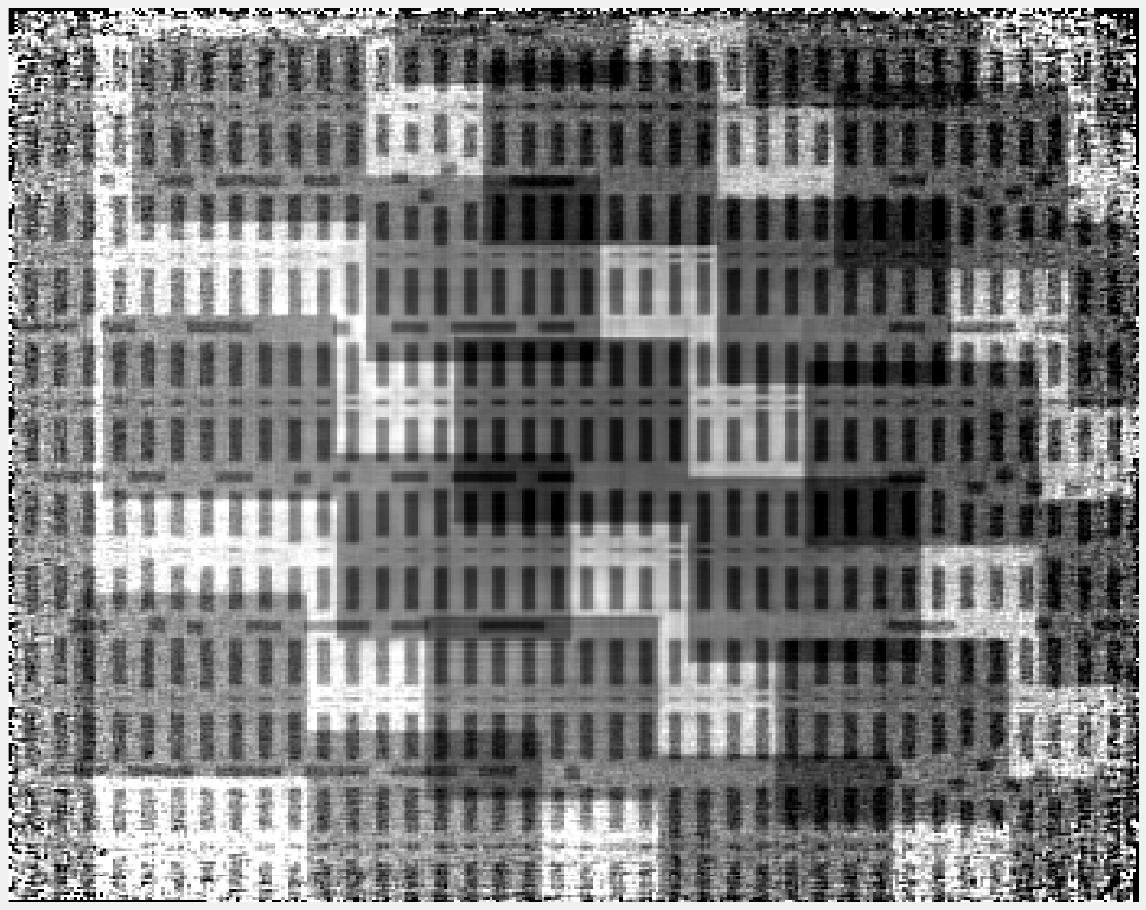}} & \subcaptionbox{rPIE}{\includegraphics[width = 1.35in]{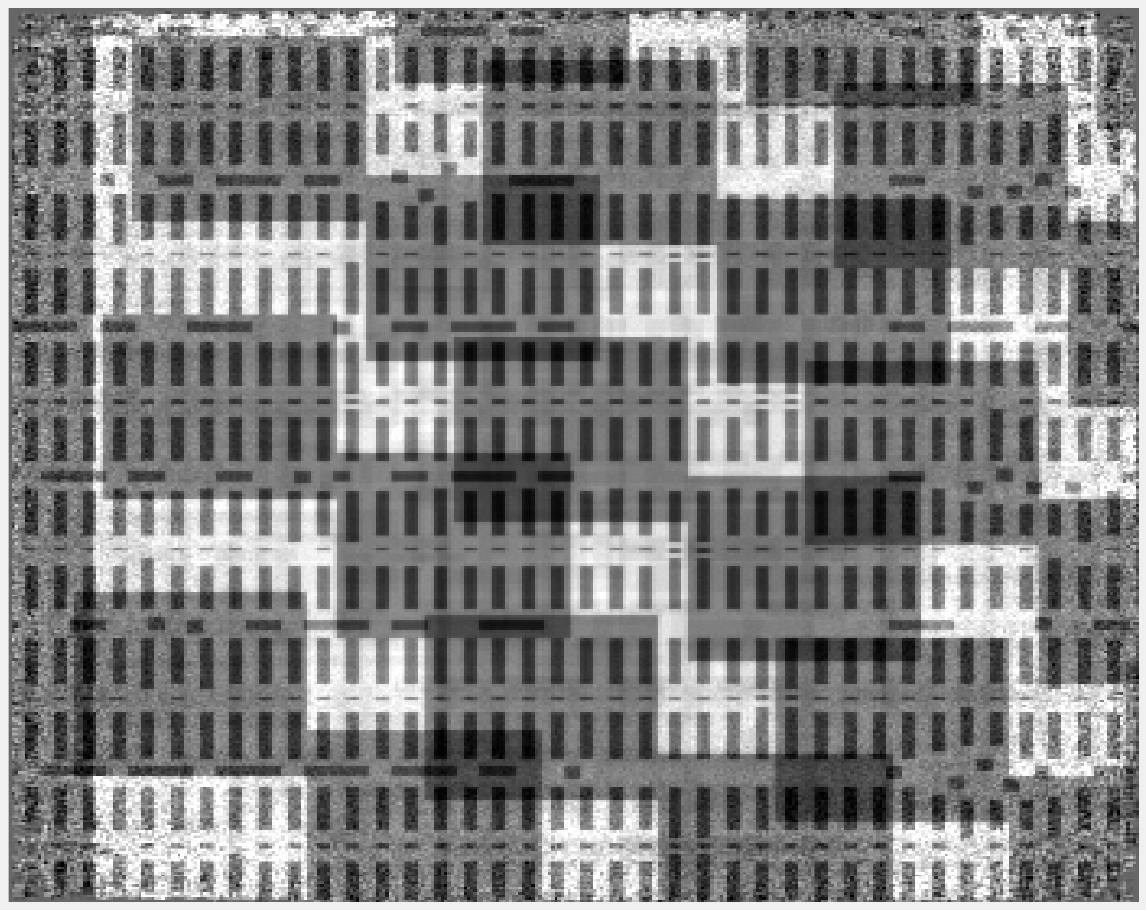}}& \subcaptionbox{PHeBIE}{\includegraphics[width = 1.35in]{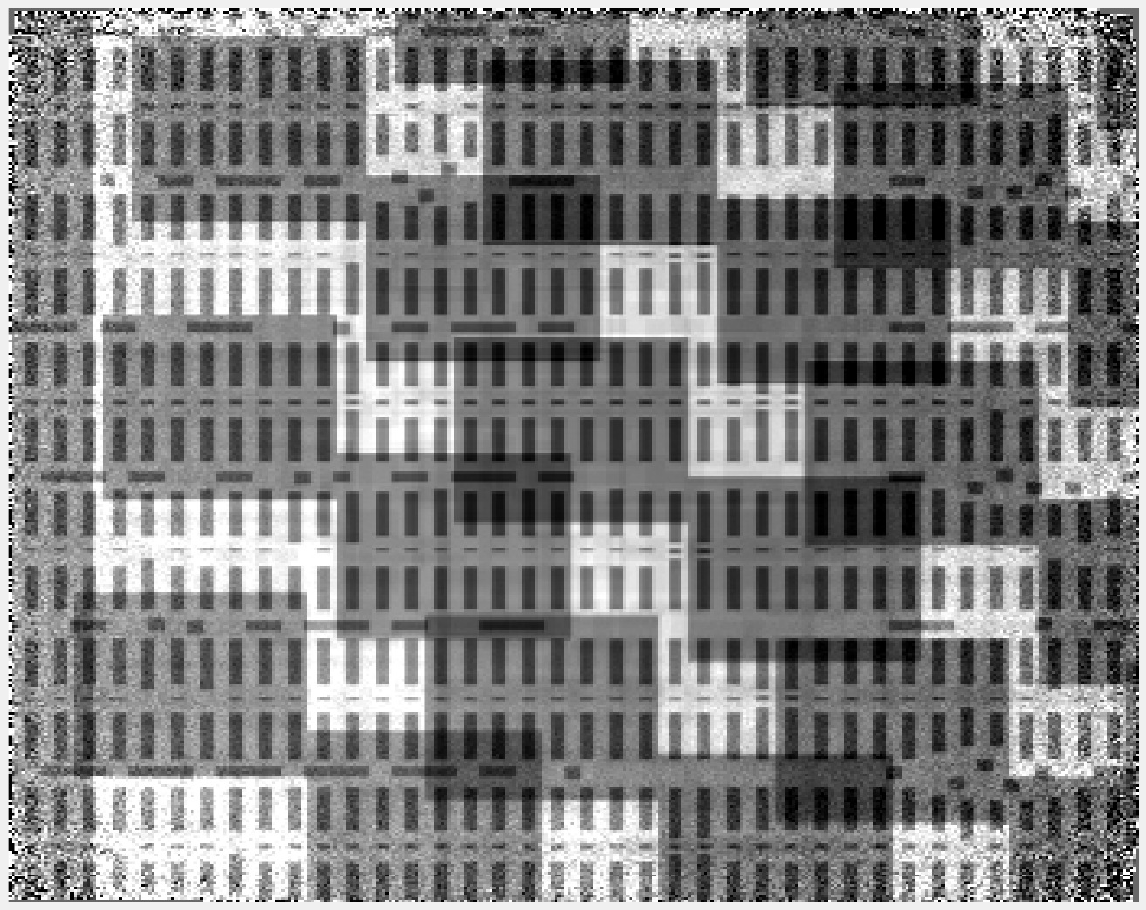}}
				\\
	\subcaptionbox{isoTV ($b=10$)}{\includegraphics[width = 1.35in]{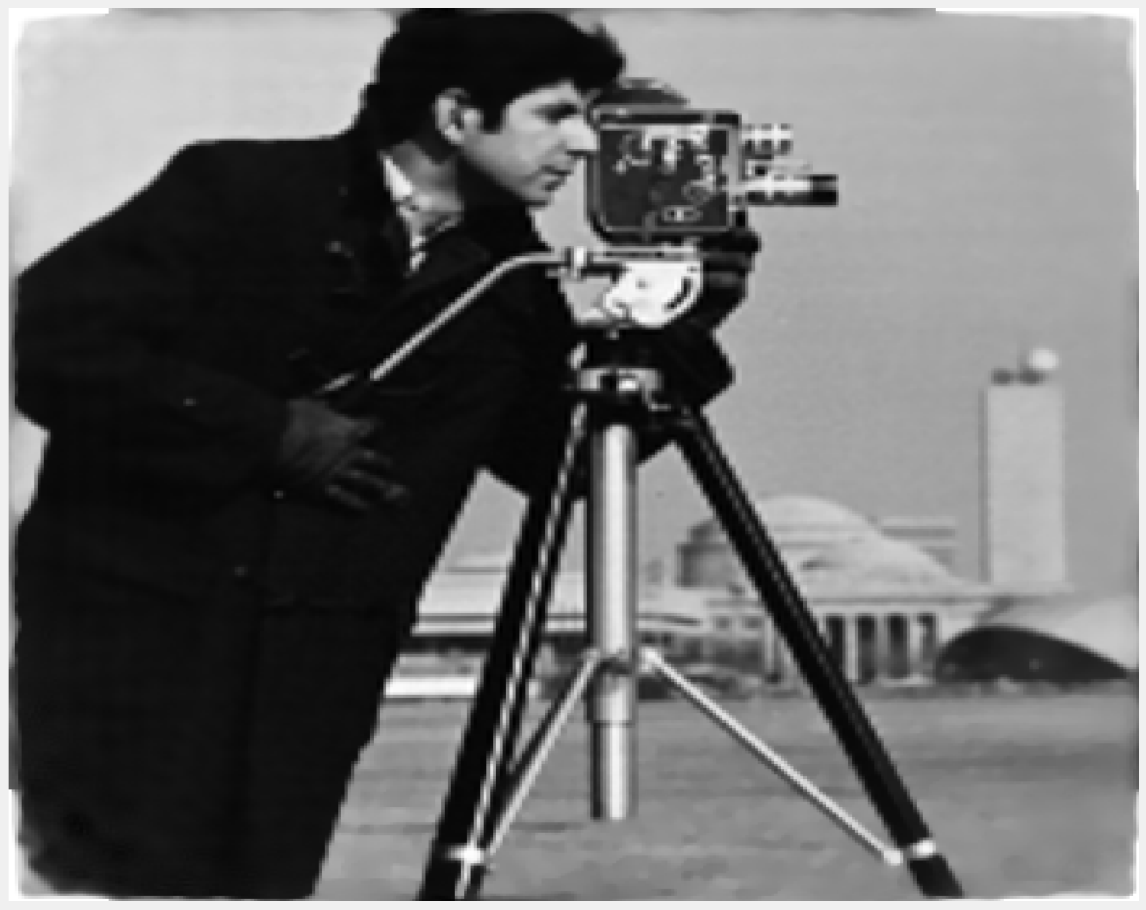}} & \subcaptionbox{AITV ($b=10$)}{\includegraphics[width = 1.35in]{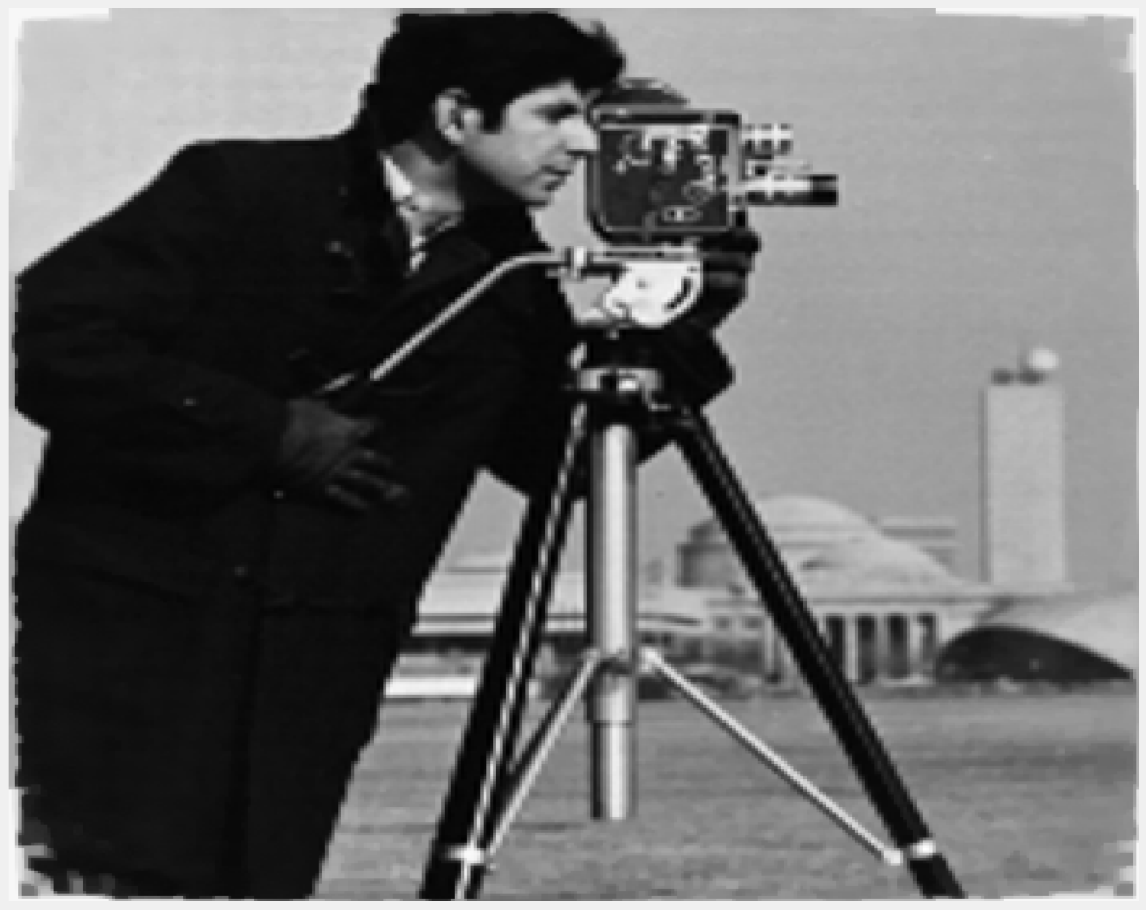}}  & \subcaptionbox{AITV (full)}{\includegraphics[width = 1.35in]{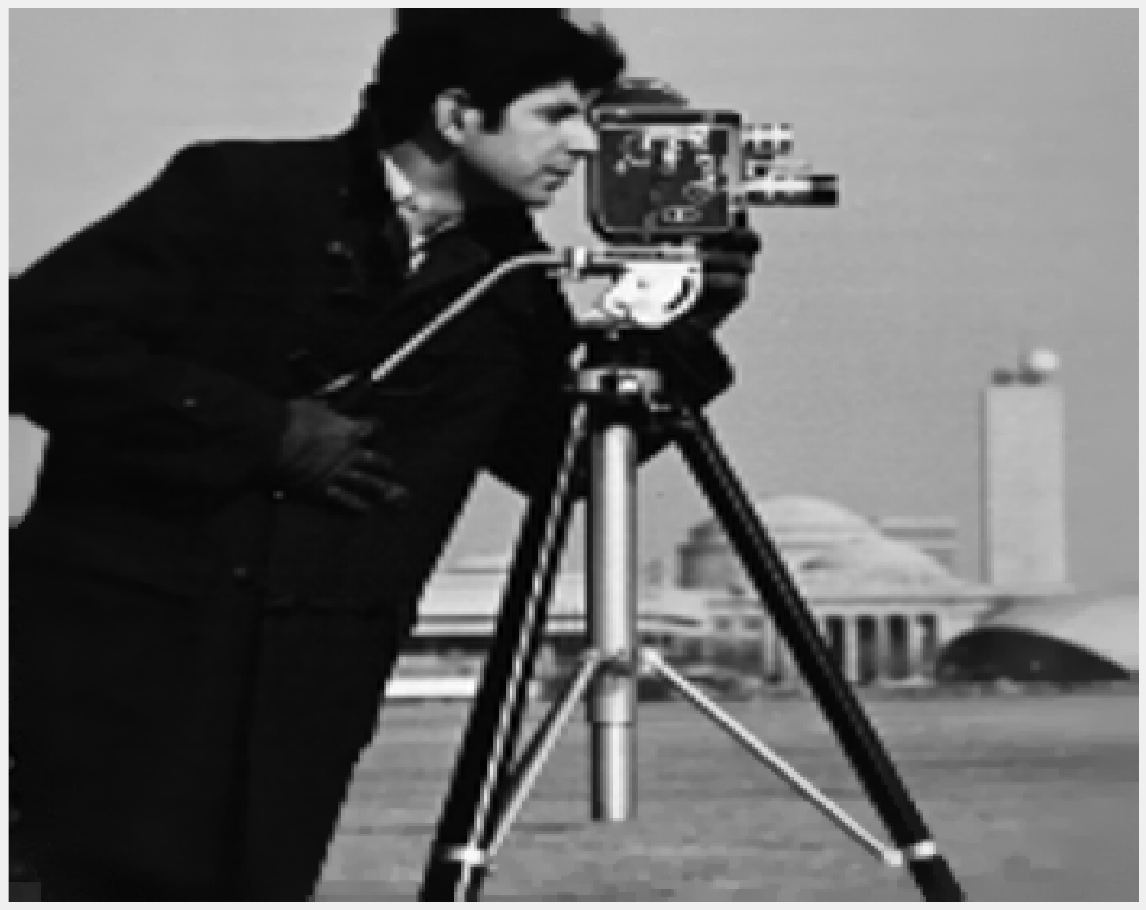}}  & \subcaptionbox{DR}{\includegraphics[width = 1.35in]{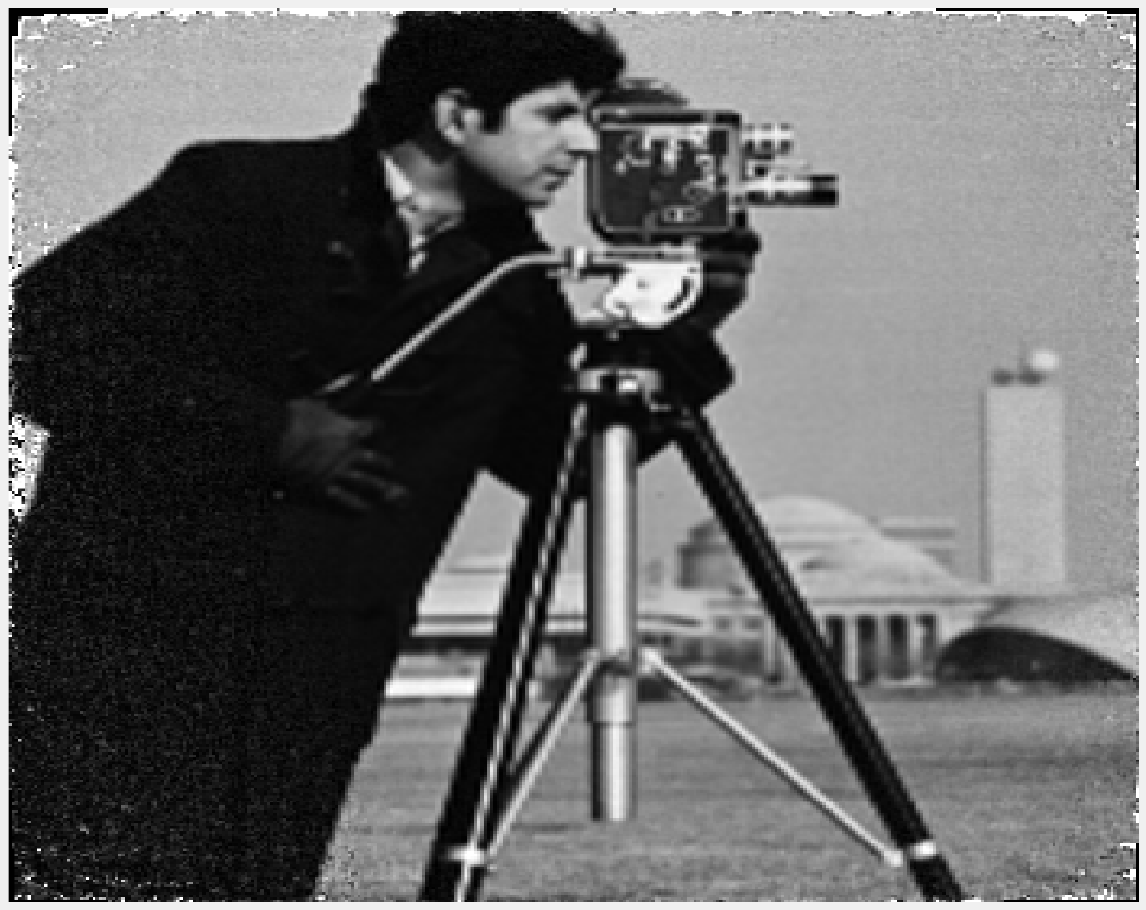}}& \subcaptionbox{rPIE}{\includegraphics[width = 1.35in]{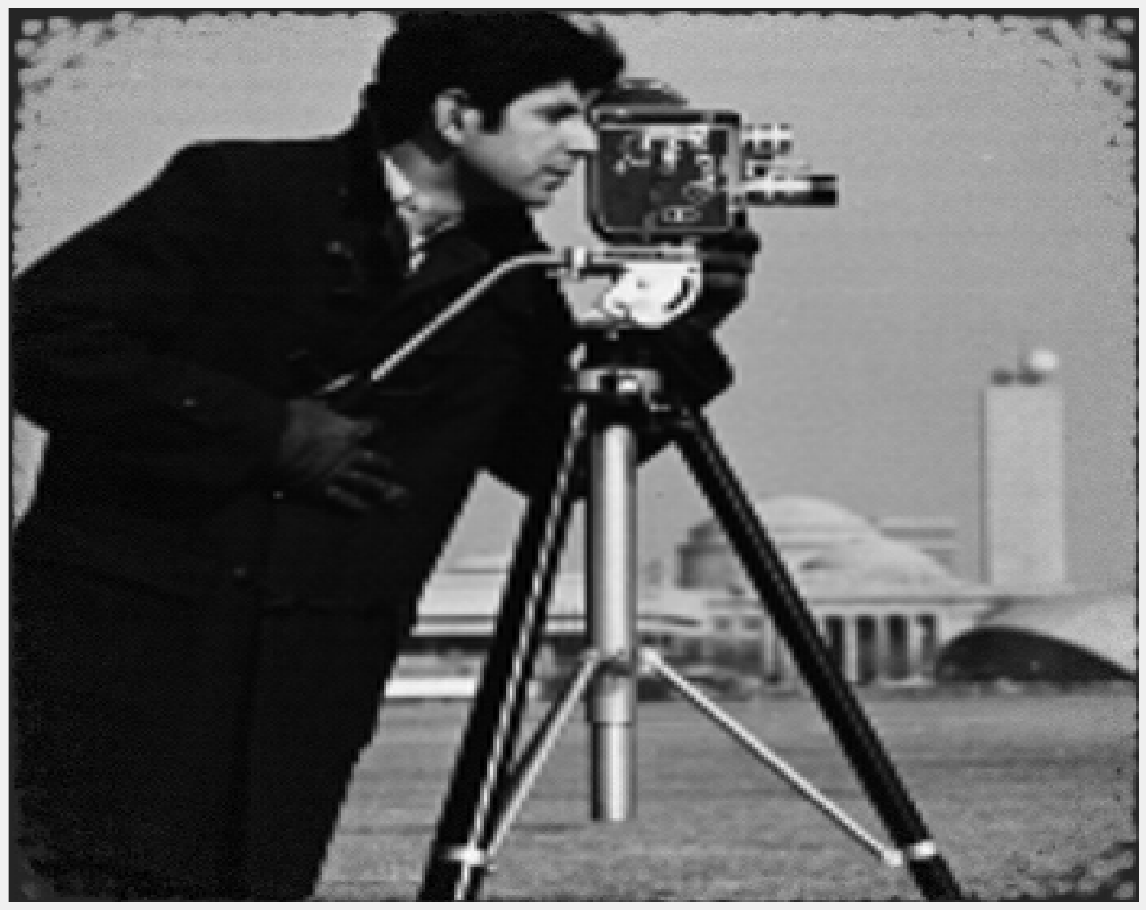}} & \subcaptionbox{PHeBIE}{\includegraphics[width = 1.35in]{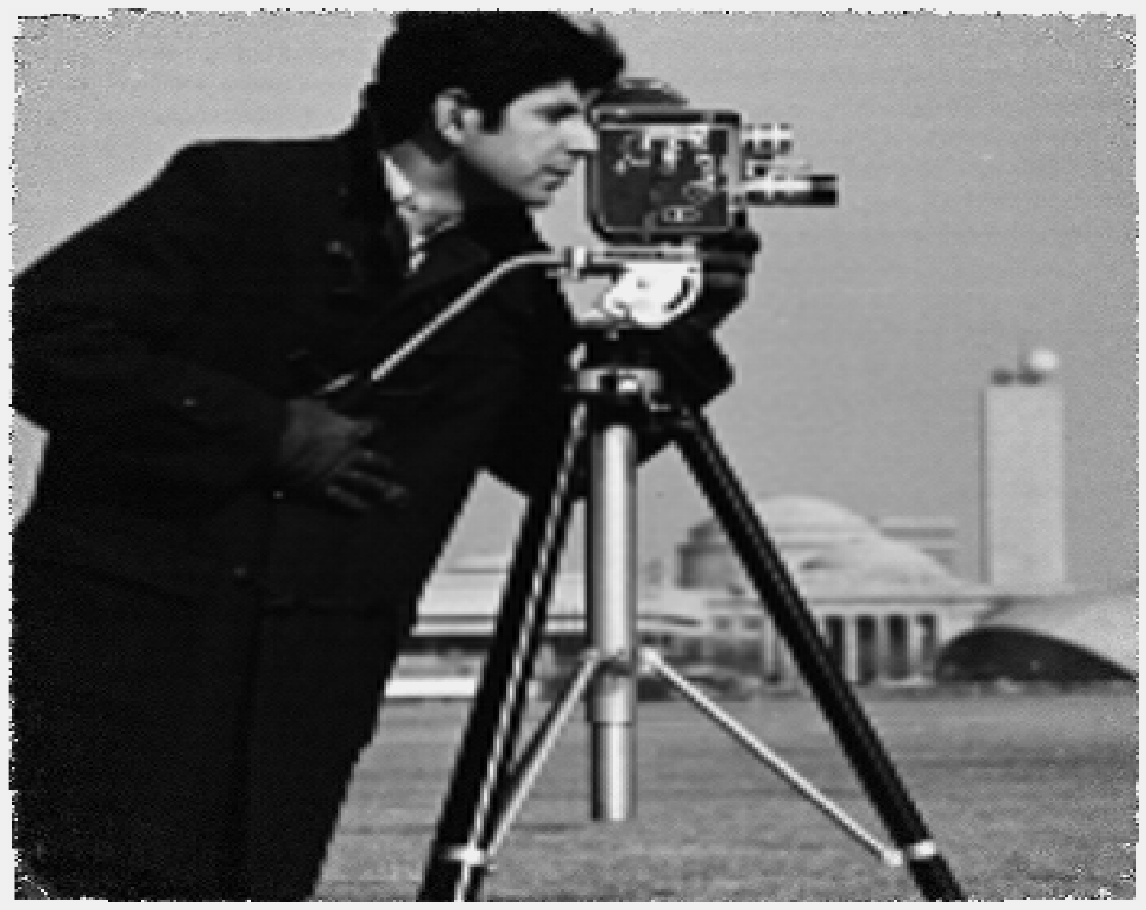}} 
				\\
				\subcaptionbox{isoTV ($b=10$)}{\includegraphics[width = 1.35in]{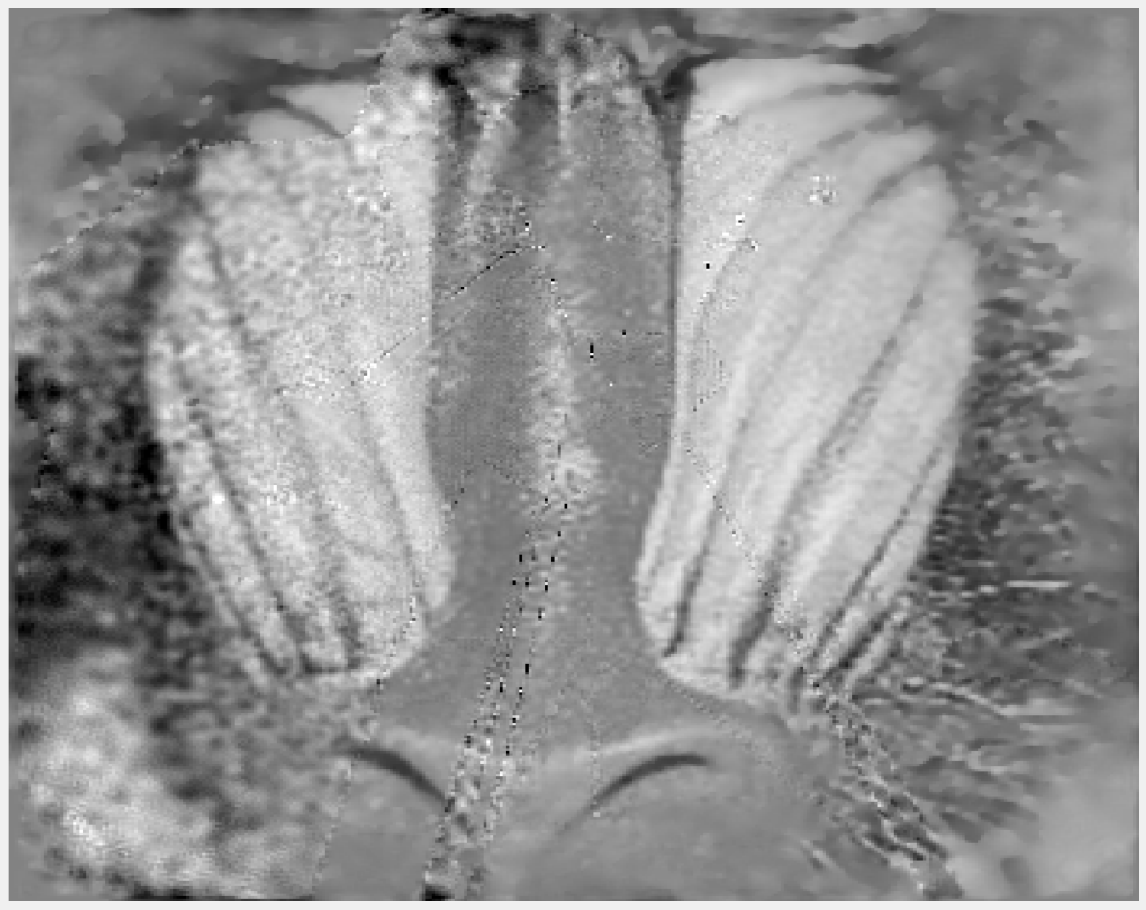}} & \subcaptionbox{AITV ($b=10$)}{\includegraphics[width = 1.35in]{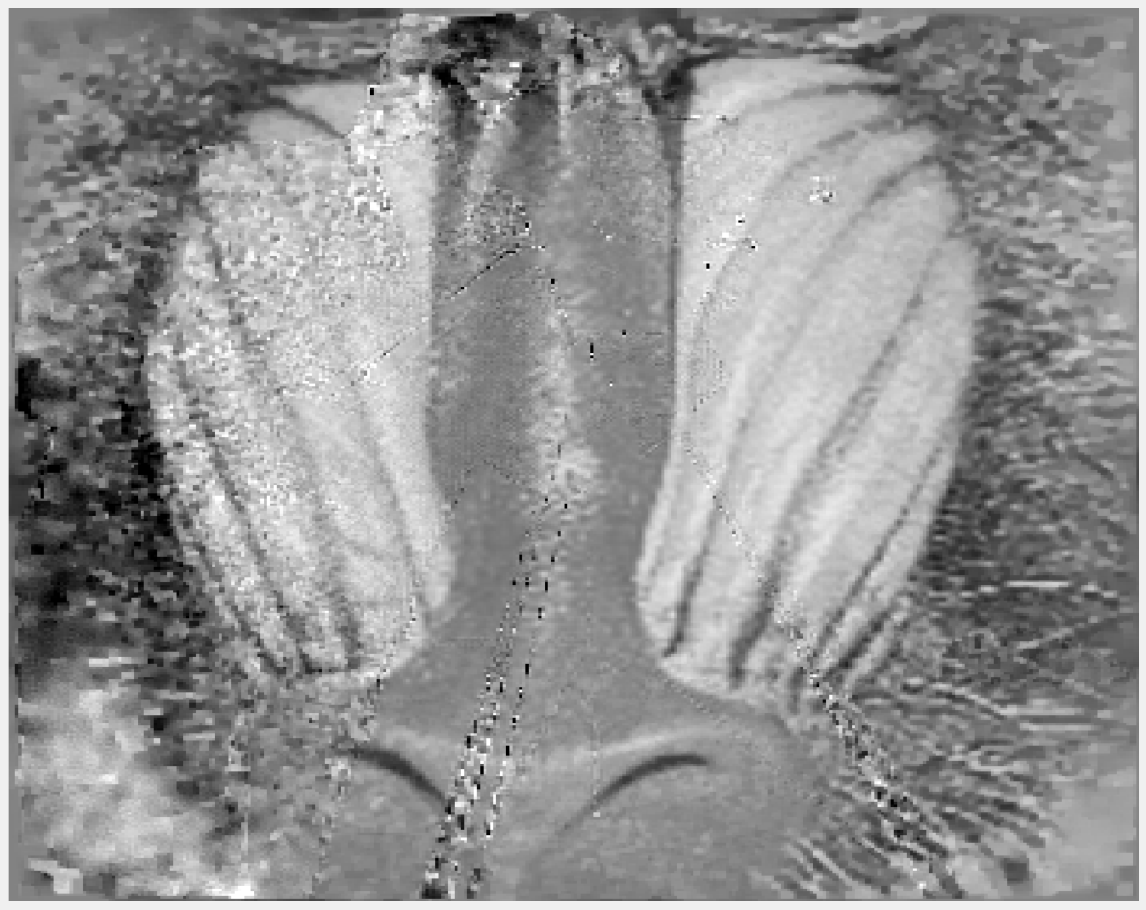}}  & \subcaptionbox{AITV (full)}{\includegraphics[width = 1.35in]{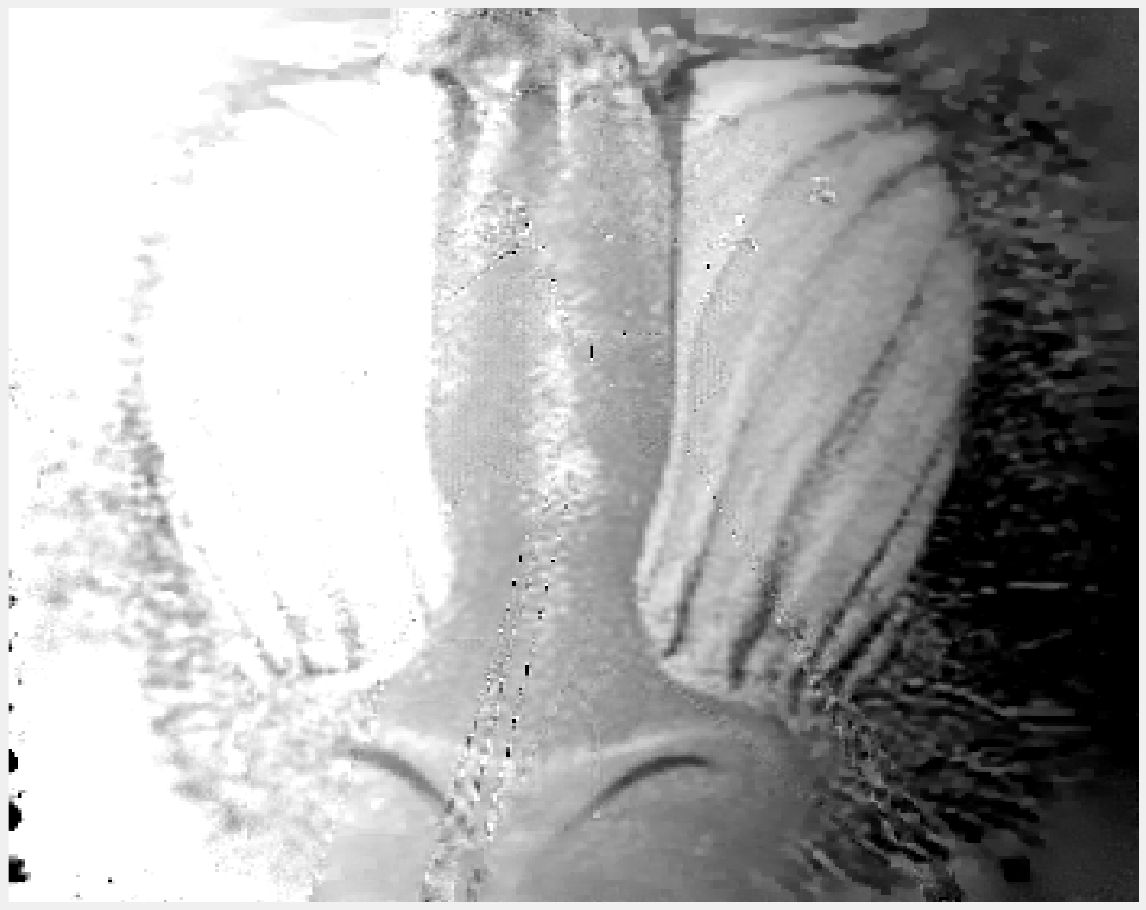}} & \subcaptionbox{DR}{\includegraphics[width = 1.35in]{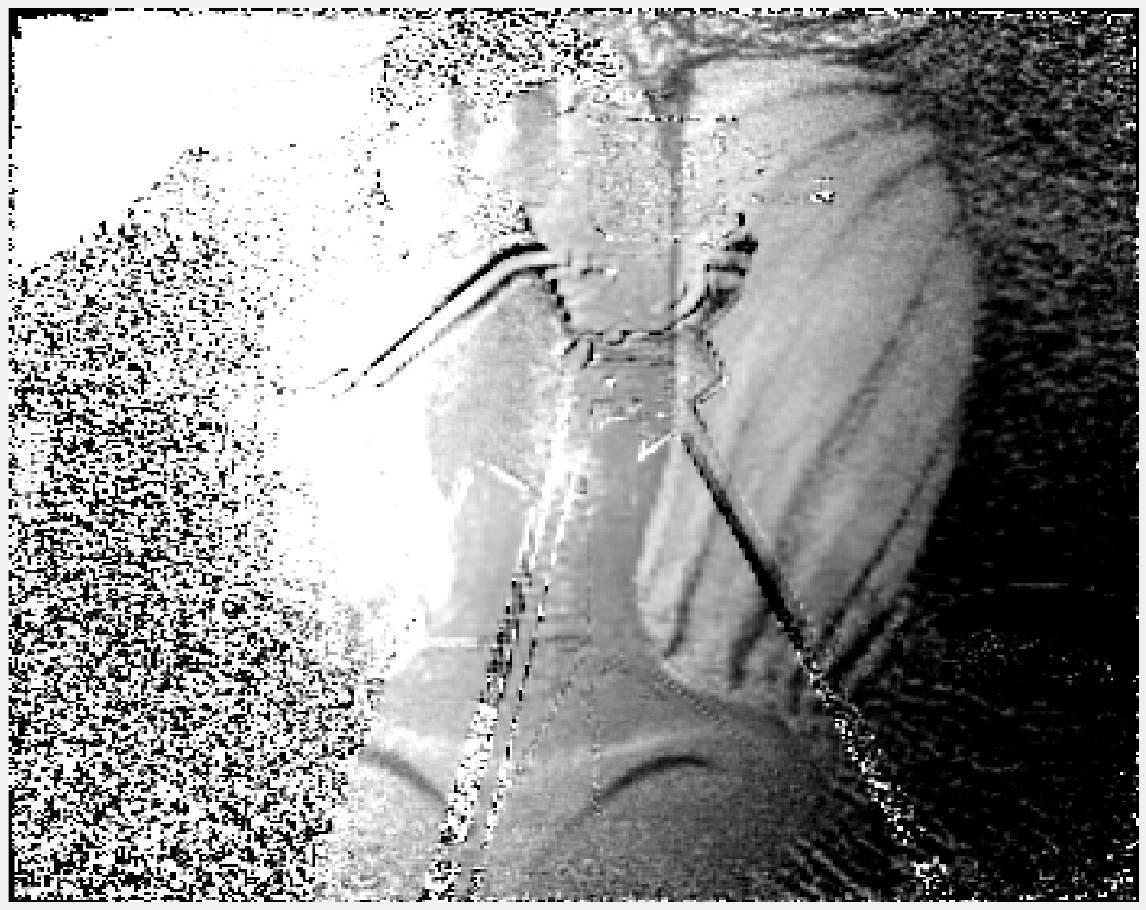}}& \subcaptionbox{rPIE}{\includegraphics[width = 1.35in]{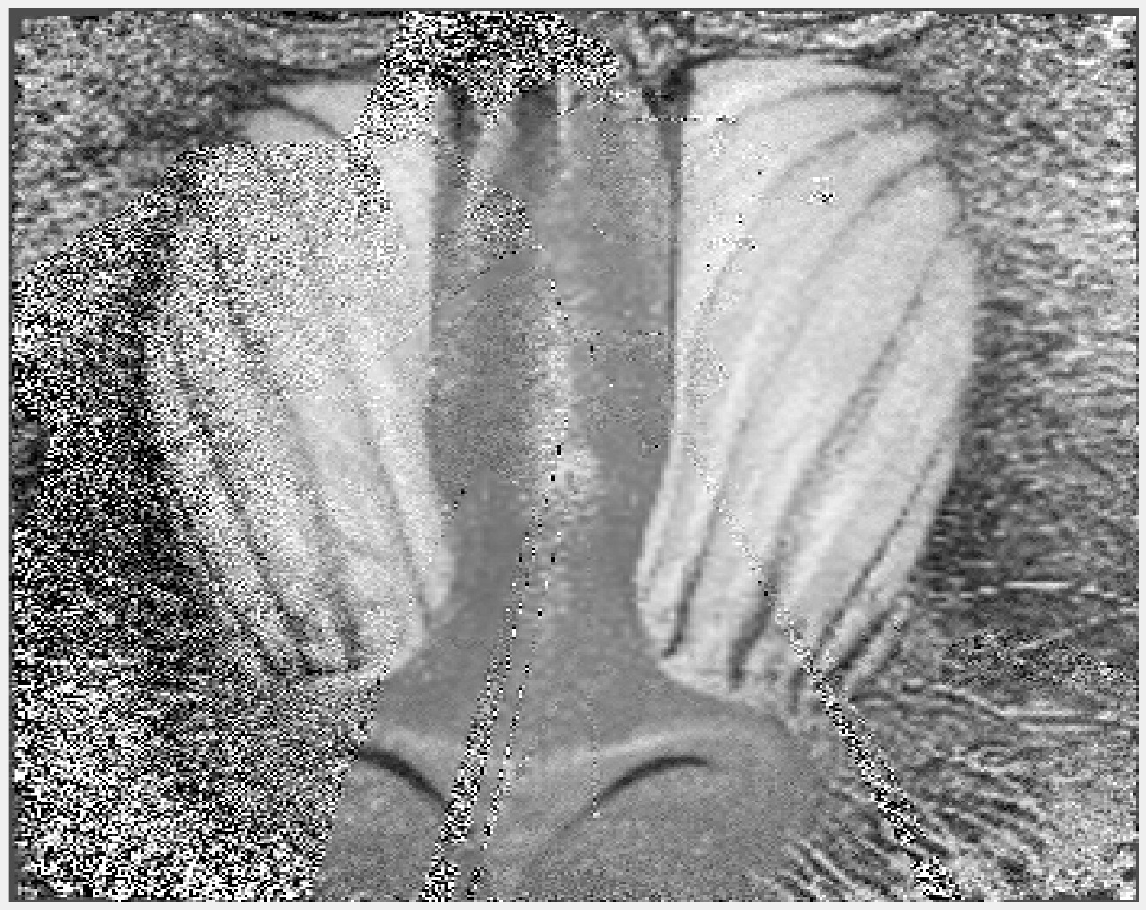}}   & \subcaptionbox{PHeBIE}{\includegraphics[width = 1.35in]{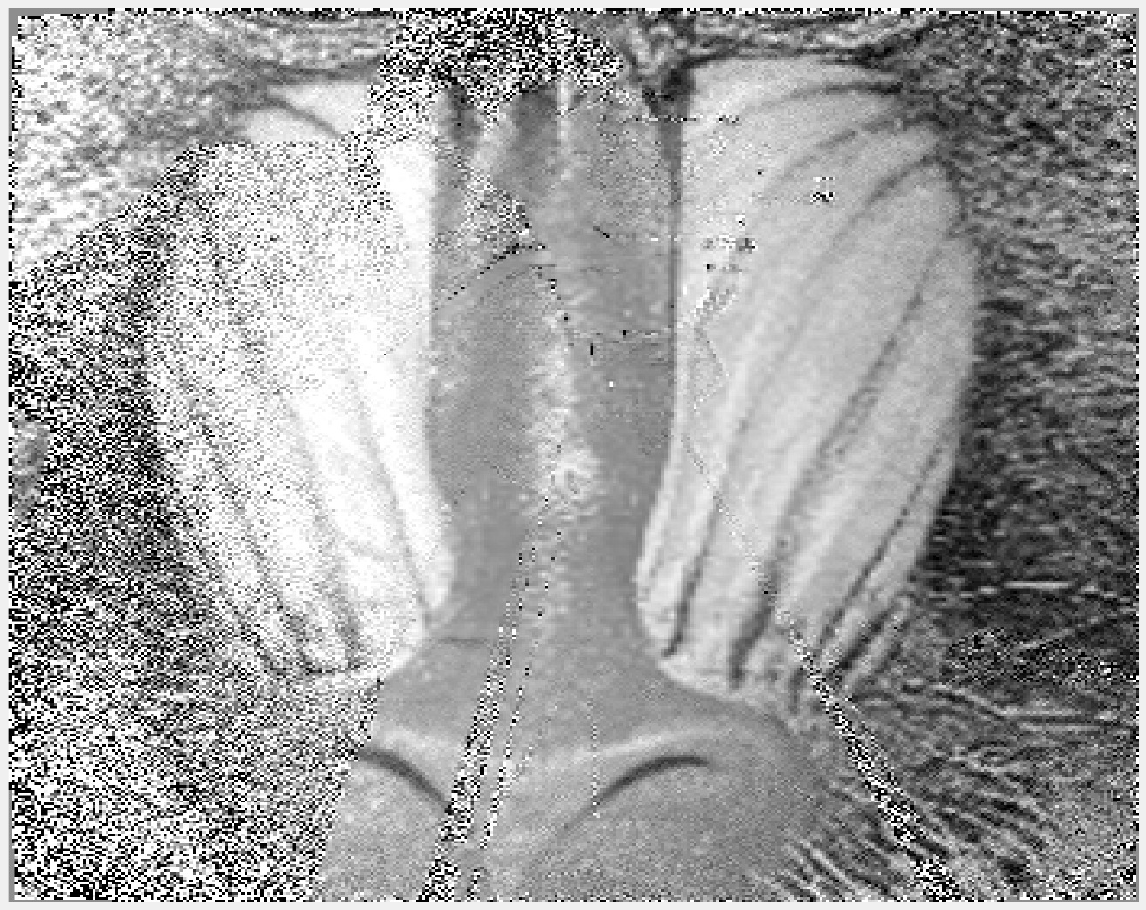}}
		\end{tabular}}
		\caption{Reconstructions of the blind case for the Gaussian noise with SNR = 40. Top two rows: reconstructions of Figures \ref{fig:chip_mag}-\ref{fig:chip_phase}; bottom two rows: reconstructions of Figs. \ref{fig:cameraman_mag}-\ref{fig:cameraman_phase}. }
		\label{fig:blind_agm}
	\end{minipage}
     \centering
         \includegraphics[width=0.5\textwidth]{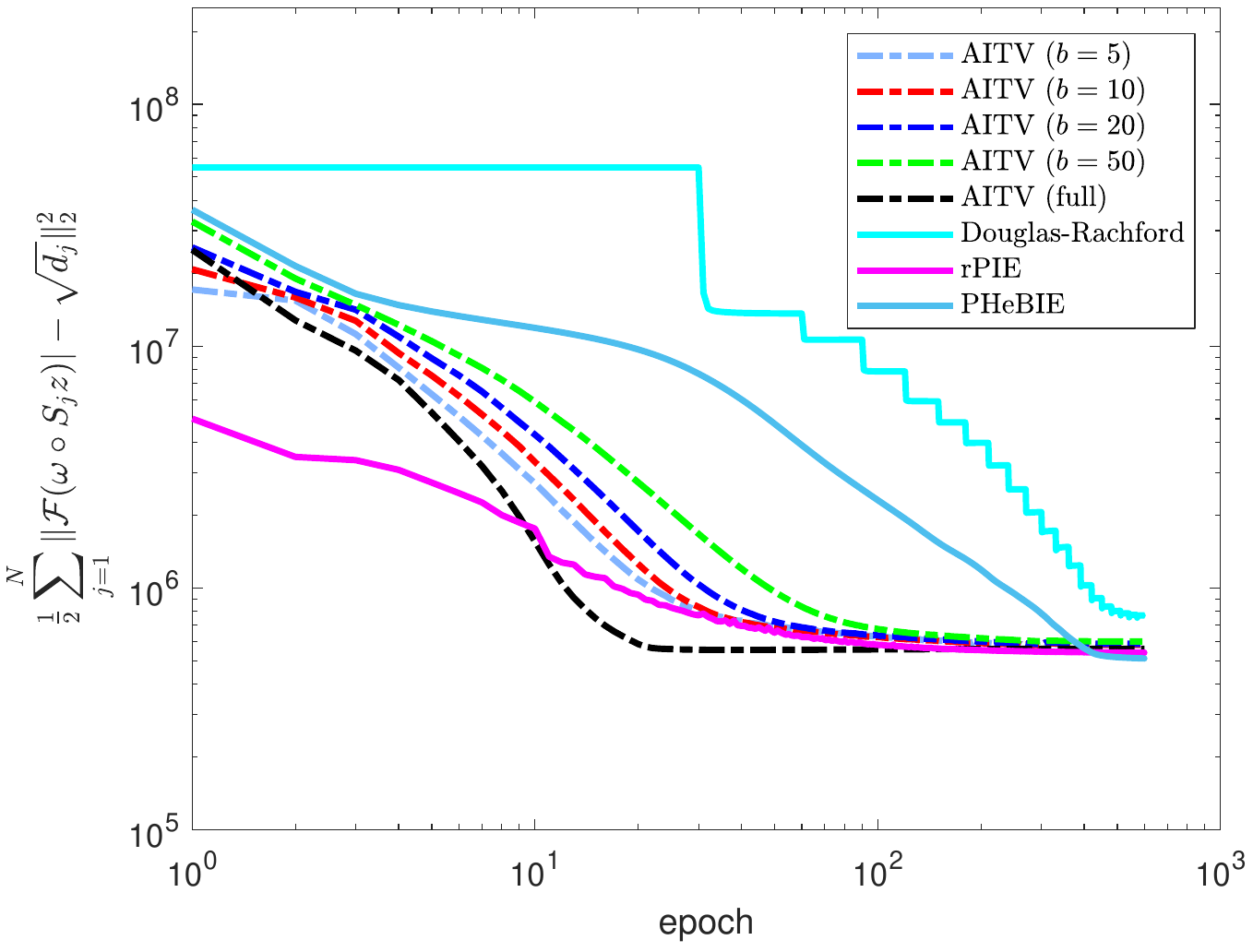}
              \caption{Amplitude Gaussian metric plotted across 600 epochs for the blind algorithms applied to the complex image given by Figures \ref{fig:cameraman_mag}-\ref{fig:cameraman_phase}, where the measurements are corrupted with Gaussian noise with SNR = 40.}
     \label{fig:agm_convergence}
\end{figure}\clearpage}%
For both the non-blind and blind cases, we examine the case when the noisy measurements have $\text{SNR} = 40$, so we set the regularization parameter $\lambda = 10.0$. Table \ref{tab:gaussian_result} records the SSIMs of the magnitude and phase components of the test images for the non-blind and blind cases. For all cases, DR, rPIE, and PHeBIE yield the worst magnitude SSIMs, and AITV attains better magnitude and phase SSIMs than its corresponding isoTV counterpart. The stochastic AITV ($b=10,20$) has slightly lower magnitude SSIMs by at most $0.03$ than the best results obtained from the deterministic, full-batch AITV. In fact, stochastic AITV attains the second best magnitude SSIMs in three out of the four cases considered. On the other hand, stochastic AITV ($b=10,20$) has the best phase SSIMs, such as outperforming their deterministic counterparts by up to 0.19 for the the blind case of the cameraman/baboon image. In general, the stochastic algorithm does best in recovering the phase components while recovering the magnitude components with comparable quality as the deterministic algorithm. 
\afterpage{
\begin{figure}[t!]
    \centering
    \includegraphics[scale=0.5]{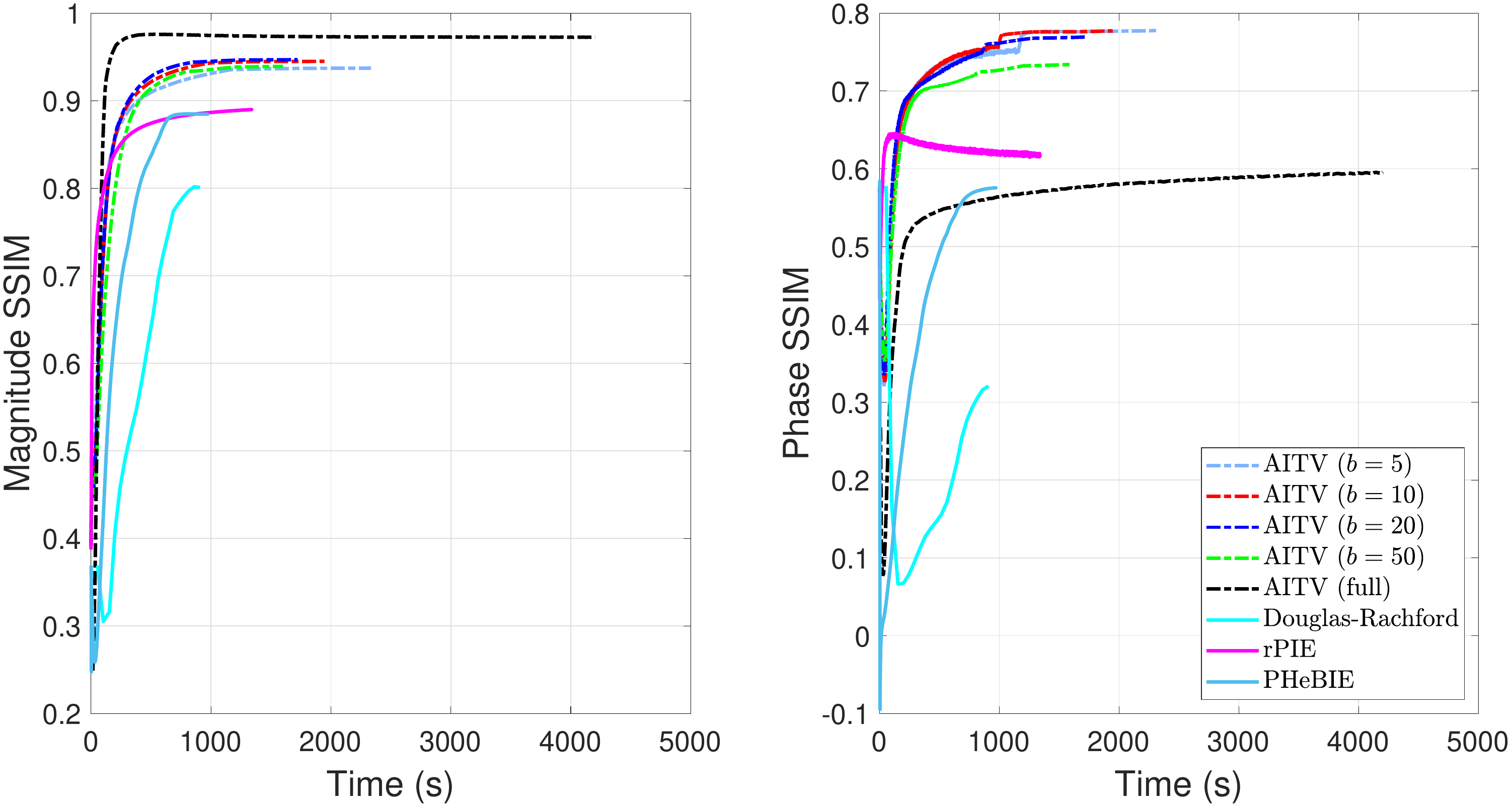}
    \caption{Computational time (seconds) vs.\ magnitude/phase SSIMs for the blind algorithms. Each algorithm is ran for 600 epochs for the complex image given by Figures \ref{fig:cameraman_mag}-\ref{fig:cameraman_phase}, where the measurements are corrupted with Gaussian noise with SNR =40.}
    \label{fig:agm_compute}
\end{figure}
\begin{figure}[t!]
    \centering
    \includegraphics[scale=0.5]{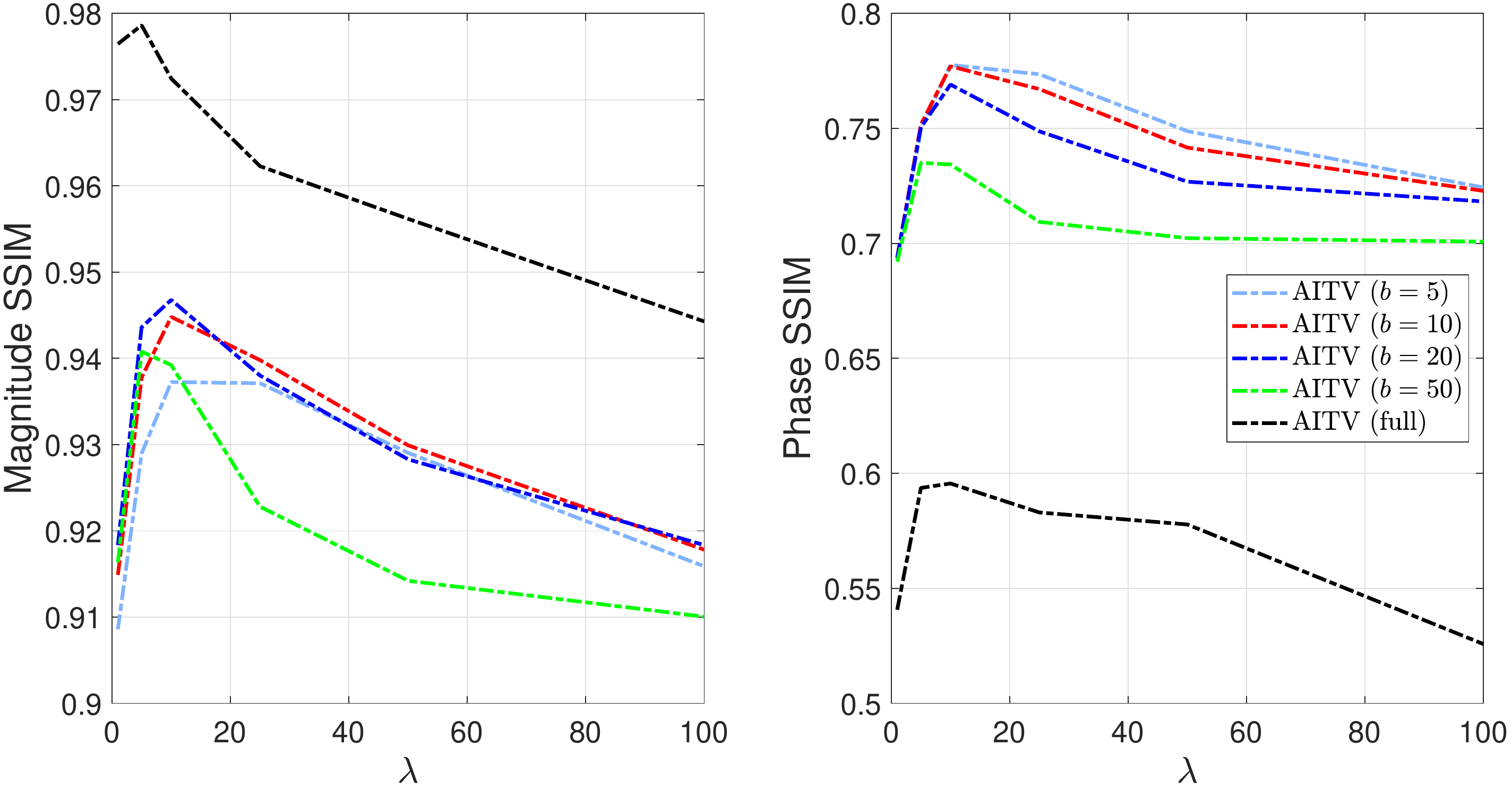}
    \caption{Changes in magnitude and phase SSIMs with respect to the AITV regularization parameter $\lambda \in \{1.0, 5.0, 10.0, 25.0, 50.0, 100.0\}$. Each AITV algorithm is ran for 600 epochs in the blind case for the complex image given by Figures {\ref{fig:cameraman_mag}}-{\ref{fig:cameraman_phase}}, where the measurements are corrupted with Gaussian noise with SNR = 40.}
    \label{fig:agm_lambda}
\end{figure}
\begin{figure}[t!]
\centering
\includegraphics[scale=0.5]{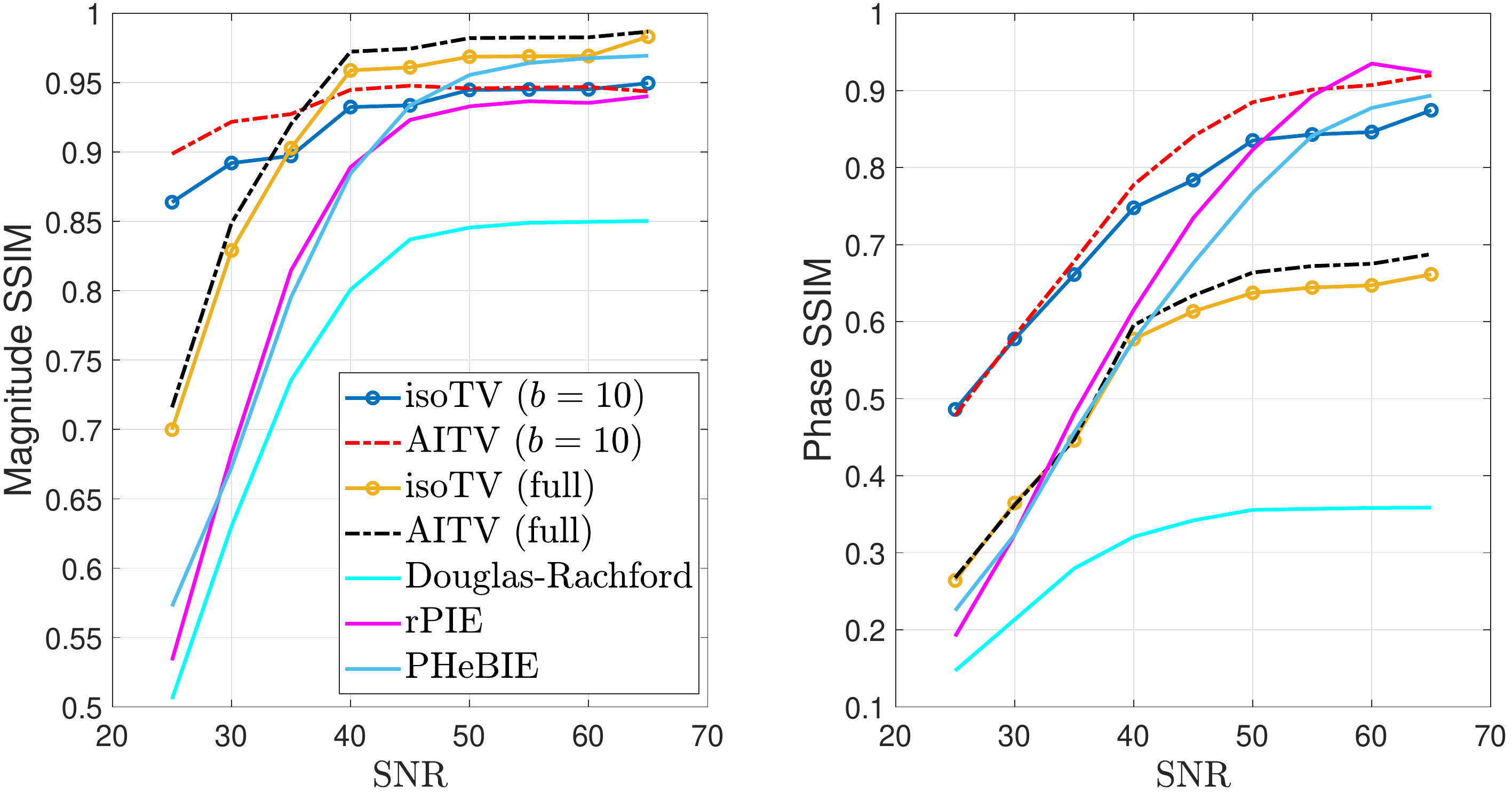}
\caption{Magnitude and phase SSIMs over different Gaussian noise level for the complex image given by Figure \ref{fig:cameraman_mag}-\ref{fig:cameraman_phase} for the blind case.}\label{fig:robust_agm}
\end{figure}}%

The reconstructed images for the non-blind experiments are presented in Figure \ref{fig:agm_nonblind}. DR, rPIE, PHeBIE, and the stochastic algorithms have artifacts in all four corners of the magnitude images  because the corners are scanned significantly less than in the middle of the image. However, the deterministic AITV has no artifacts because \eqref{eq:z_opt} is solved exactly for the image solution $z$. As a result, it has higher magnitude SSIMs than their stochastic counterparts. Nevertheless, the stochastic algorithms yield better phase images with less noise artifacts than any other algorithms. For example, the phase images of Figure \ref{fig:cameraman_phase} reconstructed from stochastic isoTV and AITV have the least amount of artifacts from the cameraman in the magnitude component.

Figure \ref{fig:blind_agm} shows the results of the blind algorithms. The phase images reconstructed by the deterministic AITV, DR, rPIE and PHeBIE are significantly worse than the stochastic algorithms. For example in Figure \ref{fig:chip_phase}, the contrasts of the reconstructed images by the deterministic AITV, DR, and PHeBIE are inconsistent as they become darker from left to right while the contrasts are more consistent with the stochastic algorithms. For Figure \ref{fig:cameraman_phase}, the stochastic algorithms perform the best in recovering the phase image while deterministic AITV is unable to recover the left half of the image and DR, rPIE, and PHeBIE have strong remnants of the cameraman present. Like in the non-blind case, stochastic AITV reconstructs the phase image the best. Overall, recovering the phase component in the blind case is more difficult than the nonblind case because of the artifacts created by the inherent ambiguities encountered when recovering both the probe and the object \mbox{\cite{bendory2019blind, fannjiang2019raster,fannjiang2020blind}}. 

In Figure \ref{fig:agm_convergence}, we examine the convergence of the blind algorithms applied to the cameraman/baboon image by recording their AGM values for each epoch. We omit the convergence curves of isoTV since they are similar to their AITV counterparts. Overall, the curves for our proposed stochastic algorithms are decreasing, validating the numerical convergence of Algorithm \ref{alg:stochastic_ADMM_blind} with AGM fidelity. However, their curves are slightly above the deterministic ADMM algorithm and rPIE. The reason why rPIE outperforms the AITV algorithms is because it seeks to only minimize AGM while the AITV algorithms minimize a larger objective function given by \eqref{eq:blind_Lagrangian}. Overall, after several hundred epochs, our proposed stochastic algorithms can give comparable  AGM values as the deterministic AITV and rPIE algorithms.

Figure \ref{fig:agm_compute} shows the magnitude/phase SSIMs of the cameraman/baboon image with respect to computational time in seconds when running each blind algorithm for 600 epochs. Although DR, rPIE, and PHeBIE finish running 600 epochs faster than the AITV algorithms, the stochastic AITV algorithms can attain magnitude SSIM higher than 0.90 and phase SSIM higher than 0.70 in about 500 seconds, thereby outperforming DR, rPIE, and PHeBIE in faster time. Moreover, because the stochastic AITV algorithms solve  \eqref{eq:full_omega_subprob} and \eqref{eq:z_opt} inexactly by SGD, they finish running 600 epoches faster than the deterministic AITV by at least 2000 seconds or 33.33 minutes. Apparently, larger batch size $b$ allows the stochastic algorithm to finish faster.

We examine the effect of the AITV regularization parameter $\lambda$ in \eqref{eq:AITV_blind_prob} on the image reconstruction quality. The value of $\lambda$ is adjusted depending on the amount of noise corrupting the measurements. If $\lambda$ is chosen to be small, AITV regularization will have a minimal impact on the image reconstruction, resulting in a potentially noisy image. However, if $\lambda$ is chosen to be large, AITV will oversmooth the reconstructed image. Therefore, a moderate value of $\lambda$ is recommended. For the AITV algorithms, we test $\lambda \in \{1.0, 5.0, 10.0, 25.0, 50.0, 100.0\}$ for the blind case on the complex image given by Figures \ref{fig:cameraman_mag}-\ref{fig:cameraman_phase}, whose measurements are corrupted by Gaussian noise of SNR = 40. According to Figure \ref{fig:agm_lambda}, the magnitude/phase SSIMs have a concave relationship with the AITV regularization parameter $\lambda$, where the SSIMs appear to peak at $\lambda = 10$ for most AITV algorithms.

Lastly, we analyze the robustness of the blind algorithms applied to the cameraman/baboon image whose measurements are corrupted by different levels of Gaussian noise, from SNR 25 to 65. The regularization parameter $\lambda$ is adjusted for different noise level of the measurements: $\lambda = 100$ for SNR = 25; $\lambda  =50$ for SNR = 30, 35; $\lambda = 10$ for SNR =40, 45; $\lambda =5$ for SNR = 50, 55, 60; and $\lambda = 3$ for SNR = 65. The algorithms ran for 600 epochs and the magnitude and phase SSIMs across different SNRs are plotted in Figure \ref{fig:robust_agm}. For $\text{SNR} \geq 40$, the deterministic isoTV and AITV algorithms have the best magnitude SSIMs than the other algorithms while their stochastic counterparts have slightly lower SSIMs. When $\text{SNR} < 40$, the stochastic algorithms perform the best. In fact, stochastic AITV $(b=10)$ has magnitude SSIM at least 0.90 across all noise levels considered. For the phase image, the stochastic algorithms have the highest SSIMs up to SNR = 55. For $\text{SNR} \geq 60$, the rPIE algorithm has the best phase SSIM while stochastic AITV has the second best. Overall, stochastic AITV is the most stable across different levels of Gaussian noise.

\subsection{Poisson noise}
\begin{table}[t!!]
      \centering
    \caption{SSIM results of the algorithms applied to the Poisson corrupted measurements with $\zeta = 0.01$ (SNR $\approx 44$ for Figures \mbox{\ref{fig:chip_mag}-\ref{fig:chip_phase}}; SNR $\approx 40$ for Figures \mbox{\ref{fig:cameraman_mag}-\ref{fig:cameraman_phase}}). The stochastic algorithms (e.g., AITV and isoTV, $b \in \{5, 10, 20, 25\}$) are ran three times to obtain the average SSIM values. \textbf{Bold} indicates best value; \underline{underline} indicates second best value.}
    \label{tab:ipm_result}
   \scriptsize
    \begin{tabular*}{\textwidth}{l|cc||cc|cc||cc}
    \hline
    & \multicolumn{4}{c|}{Non-blind} &  \multicolumn{4}{c}{Blind}\\ \hline
    &  \multicolumn{2}{c||}{Chip} &  \multicolumn{2}{c|}{Cameraman/Baboon}  &  \multicolumn{2}{c||}{Chip} &  \multicolumn{2}{c}{Cameraman/Baboon}  \\
        ~ & \makecell{mag.\\ SSIM} & \makecell{phase\\ SSIM} & \makecell{mag.\\ SSIM} & \makecell{phase\\ SSIM} &  \makecell{mag.\\ SSIM} & \makecell{phase\\ SSIM}& \makecell{mag.\\ SSIM} & \makecell{phase\\ SSIM} \\ \hline
        DR & 0.8523	&0.8455&	0.8704	&0.5043& 0.8431&	0.7387&	0.7630	&0.2529
 
 \\ \hline
 rPIE & 0.0206 & 0.1400&0.0701	&0.2136& 0.0213 & 0.1416	&	0.0700	& 0.2137
 
 \\ \hline
         PHeBIE & 	0.9404&0.9398	&0.9271&0.6791
         & 0.9280&	0.9082&	0.8678& 0.5470
 
 \\ \hline
        isoTV $(b=5)$ &0.9491&	0.9130&	0.9364&	0.7063&	0.9402&	0.8964&	0.9256&	0.7075

 \\ \hline
        isoTV $(b=10)$ & 0.9411	&0.9164&	0.9338&	0.6892	&0.9343&0.9029&0.9221	&0.6911
 \\ \hline
        isoTV $(b=20)$ & 0.9377&	0.9319	&0.9321&	0.6788&	0.9339&	0.9236	&0.9186&	0.6810

 \\ \hline
        isoTV $(b=25)$ & 0.9379&	0.9382	&0.9313&	0.6768&	0.9355&	0.9306	&0.9175&	0.6788

 \\ \hline
        isoTV (full batch) & \underline{0.9767}&	0.9590&	\underline{0.9773}&	0.7093
 & \underline{0.9655}&	0.9192&	\underline{0.9588}&	0.4920
\\ \hline
        AITV $(b=5)$ &0.9594&	0.9598&	0.9384&	\textbf{0.7237}&	0.9484&	\textbf{0.9551}&	0.9298&	\textbf{0.7384}

 \\ \hline
        AITV $(b=10)$ & 0.9655&	0.9634&	0.9447&	\underline{0.7168}&	0.9551&	\underline{0.9548}&	0.9381&	\underline{0.7303}

 \\ \hline
        AITV $(b=20)$ & 0.9645&	\underline{0.9641}&	0.9473&	0.7012&	0.9559&	0.9543&	0.9373&	0.7082
 \\ \hline
        AITV $(b=25)$ & 0.9629&	0.9634&	0.9462&	0.6957&	0.9549&	0.9535&	0.9347&	0.7012

 \\ \hline
        AITV (full batch) & \textbf{0.9803}	&\textbf{0.9644}&	\textbf{0.9782}	&0.7084
 & \textbf{0.9741}&	0.9354&	\textbf{0.9671}&	0.4975 \\
 \hline
    \end{tabular*}
\end{table}

\begin{figure}[t!]
    \centering
    \includegraphics[scale=0.5]{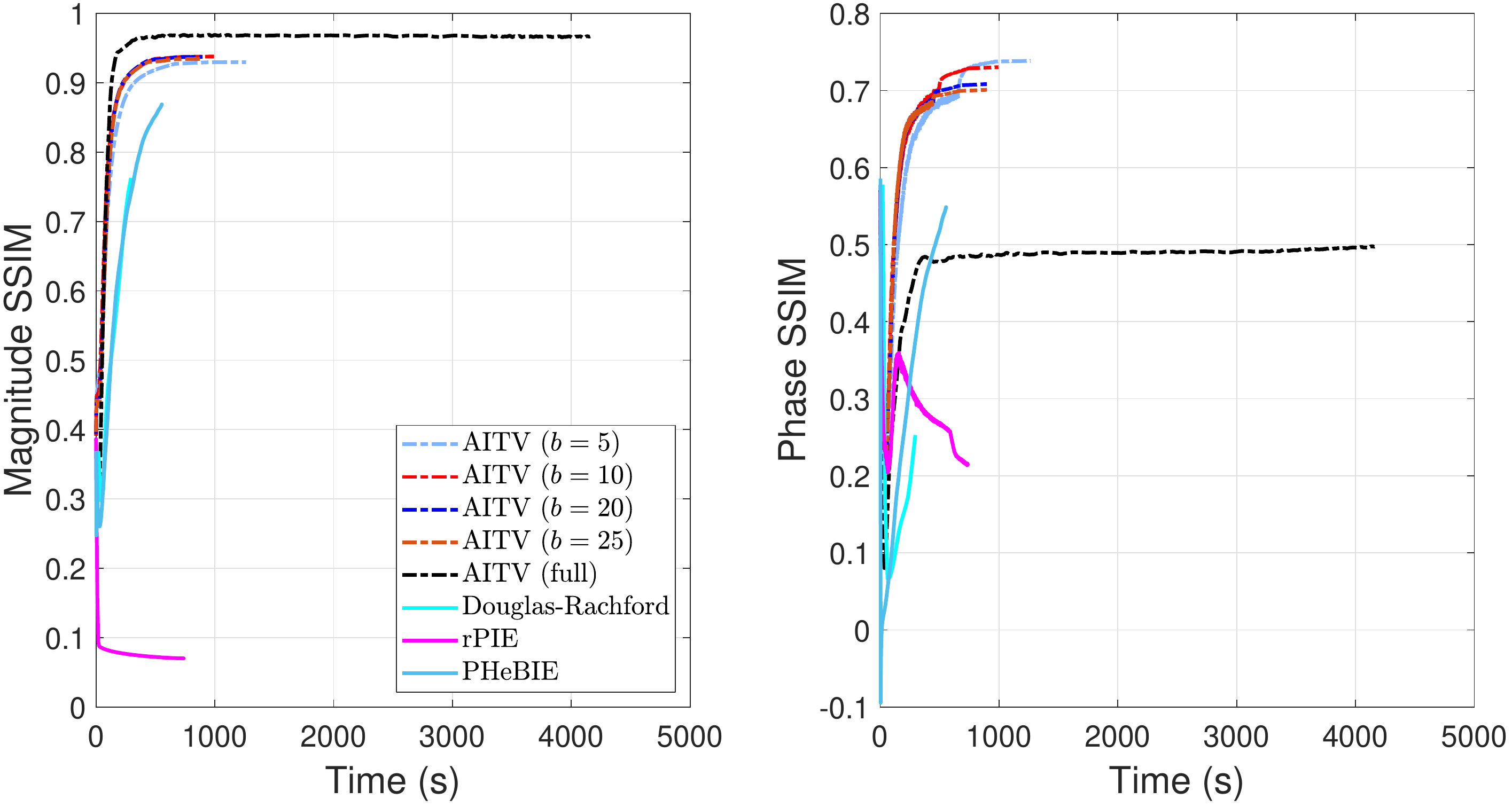}
    \caption{Total computational time (seconds) vs.\ magnitude/phase SSIMs for the blind algorithms. Each algorithm is ran for 300 epochs for the complex image given by Figures \ref{fig:cameraman_mag}-\ref{fig:cameraman_phase}, where the measurements are corrupted with Poisson noise with $\zeta = 0.01$ (SNR $\approx 40$)}.
    \label{fig:ipm_compute}
\end{figure}

\begin{figure}[t!]
    \centering
    \includegraphics[scale=0.5]{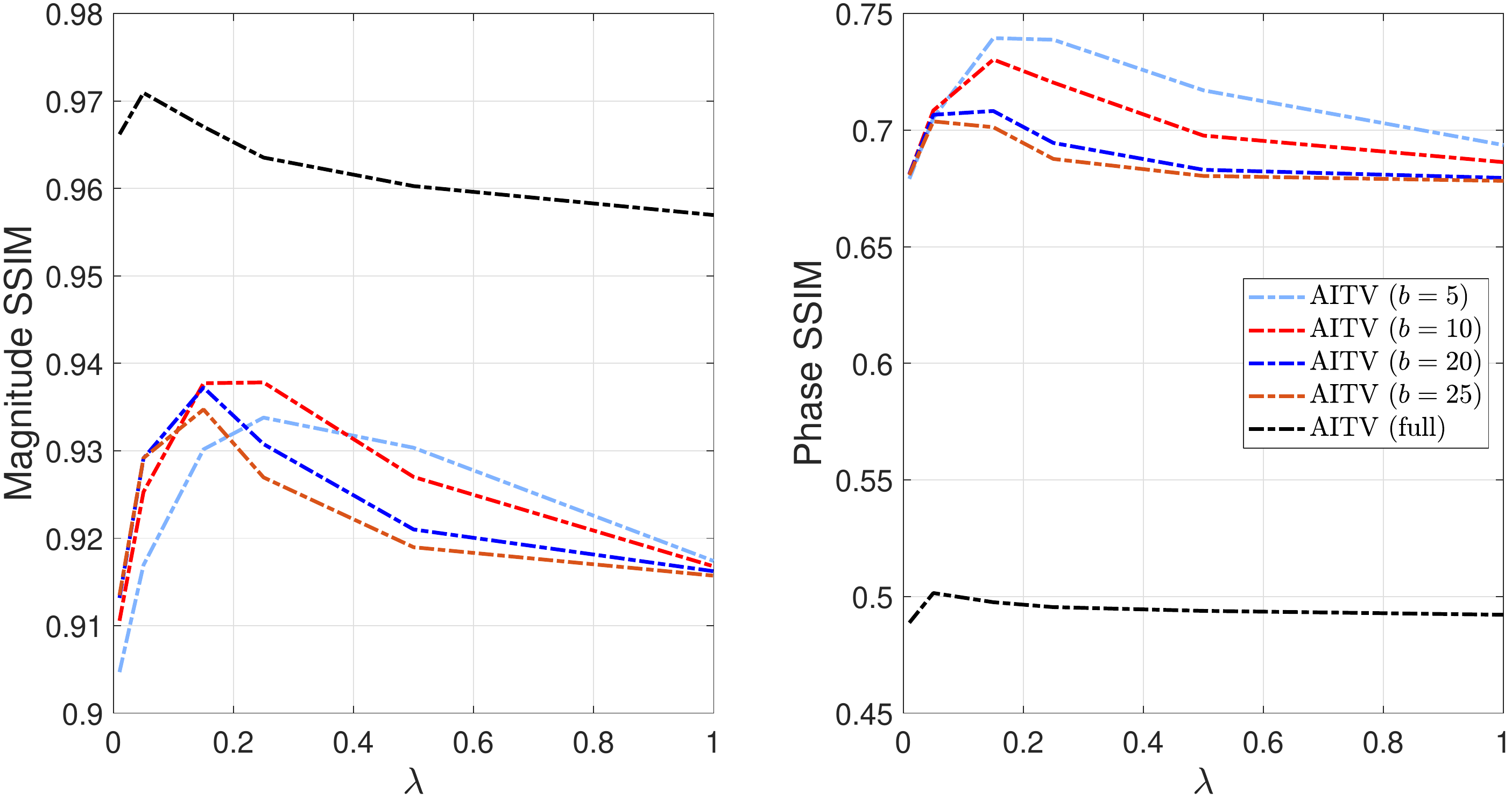}
    \caption{Changes in magnitude and phase SSIMs with respect to the AITV regularization parameter $\lambda \in \{0.01, 0.05, 0.15, 0.25, 0.50, 1.00\}$. Each AITV algorithm is ran for 300 epochs in the blind case for the complex image given by Figures {\ref{fig:cameraman_mag}}-{\ref{fig:cameraman_phase}}, where the measurements are corrupted with Poisson noise with $\zeta = 0.01$ (SNR $\approx 40$).}
    \label{fig:ipm_lambda}
\end{figure}

For both the non-blind and blind case, we examine the measurements corrupted with Poisson noise with $\zeta = 0.01$ (SNR $\approx 44$ for Figures \ref{fig:chip_mag}-\ref{fig:chip_phase}; SNR $\approx 40$ for Figures \ref{fig:cameraman_mag}-\ref{fig:cameraman_phase}) according to 
\eqref{eq:noise_d_j}. We set the regularization parameter $\lambda  = 0.15$. The numerical results are recorded in Table \ref{tab:ipm_result}.

Almost all cases, deterministic AITV attains the highest magnitude SSIMs and stochastic AITV attains the highest phase SSIMs while rPIE performs the worst in reconstructing images from Poisson-corrupted measurements. Using AITV over isoTV, we observe improvement in both magnitude and phase SSIMs. Although the stochastic algorithms have lower magnitude SSIMs than their deterministic counterparts, the difference is at most 0.04 for AITV and at most 0.05 for isoTV. Moreover, the magnitude SSIMs from the stochastic algorithms remain above 0.91. Similar to the Gaussian noise case, stochastic AITV reconstructs the phase image well while recovering the magnitude image with satisfactory quality.  

Figure \ref{fig:ipm_compute} shows the magnitude/phase SSIMs of the cameraman/baboon image with respect to computational time in seconds when running each blind algorithm for 300 epochs. Like for the Gaussian case, the stochastic AITV algorithms can recover faster the complex image with better quality in both the magnitude and phase components than DR, rPIE, and PHeBIE. Moreover, the stochastic AITV algorithms finish running 300 epochs significantly faster than the deterministic AITV algorithm. In fact, with stochastic AITV, one can recover within 1000 seconds the complex image with magnitude SSIM higher than 0.90 and phase SSIM higher than 0.70. Overall, the stochastic algorithms are computationally efficient in recovering complex images with satisfactory quality. 

We perform sensitivity analysis of the AITV regularization parameter $\lambda$ on the magnitude/phase SSIMs, where we test $\lambda \in \{0.01, 0.05, 0.15, 0.25, 0.50, 1.00\}$. The AITV algorithms ran for 300 epochs in the blind case on the complex image given by Figures \ref{fig:cameraman_mag}-\ref{fig:cameraman_phase}. The measurements are corrupted by Poisson noise with $\zeta = 0.01$ (SNR $\approx 40$). Figure {\ref{fig:ipm_lambda}} shows the impact of the AITV regularization parameter $\lambda$ on the magnitude/phase SSIMs. Like for the AGM case, the SSIMs are concave with respect to $\lambda$, where they peak at $\lambda = 0.15$.

\begin{figure}[t!!!]
\centering
\includegraphics[scale=0.5]{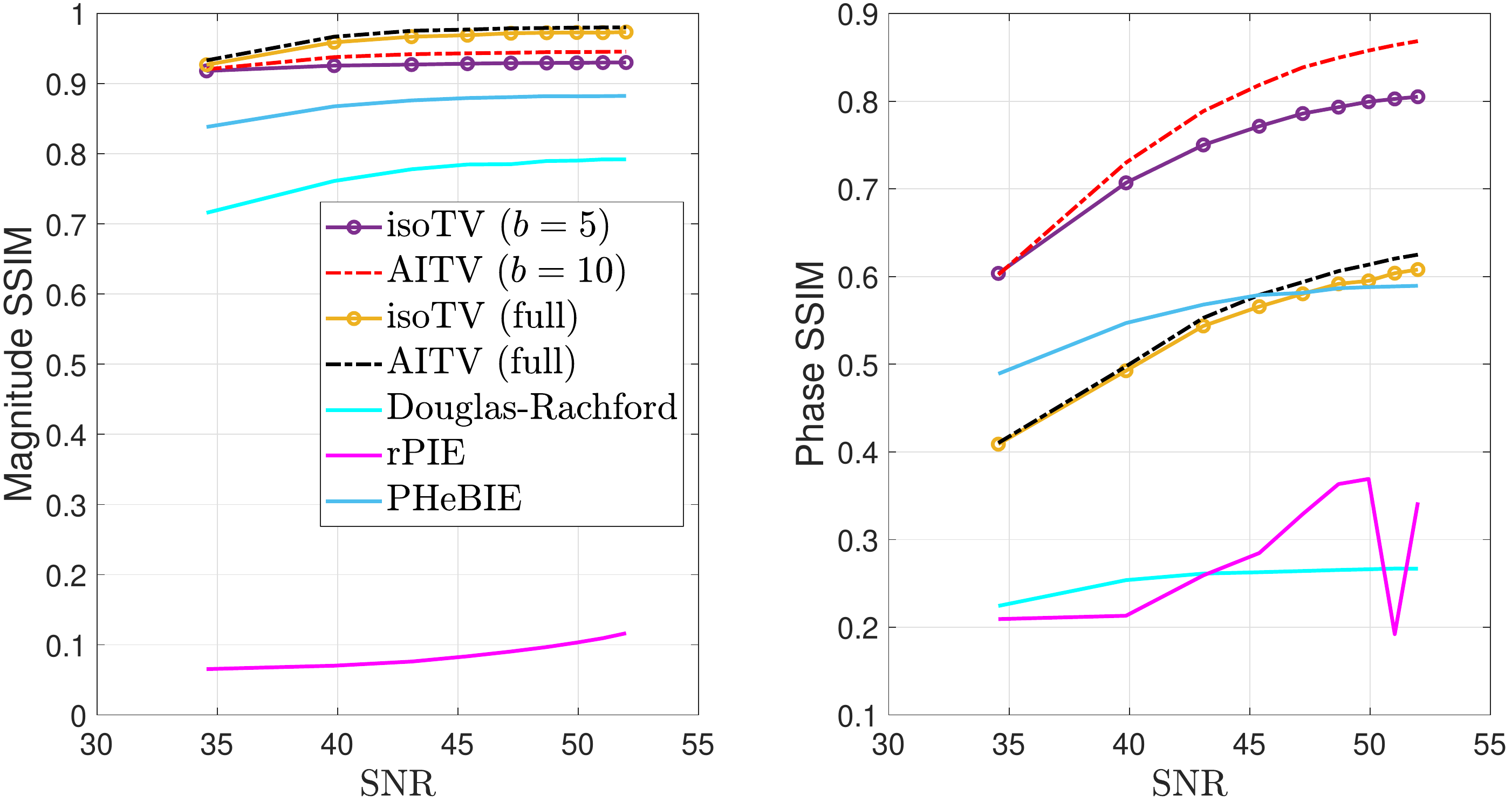}
\caption{Magnitude and phase SSIMs over different Poisson noise level for the complex image given by Figure \ref{fig:cameraman_mag}-\ref{fig:cameraman_phase} for the blind case.} \label{fig:robust_ipm}
\end{figure}

We examine the robustness of the blind algorithms on Figures \ref{fig:cameraman_mag}-\ref{fig:cameraman_phase} with different level of Poisson noise. The noise levels we examine are $\zeta \in \{0.005k\}_{k=1}^{9}$ (SNR $\approx 35, 40, 43, 45, 47, 49, 50, 51, 52$, respectively). We set the regularization parameter to be $\lambda = 15\times \zeta$. The magnitude and phase SSIMs across different Poisson noise levels are plotted in Figure \ref{fig:robust_ipm}. We observe that the deterministic algorithms yield the best magnitude SSIMs while the stochastic algorithms yield the best phase SSIMs. Both DR and rPIE yield the worst results for both magnitude and phase components. Although stochastic AITV yields the third best SSIMs for the magnitude image, its SSIMs are at least 0.90. Moreover, it has the best phase SSIMs, significantly more than its deterministic counterpart by about 0.20. In summary, stochastic AITV is a robust method across different levels of Poisson noise.

\section{Conclusion}\label{sec:conclude}
In this study, we present a novel approach for image ptychography utilizing AITV-regularized variational models. These models effectively handle measurements corrupted by Gaussian or Poisson noise. To address the challenges posed by a large number of measurements, we develop a stochastic ADMM algorithm with adaptive step sizes inspired by the inherent structure of ptychography (e.g., per-pixel illumination strength). Our proposed method demonstrates the ability to reconstruct high-quality images from severely corrupted measurements, with particular emphasis on accurately recovering the phase component. Theoretical convergence of the stochastic ADMM algorithm is established under certain conditions, and our numerical experiments confirm both convergence and computational efficiency.

For future research, we aim to design a globally convergent algorithm for the AITV-regularized ptychography model with box constraints and provide convergence analysis with weaker assumptions than in Theorem \ref{thm:kkt}. Additionally, we plan to explore the incorporation of variance-reduced stochastic gradient estimators such as SVRG \cite{johnson2013accelerating} and SARAH \cite{nguyen2017sarah}. These techniques have the potential to further enhance computational efficiency and improve the quality of image reconstruction. Lastly, we plan to develop a library of ptychography models regularized with other nonconvex variants of total variation, such as TV$^p (0<p<1)$ regularization \mbox{\cite{hintermuller2013nonconvex, lanza2016constrained}}, log total variation {\cite{zhang2020tv}}, $\ell_1/\ell_2$ on the gradient \mbox{\cite{wang2022minimizing, wang2021limited}} and transformed total variation {\cite{huo2022stable}}.

\section*{Acknowledgments}
This material is based upon work supported by the U.S. Department of Energy, Office of Science, under contract number DE-AC02-06CH11357. We thank the two reviewers for their valuable feedback in improving the manuscript.

\appendix
\section{Proof of Lemma \ref{lemma:prox_l1l2}} \label{sec:prox}
\begin{proof}If $x' = 0$, then it is trivial, so for the rest of the proof, we assume that $x' \neq 0$. 

Suppose $x^*$ is the optimal solution to \eqref{eq:prox_l1_l2}.  Because $\|\cdot\|_1-\alpha\|\cdot\|_{2}$ is rotation invariant, we only need to examine and expand the quadratic term in \eqref{eq:prox_l1_l2}. We see that
\begin{align*}
      \|x^*-x'\|_2^2 &=  \sum_{i=1}^{n^2} |(x^*)_i - (x')_i|^2 =  \sum_{i=1}^{n^2} \left( |(x^*)_i|^2 + |(x')_i|^2 -2 |(x^*)_i| |(x')_i| \cos \theta_i \right),
\end{align*}
where $\theta_i$ is the angle between the components $(x^*)_i$ and $(x')_i$. This term is minimized when $\theta_i = 0$ for all $i$. This means that $\text{sgn}(x^*)_i = \text{sgn}(x')_i$ for all $i$ such that $(x^*)_i, (x')_i \neq 0$, or otherwise $x^*$ would not be an optimal solution to \eqref{eq:prox_l1_l2}. When $(x^*)_i = 0$ or $(x')_i = 0$, we choose $c_i \in \{c' \in \mathbb{C}: |c'|\leq 1\}$ such that $\sign(x^*)_i = c_i = \sign(x')_i$. Overall, we have $\sign(x^*) = \sign(x')$.

Since $\sign(x^*) = \sign(x')$, we can simplify \eqref{eq:prox_l1_l2} to an optimization problem with respect to $|x|$ given by
\begin{align*} 
    |x^*| = \argmin_{\rho \in \mathbb{R}^n, \;\rho_i \geq 0} \|\rho\|_1 - \alpha \|\rho\|_2 + \frac{1}{2\lambda}\|\rho - |x'|\|_2^2.
\end{align*}
Next, we show that nonnegativity constraint is redundant, i.e., 
\begin{align}\label{eq:opt_equality}
    \min_{\rho \in \mathbb{R}^n, \;\rho_i \geq 0} \|\rho\|_1 - \alpha \|\rho\|_2 + \frac{1}{2\lambda}\|\rho - |x'|\|_2^2 = \min_{\rho \in \mathbb{R}^n} \|\rho\|_1 - \alpha \|\rho\|_2 + \frac{1}{2\lambda}\|\rho - |x'|\|_2^2.
\end{align}
Let $\rho^*$ be the optimal solution to the right-hand optimization problem in {\eqref{eq:opt_equality}} and $\mathcal{I} = \{i: (\rho^*)_i < 0\}$. Comparing the optimal values between $\rho^*$ and $|\rho^*|$, we have
\begin{align*}
\|\rho^*\|_1 - \alpha \|\rho^*\|_2 + \frac{1}{2\lambda}\|\rho^* - |x'|\|_2^2 \leq \|\rho^*\|_1 - \alpha \|\rho^*\|_2 + \frac{1}{2\lambda}\||\rho^*| - |x'|\|_2^2,
\end{align*}
which reduces to
\begin{align} \label{eq:contradict_ineq}
\sum_{i \in \mathcal{I}} \rho_i |(x')_i| \geq \sum_{i \in \mathcal{I}} |\rho_i| |(x')_i|.
\end{align}
However, if $\mathcal{I} \neq \emptyset$, then 
\eqref{eq:contradict_ineq} is invalid because the left-hand side is negative. As a result, $\mathcal{I} = \emptyset$, implying that $\rho_i \geq 0$ for all $i$, so the nonnegativity constraint of the right-hand optimization problem in \eqref{eq:opt_equality} is redundant. This follows that we have
\begin{align*}
|x^*| = \argmin_{\rho \in \mathbb{R}^n} \|\rho\|_1 - \alpha \|\rho\|_2 + \frac{1}{2\lambda}\|\rho - |x'|\|_2^2.
\end{align*}
Hence, by applying \mbox{\cite[Lemma 1]{lou2018fast}} to the optimization problem, we establish the following:
\begin{enumerate}
\item[(1)] When $\|x'\|_\infty >\lambda$, 
\begin{align*}
|x^*| = (\|\xi\|_2 + \alpha \lambda) \frac{|\xi|}{\|\xi\|_2}.
\end{align*}
\item[(2)] When $\|x'\|_\infty =\lambda$, $|x^*|$ is an optimal solution if and only if it satisfies $x^*_i=0$ if $|(x')_i|<\lambda$, and $\|x^*\|_2=\alpha\lambda$. 
\item[(3)] When $(1-\alpha)\lambda<\|x'\|_\infty <\lambda$, $|x^*|$ is an optimal solution if and only if it is a 1-sparse vector satisfying $(x^*)_i=0$ if $|(x')_i|<\|x'\|_\infty$, and $\|x^*\|_2=\|x'\|_\infty+(\alpha-1)\lambda$. 
\item[(4)] When $\|x'\|_\infty\leq(1-\alpha)\lambda$, $|x^*|=0$.
\end{enumerate}
When $(1-\alpha) \lambda < \|x'\|_{\infty} \leq \lambda$, we show that for a selected index $i \in \argmax_j (|(x')_j|)$  we have
\begin{align*}
|(x^*)_j| = \begin{cases}
|(x')_j| + (\alpha-1)\lambda, & \text{ if } j = i,\\
0, & \text{ if } j \neq i
\end{cases}
\end{align*}
as an optimal solution. By construction in the statement of Lemma \ref{lemma:prox_l1l2}, $|x^*|$ is a 1-sparse vector satisfying $(x^*)_j = 0$ when $|(x')_j| <  \|x'\|_{\infty} \leq \lambda$. Since $|(x')_i| = \|x'\|_{\infty}$, we have
\begin{align*}
\|x^*\|_2 = |(x^*)_i| = | |(x')_i| +(\alpha - 1) \lambda | =  \begin{cases}
    \alpha \lambda, &\text{ if }
 \|x'\|_{\infty} = \lambda, \\
 \|x'\|_{\infty} + (\alpha - 1 )\lambda, & \text{ if }  \|x'\|_{\infty} < \lambda.
 \end{cases}
\end{align*}
By satisfying the conditions in (2) and (3) above, $|x^*|$ is an optimal solution. Lastly, by multiplying $|x^*|$ by $\sign(x^*)$, we obtain the desired results. 
\end{proof}
\section{Proofs of Section \ref{sec:convergence}}
\label{sec:appendix}
 Before proving our main results, we present preliminary tools necessary for the convergence analysis.
 \begin{definition}[\cite{li2023solving,rockafellar2009variational}]
Let $h:\mathbb{C}^{n^2} \rightarrow (-\infty, +\infty]$ be a proper and lower semicontinuouous function and $\text{dom h} \coloneqq \{x \in \mathbb{C}^{n^2}: h(x) < \infty\}$. 
\begin{enumerate}[label=(\alph*)]
    \item The Fr\'echet subdifferential of $h$ at the point $x \in \text{dom } h$ is the  set
    \begin{align*}
        \hat{\partial}h(x) = \left\{ v \in \mathbb{C}^{n^2}: \liminf_{y \neq x, y \rightarrow x} \frac{h(y)-h(x) - \mathbb{R}(\langle v, y-x \rangle)}{\|y-x\|} \geq 0 \right\}.
    \end{align*}
\item The limiting subdifferential of $h$ at the point $x \in \text{dom } h$ is the set
\begin{align*}
    \partial h(x) = \left\{ v \in \mathbb{C}^{n^2}: \exists \{(x^t,v^t)\}_{t=1}^{\infty} \text{ s.t. } x^t \rightarrow x, h(x^t) \rightarrow h(x),\; \hat{\partial}h(x^t) \ni v^t \rightarrow v \right\}.
\end{align*}
\end{enumerate}
\end{definition}
We note that the limiting subdifferential is closed \cite{rockafellar2009variational}:
\begin{align} \label{eq:subdiff_closed}
     (x^t,v^t) \rightarrow (x,v), h(x^t) \rightarrow h(x), v^t \in \partial h(x^t) \implies v \in \partial h(x).
\end{align}

\begin{lemma}[\cite{chung2001course}, pg.129, Exercise 7] \label{lemma:sum_converge}
    Let $\{X_j\}_{j=1}^{\infty}$ be a sequence of random variables. If $\sum_{j=1}^{\infty} \mathbb{E}[|X_j|] < \infty$, then $\sum_{j=1}^{\infty} X_j$ converges absolutely a.s.
\end{lemma}
\begin{lemma} \label{lemma:sampling}
Let $n^t$ be an index set chosen uniformly at random from all subsets of $\{1, \ldots, N\}$ with cardinality $b \leq N$ at iteration $t$. For any $j \in \{1, \ldots, N\}$, we have $j \in n^t $ for infinitely many $t$ a.s.
\end{lemma}
\begin{proof}
Since sampling $n^t$ from $\{1, \ldots, N\}$  for each $t$ is independent and $\sum_{t=1}^{\infty} \mathbb{P}(j \in n^t) = \sum_{t=1}^{\infty}\frac{b}{n} = \infty$, we apply Second Borel-Cantelli Lemma to obtain the desired result. 
\end{proof}
\subsection{Proof of Lemma  \ref{lemma:bounded_sequence}}
\begin{proof}
Because $\sum_{t=1}^{\infty} \|\Omega^{t+1} - \Omega^t\|_2^2 < \infty$, we have
    $\displaystyle \lim_{t \rightarrow \infty} \Lambda^{t+1} - \Lambda^t = 0$ and
    $\displaystyle \lim_{t \rightarrow \infty} y^{t+1} - y^t = 0$. From \eqref{eq:y_update}, we have
    \begin{align}\label{eq:v_equal_nabla_z}
        \lim_{t \rightarrow \infty} v^t - \nabla z^t &= 0 ,
    \end{align}
    so $\{v^t\}_{t=1}^{\infty}$ is bounded since $ \{z^t\}_{t=1}^{\infty}$ is bounded. Because  $\displaystyle \lim_{t \rightarrow \infty} \Lambda^{t+1} - \Lambda^t = 0$, it follows that $\{\Lambda^{t+1} - \Lambda^t\}_{t=1}^{\infty}$ is bounded. Combining the first case of \eqref{eq:lambda_update2} with the second case of \eqref{eq:u_closed_form}, we have
    \begin{align} 
    u_j^{t+1} = \begin{cases}
        \displaystyle \frac{1}{\beta_1} \left( \Lambda_j^{t+1} - \Lambda_j^{t} \right) + \mathcal{F}(\omega^{t+1} \circ S_j z^{t+1}), & \text{ if } j \in n^{t}, \\
        u_j^{t}, & \text{ if } j \not \in n^{t}.\label{eq:u_update3}
    \end{cases}
    \end{align}
    If $j \in n^{t}$, then we have
    \begin{align*}
        \|u_j^{t+1}\|_{2} \leq \frac{1}{\beta_1} \|\Lambda^{t+1} - \Lambda^t\|_2 + \|\mathcal{F}(\omega^{t+1} \circ S_j z^{t+1})\|_2.
    \end{align*}
    If $j \not \in n^t$, we have two cases. For one case, if we have $j \not \in n^{\tau}$ for all $0 \leq \tau < t$, then by \eqref{eq:u_update3}, we have $u_j^{t+1} = u_j^t = \ldots = u_j^0$. For the other case,
   there exists $t^* \coloneqq \max \{ \tau: \tau < t, j\in n^{\tau}\}$. Then by \eqref{eq:u_update3}, we have
    \begin{align*}
        \|u_j^{t+1}\|_2 = \|u_j^t\|_2 = \ldots = \|u_j^{t^*+1}\|_2 \leq  \frac{1}{\beta_1} \|\Lambda^{t^*+1} - \Lambda^{t^*}\|_2 + \|\mathcal{F}(\omega^{t^*+1} \circ S_j z^{t^*+1})\|_2.
    \end{align*}
    Altogether,  $\{u^t\}_{t=1}^{\infty}$ is bounded because $\{(\omega^t, z^t, \Lambda^{t+1} - \Lambda^t)\}_{t=1}^{\infty}$ is bounded.
By \eqref{eq:u_closed_form}, when $\mathcal{B}(\cdot, \cdot)$ is AGM, if $j \in n^t$, we have
\begin{align*}
    \|u_j^{t+1}\|_2 = \left\|\frac{\sqrt{d_j}+\beta_1\left| \mathcal{F}(\omega^t \circ S_j z^t) - \frac{1}{\beta_1} \Lambda_j^t \right|}{1+ \beta_1} \right\|_2 \geq \frac{\beta_1}{1+\beta_1}\left(\frac{1}{\beta_1} \left \| \Lambda_j^t\right\|_2 -  \|\mathcal{F}(\omega^t \circ S_j z^t)\|_2 \right),
\end{align*}
or equivalently,
\begin{align}\label{eq:u_bound_lambda}
    (1+\beta_1) \|u_j^{t+1}\|_2 + \beta_1 \|\mathcal{F}(\omega^t \circ S_j z^t)\|_2 \geq \|\Lambda_j^t\|_2.
\end{align}
Similarly, when $\mathcal{B}(\cdot, \cdot)$ is IPM, we have the same inequality as \eqref{eq:u_bound_lambda}. If $j \not \in n^{t}$, then as before, we have two cases. For the first case, if $j \not \in n^{\tau}$ for all $0 \leq \tau < t$, then by \eqref{eq:lambda_update2}, we have $\Lambda_j^{t} = \Lambda_j^{t-1} = \ldots= \Lambda_j^0$. Otherwise, there exists $t^*$, so
\begin{align*}
    \|\Lambda_j^t\|_2 =  \|\Lambda_j^{t-1}\|_2 = \ldots= \left\|\Lambda_j^{t^*+1} \right\|_2 
    &= \left \|\Lambda_j^{t^*} + \beta_1 \left(u_j^{t^*+1} - \mathcal{F}(\omega^{t^*+1} \circ S_jz^{t^*+1})\right) \right\|_2 \\
    &\leq \|\Lambda_j^{t^*}\|_2 + \beta_1 \left(\|u_j^{t^*+1}\|_2 + \|\mathcal{F}(\omega^{t^*} \circ S_j z^{t^*})\|_2  \right) \\
    &\leq (1+ 2 \beta_1)\|u_j^{t^*+1}\|_2 +2 \beta_1\|\mathcal{F}(\omega^{t^*} \circ S_j z^{t^*})\|_2. 
\end{align*}
Hence, we conclude that $\{\Lambda^t\}_{t=1}^{\infty}$ is bounded. Finally, we show that $\{y^t\}_{t=1}^{\infty}$ is bounded by examining \eqref{eq:v_sub3}. For each $i=1,\ldots, n^2$, we have two cases. For the first case, if there exists some $d \in \{x,y\}$ such that
\begin{align*}
 \left|(\nabla_d z^t)_i - \frac{(y_d^t)_i}{\beta_2} \right| \leq \frac{\lambda}{\beta_2},
\end{align*}
then by reverse triangle inequality, we have
\begin{align*}
|(y_d^t)_i| \leq \lambda + \beta_2 \left|(\nabla_d z^t)_i \right| \leq \lambda + \beta_2 \|\nabla z^t\|_{\infty} \leq \lambda + \beta_2  \left(\|v^{t+1}\|_{\infty} + \|\nabla z^t\|_{\infty}\right).
\end{align*}
Let
\begin{align*}
 \xi_i' = \text{sgn}\left((\nabla z^t)_i - \frac{(y^t)_i}{\beta_2} \right) \circ \max\left( \left|( \nabla z^t)_i - \frac{(y^t)_i}{\beta_2} \right| - \frac{\lambda}{\beta_2}, 0\right).
\end{align*}
For the other case, if there exists $d \in \{x,y\}$ such that 
\begin{align*}
 \left|(\nabla_d z^t)_i - \frac{(y_d^t)_i}{\beta_2} \right| > \frac{\lambda}{\beta_2},
\end{align*}
then by Lemma \ref{lemma:prox_l1l2}(i) and reverse triangle inequality, we have 
\begin{align*}
    \|(v^{t+1})_i\|_{\infty} & = \left(\| \xi'_i\|_2 + \frac{\alpha \lambda}{\beta_2} \right) \frac{\|\xi_i'\|_{\infty}}{\|\xi_i'\|_2}  \geq  \|\xi_i'\|_{\infty} \geq \left|(\nabla_d z^t)_i - \frac{(y_d^t)_i}{\beta_2} \right| - \frac{\lambda}{\beta_2} 
   \\ &\geq 
    \frac{1}{\beta_2}\left| (y_d^t)_i  \right|- \left| (\nabla_d z^t)_i  \right| - \frac{\lambda}{\beta_2}, 
\end{align*}
which follows that
\begin{align*}
 \left| (y_d^t)_i  \right| \leq \lambda + \beta_2  \left(\|(v^{t+1})_i\|_{\infty} + \left| (\nabla_d z^t)_i  \right|\right) \leq \lambda + \beta_2  \left(\|v^{t+1}\|_{\infty} + \|\nabla z^t\|_{\infty}\right).
\end{align*}
Altogether, we have
\begin{align*}
    \|y^t\|_{\infty} \leq  \lambda + \beta_2  \left(\|v^{t+1}\|_{\infty} + \|\nabla z^t\|_{\infty}\right),
\end{align*}
so $\{y^t\}_{t=1}^{\infty}$ is bounded because $\{(v^t, z^t)\}_{t=1}^{\infty}$ is bounded. Therefore, we establish that $\{(Z^t, \Omega^t)\}_{t=1}^{\infty}$ is bounded.

Because $\{(Z^t, \Omega^t)\}_{t=1}^{\infty}$ is bounded, it follows from {\eqref{eq:blind_Lagrangian}} that $\{(\mathcal{L}(\omega^t), \mathcal{L}(z^t))\}_{t=1}^{\infty}$ is bounded above. Then there exists a constant $B > 0$ such that $\mathcal{L}(\omega^t) \leq B$ and $\mathcal{L}(z^t) \leq B$ for each iteration $t \in \mathbb{N}$. Moreover, boundedness of $\{(Z^t, \Omega^t)\}_{t=1}^{\infty}$  implies that 
 $\{(\nabla_{\omega} \mathcal{L}(\omega^t), \nabla_z \mathcal{L}(z^t))\}_{t=1}^{\infty}$ is bounded because
\begin{align*}
    \left\| \nabla_{\omega} \mathcal{L}(\omega^t) \right\|_2 
    & \leq \beta_1 \sum_{j=1}^N \left\|(S_j z^t)^* \circ \left[\mathcal{F}^{-1}\left( u_j^{t+1} + \frac{\Lambda_j^t}{\beta_1} \right) - \omega^t \circ S_j z^t \right] \right\|_2 \\
    & \leq \beta_1 \sum_{j=1}^N \|z^t\|_{\infty} \left(\|u_j^{t+1}\|_2 + \frac{1}{\beta_1} \|\Lambda_j^t\|_2 + \|z^t\|_{\infty}\|\omega^t\|_2 \right)
\end{align*}
and
\begin{align*}
\|\nabla_z & \mathcal{L}(z^t)\|_2\\ \leq& \sum_{j=1}^N \left[\beta_1 \left\| (P_j^{t+1})^* \mathcal{F}^{-1}\left( u_j^{t+1} + \frac{\Lambda_j^t}{\beta_1} \right) -(P_j^{t+1})^* P_j^{t+1}z^{t+1}   \right\|_2\right]  + \beta_2 \left\| \nabla^{\top} \left( v^{t+1} + \frac{y^t}{\beta_2} \right) + \Delta z^t\right\|_2 \\
\leq& \sum_{j=1}^N \left[\beta_1 \|\omega^{t+1}\|_{\infty} \left(\|u_j^{t+1}\|_2 + \frac{1}{\beta_1}\|\Lambda_j^t\|_2 + \|\omega^{t+1}\|_{\infty}\|z^{t+1}\|_2 \right)\right]\\ &+ \beta_2 \left[\|\nabla\|\left( \|v^{t+1}\|_2 + \frac{1}{\beta_2}\|y^t\|_2\right) + \|\Delta\|\|z^t\|_2 \right].
\end{align*}
So, there exists a constant $C > 0$ such that $\|\nabla_{\omega} \mathcal{L} (\omega^t)\|_2^2 \leq C$ and $\|\nabla_{z} \mathcal{L} (z^t) \|_2^2 \leq C$ for each iteration $t \in \mathbb{N}$. From {\eqref{eq:expect_smooth_omega}}-{\eqref{eq:expect_smooth_z}}, we have
\begin{align*}
\mathbb{E}_t \left[ \|\tilde{\nabla}_{\omega} \mathcal{L}(\omega^t)\|_2^2 \;\middle |\; u^{t+1}\right]- \| \nabla_{\omega} \mathcal{L}(\omega^t)\|_2^2  &\leq A_1 \mathcal{L}(\omega^t) + A_2  \| \nabla_{\omega} \mathcal{L}(\omega^t)\|_2^2 +A_3 \\
& \leq A_1B + A_2C +A_3 = \sigma^2
\end{align*}
and
\begin{align*}
\mathbb{E}_t \left[ \|\tilde{\nabla}_{z} \mathcal{L}(z^t)\|_2^2\;\middle |\; u^{t+1}, \omega^{t+1}, v^{t+1}\right] -\| \nabla_{z} \mathcal{L}(z^t)\|_2^2 &\leq A_1 \mathcal{L}(z^t) + A_2 \| \nabla_{z} \mathcal{L}(z^t)\|_2^2 + A_3\\ & \leq A_1B + A_2C +A_3 = \sigma^2,
\end{align*}
where $\sigma \coloneqq \sqrt{A_1B + A_2C + A_3}$. 
\end{proof}
\subsection{Proof of Lemma \ref{lemma:sufficient_decrease}}
\begin{proof} If $\{(\omega^t, z^t)\}_{t=1}^{\infty}$ is bounded, then there exists a constant $C>0$ such that $\|\omega^t\|_{\infty}, \|z^t\|_{\infty} \leq C$ for all $t \in \mathbb{N}$. 
We establish that $\mathcal{L}(\omega) \coloneqq \mathcal{L}(u^{t+1}, \omega, v^t, z^t, \Lambda^t, y^t)$ has a Lipschitz continuous gradient with respect to $\omega$. At each iteration $t$, we have
\begin{align*}
    \nabla_{\omega} \mathcal{L} (\omega) = -\sum_{j=1}^N  \beta_1 (S_jz^t)^* \circ \left[\mathcal{F}^{-1}\left(u_j^{t+1} + \frac{\Lambda_j^t}{\beta_1} \right) - \omega \circ S_j z^t \right],
\end{align*} so for any $\omega_1, \omega_2 \in \mathbb{C}^{m^2}$, we estimate
\begin{align*}
    &\left \|\nabla_{\omega} \mathcal{L}(\omega_2) - \nabla_{\omega} \mathcal{L}(\omega_1) \right\|_2 = \left \|\sum_{j=1}^N \beta_1 (S_j z^t)^* \circ (\omega_2 - \omega_1) \circ (S_j z^t)  \right\|_2\\ &\leq \left[ \sum_{j=1}^N \beta_1 \left\|(S_j z^t)^*\circ (S_j z^t) \right\|_{\infty}\right] \|\omega_2 -\omega_1\|_2 \leq
    \beta_1NC^2 \|\omega_2- \omega_1\|_2. 
\end{align*}
Hence, we observe that $\mathcal{L}(\omega)$ has a Lipschitz continuous gradient with Lipschitz constant $L_{\omega} \coloneqq \beta_1NC^2$. By the descent property \cite[Definition 1]{chang2019blind}, at iteration $t$ we have
\begin{align*}
    \mathcal{L}(\omega^{t+1}) - \mathcal{L}(\omega^t) &\leq \mathbb{R}(\langle  \nabla_{\omega} \mathcal{L}(\omega^t), \omega^{t+1} - \omega^t\rangle) + \frac{L_{\omega}}{2} \|\omega^{t+1} - \omega^t \|_2^2\\
    &= -\delta_{\omega}^t\mathbb{R}(\langle  \nabla_{\omega} \mathcal{L}(\omega^t),\tilde{\nabla}_{\omega}\mathcal{L}(\omega^t) \rangle) + \frac{L_{\omega} (\delta^t_{\omega})^2}{2} \|\tilde{\nabla}_{\omega}\mathcal{L}(\omega^t)  \|_2^2,
\end{align*}
where the last equality is due to {\eqref{eq:sgd_omega}}. Because we assume that $\{(\omega^t, z^t)\}_{t=1}^{\infty}$ is bounded and  $\sum_{t=1}^{\infty} \|\Omega^{t+1} - \Omega^t\|_2^2 < \infty$, \eqref{eq:omega_var_bound}-\eqref{eq:z_var_bound} hold by Lemma \ref{lemma:bounded_sequence}. Taking the expectation conditioned on the first $t$ iterations and $u^{t+1}$, we obtain
\begin{align*}
&\mathbb{E}_t\left[\mathcal{L}(\omega^{t+1})\;\middle |\; u^{t+1} \right] - \mathcal{L}(\omega^t)\\
    &= -\delta_{\omega}^t\mathbb{R}\left(\mathbb{E}_t\left[ \left \langle  \nabla_{\omega} \mathcal{L}(\omega^t),\tilde{\nabla}_{\omega}\mathcal{L}(\omega^t) \right \rangle \;\middle |\; u^{t+1}\right]\right) + \frac{L_{\omega} (\delta^t_{\omega})^2}{2} \mathbb{E}_t\left[\|\tilde{\nabla}_{\omega}\mathcal{L}(\omega^t)  \|_2^2\;\middle |\; u^{t+1}\right] \\
    &= -\delta_{\omega}^t \mathbb{R} \left( \left \langle  \nabla_{\omega} \mathcal{L}(\omega^t),\mathbb{E}_t\left[\tilde{\nabla}_{\omega}\mathcal{L}(\omega^t)\;\middle |\; u^{t+1} \right] \right \rangle \right) + \frac{L_{\omega} (\delta^t_{\omega})^2}{2} \mathbb{E}_t\left[\|\tilde{\nabla}_{\omega}\mathcal{L}(\omega^t)  \|_2^2 \;\middle |\; u^{t+1}\right]\\&\leq \frac{-2 \delta_{\omega}^t + L_{\omega} (\delta_{\omega}^t)^2}{2 }\|\nabla_{\omega} \mathcal{L}(\omega^t) \|_2^2 + \frac{L_{\omega} (\delta_{\omega}^t)^2\sigma^2}{2}.
\end{align*}
 The last inequality is due to Assumption \ref{assume:gradient_bound}(a) and {\eqref{eq:omega_var_bound}}. Taking total expectation gives us \eqref{eq:lemma4_5_omega_ineq}.

Similarly, we can estimate \eqref{eq:lemma4_5_z_ineq} because we can compute that $\mathcal{L}(z)\coloneqq \mathcal{L}(u^{t+1}, \omega^{t+1}, v^{t+1}, z, \Lambda^t, y^t)$ has a Lipschitz continuous gradient with Lipschitz constant $L_z \coloneqq \beta_1 N C^2 + \beta_2 \|\Delta\|$ and follow the same steps as above by taking expectation conditioned on the first $t$ iterations and $(u^{t+1}, \omega^{t+1}, v^{t+1})$.
\end{proof}
\subsection{Proof of Proposition \ref{prop:finite_length}}
\begin{proof}
By Lemma {\ref{lemma:bounded_sequence}},  $\{(Z^t, \Omega^t)\}_{t=1}^{\infty}$ is bounded.
We see that
\begin{align*}
    \mathcal{L}(Z, \Omega) = &\sum_{j=1}^N \left[ \mathcal{B}(|u_j|^2, d_j) + \frac{\beta_1}{2} \left\|u_j - \mathcal{F}(\omega \circ S_j z ) + \frac{\Lambda_j}{\beta_1} \right \|_2^2 - \frac{1}{2 \beta_1} \|\Lambda_j\|_2^2 \right]\\
    &+\lambda(\|v\|_1 - \alpha \|v\|_{2,1}) + \frac{\beta_2}{2} \left\| v - \nabla z + \frac{y}{\beta_2} \right\|_2^2 - \frac{1}{2 \beta_2}\| y\|_2^2 \\
    &\geq \sum_{j=1}^N \left[\mathcal{B}(|u_j|^2, d_j)- \frac{1}{2 \beta_1} \|\Lambda_j\|_2^2 \right]- \frac{1}{2 \beta_2}\| y\|_2^2.
\end{align*}
Because $\mathcal{B}(\cdot, \cdot)$ is bounded below according to \eqref{eq:phase_fidelity} and $\{(Z^t, \Omega^t)\}_{t=1}^{\infty}$ is bounded, $\{\mathcal{L}(Z^t, \Omega^t)\}_{t=1}^{\infty}$ is bounded below by some constant $\mathcal{L}_{\inf}$. At iteration $t$, we denote the following for the sake of brevity:
\begin{align*}
    \mathcal{L}(u) &\coloneqq \mathcal{L}(u, \omega^t, v^t, z^t, \Lambda^t, y^t), \\
    \mathcal{L}(v) &\coloneqq \mathcal{L}(u^{t+1}, \omega^{t+1}, v, z^t, \Lambda^t, y^t), \\
    \mathcal{L}(\Lambda) &\coloneqq \mathcal{L}(u^{t+1}, \omega^{t+1}, v^{t+1}, z^{t+1}, \Lambda, y^t),\\
    \mathcal{L}(y) &\coloneqq  \mathcal{L}(u^{t+1}, \omega^{t+1}, v^{t+1}, z^{t+1}, \Lambda^{t+1}, y).
\end{align*}By \eqref{eq:u_min}, we have
\begin{align*}
    &\mathcal{B}(|u_j^{t+1}|^2, d_j) +\mathbb{R}\left(\langle \Lambda_j^t, u_j^{t+1}- \mathcal{F}(P_j^t z^t) \rangle \right) + \frac{\beta_1}{2}  \|u_j^{t+1} - \mathcal{F}(P_j^t
    z^t) \|_2^2\\ &\leq \mathcal{B}(|u_j^{t}|^2, d_j) +\mathbb{R}\left(\langle \Lambda_j^t, u_j^{t}- \mathcal{F}(P_j^t z^t) \rangle \right) + \frac{\beta_1}{2}  \|u_j^{t} - \mathcal{F}(P_j^t
    z^t) \|_2^2
\end{align*}
for each $j \in n^t$, and by \eqref{eq:u_closed_form}, we have
\begin{align*}
 &\mathcal{B}(|u_j^{t+1}|^2, d_j) +\mathbb{R}\left(\langle \Lambda_j^t, u_j^{t+1}- \mathcal{F}(P_j^t z^t) \rangle \right) + \frac{\beta_1}{2}  \|u_j^{t+1} - \mathcal{F}(P_j^t
    z^t) \|_2^2\\ &= \mathcal{B}(|u_j^{t}|^2, d_j) +\mathbb{R}\left(\langle \Lambda_j^t, u_j^{t}- \mathcal{F}(P_j^t z^t) \rangle \right) + \frac{\beta_1}{2}  \|u_j^{t} - \mathcal{F}(P_j^t
    z^t) \|_2^2
\end{align*}
for each $j \not \in n^t$. Summing over all $j=1, \ldots, N$ and adding the term $\lambda \left( \|v^t\|_1 - \alpha \|v^t\|_{2,1} \right )+ \mathbb{R} \left(\langle y^t, v^t - \nabla z^t \rangle \right) + \frac{\beta_2}{2} \|v^t - \nabla z^t \|_2^2$ to both sides of the inequality, we obtain  $\mathcal{L}(u^{t+1}) \leq \mathcal{L}(u^{t})$. By \eqref{eq:v_sub}, we have $\mathcal{L}(v^{t+1}) \leq \mathcal{L}(v^{t})$, so taking expectation, we obtain
\begin{align}\label{eq:expect_u_update}
\mathbb{E}[\mathcal{L}(\omega^{t})]=\mathbb{E}[\mathcal{L}(u^{t+1})] &\leq \mathbb{E}[\mathcal{L}(u^{t})],\\
\mathbb{E}[\mathcal{L}(z^{t})] =\mathbb{E}[\mathcal{L}(v^{t+1})] &\leq \mathbb{E}[\mathcal{L}(v^{t})] = \mathbb{E}[\mathcal{L}(\omega^{t+1})]. 
\end{align}
In addition, we have
\begin{align*}
  \mathcal{L}(\Lambda^{t+1})  - \mathcal{L}(\Lambda^{t}) =& \sum_{j=1}^N \mathbb{R}( \langle \Lambda_j^{t+1} - \Lambda_j^t, u_j^{t+1} - \mathcal{F}(\omega^{t+1} \circ S_j z^{t+1}) \rangle) \\
  =& \sum_{j\in n^t} \mathbb{R}( \langle \Lambda_j^{t+1} - \Lambda_j^t, u_j^{t+1} - \mathcal{F}(\omega^{t+1} \circ S_j z^{t+1}) \rangle)\\ &+\sum_{j \not \in n^t} \mathbb{R}( \langle \Lambda_j^{t+1} - \Lambda_j^t, u_j^{t+1} - \mathcal{F}(\omega^{t+1} \circ S_j z^{t+1}) \rangle) \\
  =&  \frac{1}{\beta_1} \sum_{j\in n^t}\left \|\Lambda_j^{t+1} -\Lambda_j^{t}  \right\|_2^2 =\frac{1}{\beta_1}\|\Lambda^{t+1}-\Lambda^t\|_2^2,
\end{align*}
due to \eqref{eq:lambda_update2}. Taking expectation gives
\begin{align}\label{eq:expect_lambda_update}
    \mathbb{E}[  \mathcal{L}(\Lambda^{t+1})] - \mathbb{E}[  \mathcal{L}(\Lambda^{t})] =\frac{1}{\beta_1} \mathbb{E}\left[\left \|\Lambda^{t+1} -\Lambda^{t}  \right\|_2^2\right].
\end{align}
Similarly, we obtain
\begin{align}\label{eq:expect_y_update}
    \mathbb{E}[  \mathcal{L}(y^{t+1})] - \mathbb{E}[  \mathcal{L}(y^{t})] &= \frac{1}{\beta_2} \mathbb{E}[\|y^{t+1} -y^t\|_2^2].
\end{align}
Summing up \eqref{eq:lemma4_5_omega_ineq}-\eqref{eq:lemma4_5_z_ineq}, \eqref{eq:expect_u_update}-\eqref{eq:expect_y_update}, we have
\begin{gather}
\begin{aligned}
    \mathbb{E} [\mathcal{L}(Z^{t+1}, \Omega^{t+1})] &-  \mathbb{E} [\mathcal{L}(Z^{t}, \Omega^t)] \leq \frac{1}{\beta_1}\mathbb{E}\left[\left \|\Lambda^{t+1} -\Lambda^{t}  \right\|_2^2\right] + \frac{1}{\beta_2} \mathbb{E}[\|y^{t+1} -y^t\|_2^2]\\
    &+\frac{-2\delta_{\omega}^t +L_{\omega} (\delta_{\omega}^t)^2}{2}  \mathbb{E}\left[\|\nabla_{\omega}\mathcal{L}(\omega^t)\|_2^2\right] + \frac{L_{\omega}(\delta_{\omega}^t)^2 \sigma^2}{2}\\
    &+\frac{-2\delta_{z}^t +L_{z} (\delta_{z}^t)^2}{2} \mathbb{E}\left[\|\nabla_{z} \mathcal{L}(z^t)\|_2^2\right] +  \frac{L_{z}(\delta_{z}^t)^2 \sigma^2}{2}.
\end{aligned}
\end{gather}
Summing up $t=1, \ldots, T$, we obtain
\begin{align*}
    \mathcal{L}_{\inf} - \mathcal{L}(Z^0, \Omega^0) &\leq \mathbb{E}[\mathcal{L}(Z^{T+1}, \Omega^{T+1})] - \mathcal{L}(Z^0, \Omega^0) \\
    &\leq \sum_{t=1}^T \Bigg( C_1\| \Omega^{t+1} - \Omega^t\|_2^2+ 
    \frac{-2\delta_{\omega}^t +L_{\omega} (\delta_{\omega}^t)^2}{2}  \mathbb{E}\left[\|\nabla_{\omega}\mathcal{L}(\omega^t)\|_2^2\right]\\&+\frac{-2\delta_{z}^t +L_{z} (\delta_{z}^t)^2}{2} \mathbb{E}\left[\|\nabla_{z} \mathcal{L}(z^t)\|_2^2\right]+ \frac{ \sigma^2(L_{\omega}(\delta_{\omega}^t)^2+L_{z}(\delta_{z}^t)^2)}{2}\Bigg),
\end{align*}
where $C_1 = \max \{ \frac{1}{\beta_1}, \frac{1}{\beta_2}\}$. Note that $\mathbb{E}[\mathcal{L}(Z^0, \Omega^0)] = \mathcal{L}(Z^0, \Omega^0)$ since $(Z^0, \Omega^0)$ serves as initialization. By \eqref{eq:step_size_condition}, we have $\delta_{\omega}^t < \frac{2}{L_{\omega}}$ and $\delta_{z}^t < \frac{2}{L_z}$ for all $t$, which follows that $\frac{-2\delta_{\omega}^t +L_{\omega} (\delta_{\omega}^t)^2}{2}, \frac{-2\delta_{z}^t +L_{z} (\delta_{z}^t)^2}{2} <0$.
Rearranging the inequality and letting $T \rightarrow \infty$ give us
\begin{align*}
    0\leq&\sum_{t=1}^{\infty}  \left( \frac{2\delta_{\omega}^t -L_{\omega} (\delta_{\omega}^t)^2}{2} \mathbb{E} \left[ \left\| \nabla_{\omega} \mathcal{L}(\omega^t) \right\|_2^2 \right]+ \frac{2\delta_{z}^t -L_{z} (\delta_{z}^t)^2}{2} \mathbb{E} \left[ \left\| \nabla_{z} \mathcal{L}(z^t) \right\|_2^2 \right] \right) \leq \\ &\mathcal{L}(Z^0, \Omega^0)-\mathcal{L}_{\inf} + \sum_{t=1}^{\infty} \Bigg(C_1\| \Omega^{t+1} - \Omega^t\|_2^2+\frac{ \sigma^2(L_{\omega}(\delta_{\omega}^t)^2+L_{z}(\delta_{z}^t)^2)}{2}\Bigg).
\end{align*}
By the summability assumption, the right-hand side is bounded. As a result, there exists a positive constant $C_2$ such that 
\begin{align*}
    \sum_{t=1}^{\infty} \left( \frac{2\delta_{\omega}^t -L_{\omega} (\delta_{\omega}^t)^2}{2} \mathbb{E} \left[ \left\| \nabla_{\omega} \mathcal{L}(\omega^t) \right\|_2^2 \right] \right) \leq C_2. 
\end{align*}
By \eqref{eq:step_size_condition}, we have
\begin{align*}
    \sum_{t=1}^{\infty}  \frac{2\delta_{\omega}^t -L_{\omega} (\delta_{\omega}^t)^2}{2}  = \infty,
\end{align*}
which implies that for any $\tau > 0$,
\begin{align*}
 \sum_{t=\tau}^{\infty}  \frac{2\delta_{\omega}^t -L_{\omega} (\delta_{\omega}^t)^2}{2}  = \infty.
\end{align*}
For any $T > \tau > 0$, we have
\begin{align*}
    \min_{t =\tau, \ldots, T-1} \mathbb{E} \left[ \left\| \nabla_{\omega} \mathcal{L}(\omega^{t}) \right\|_2^2 \right] & \leq \left[ \sum_{t=\tau}^{T-1} \left( \frac{2\delta_{\omega}^t -L_{\omega} (\delta_{\omega}^t)^2}{2} \mathbb{E} \left[ \left\| \nabla_{\omega} \mathcal{L}(\omega^t) \right\|_2^2 \right] \right) \right] \cdot \left[ \sum_{t=\tau}^{T-1} \left( \frac{2\delta_{\omega}^t -L_{\omega} (\delta_{\omega}^t)^2}{2} \right)\right]^{-1} \\
    &\leq C_2  \left[ \sum_{t=\tau}^{T-1} \left( \frac{2\delta_{\omega}^t -L_{\omega} (\delta_{\omega}^t)^2}{2} \right)\right]^{-1}.
\end{align*}
As $T \rightarrow \infty$, we obtain $\inf_{\tau \geq t} \mathbb{E} \left[ \left\| \nabla_{\omega} \mathcal{L}(\omega^{\tau}) \right\|_2^2 \right]  =0$
for any $t > 0$. Therefore, taking the limit as $t \rightarrow \infty$, we obtain
\begin{align*}
    \liminf_{t \rightarrow \infty} \mathbb{E} \left[ \left\| \nabla_{\omega} \mathcal{L}(\omega^t) \right\|_2^2 \right] = \lim_{t \rightarrow \infty} \inf_{\tau \geq t} \mathbb{E} \left[ \left\| \nabla_{\omega} \mathcal{L}(\omega^{\tau}) \right\|_2^2 \right]  =0.
\end{align*}
\eqref{eq:liminf2} is proven similarly by following the same steps.
\end{proof}
\subsection{Proof of Theorem \ref{thm:kkt}}
\begin{proof}
By Lemma {\ref{lemma:bounded_sequence}}, $\{(Z^t, \Omega^t)\}_{t=1}^{\infty}$ and  $\{(\nabla_{\omega} \mathcal{L}(\omega^t), \nabla_z \mathcal{L}(z^t))\}_{t=1}^{\infty}$ are bounded and \eqref{eq:omega_var_bound}-\eqref{eq:z_var_bound} hold. So, there exists $C_1 > 0$ such that $\|\nabla_{\omega}\mathcal{L}(\omega^t)\|_2, \| \nabla_z \mathcal{L}(z^t)\|_2 < C_1$ for all $t \in \mathbb{N}$. Recall that the SGD steps for $(\omega,z)$ are
\begin{align*}
    \omega^{t+1} = \omega^t - \delta_{\omega}^t \tilde{\nabla}\mathcal{L}(\omega^t), \;
    z^{t+1} = z^t - \delta_z^t \tilde{\nabla}\mathcal{L}(z^t).
\end{align*}
By \eqref{eq:omega_var_bound}-\eqref{eq:z_var_bound},
\begin{align*}
    \mathbb{E}_t \left[ \left\|\omega^{t+1} - \omega^{t}\right\|_2^2 \; \middle | \; u^{t+1} \right] &=(\delta_{\omega}^t)^2 \mathbb{E}_t \left [\|\tilde{\nabla}_{\omega} \mathcal{L}(\omega^t)\|_2^2\; \middle | \;u^{t+1} \right]\\ &  \leq(\delta_{\omega}^t)^2 \left(\sigma^2 + \|\nabla_{\omega} \mathcal{L}(\omega^t)\|_2^2\right)\leq (\delta_{\omega}^t)^2 \left( \sigma^2 + C_1^2 \right),\\
     \mathbb{E}_t \left [ \left\|z^{t+1} - z^{t}\right\|_2^2\; \middle | \;u^{t+1}, \omega^{t+1}, v^{t+1} \right] &= (\delta_{z}^t)^2 \mathbb{E}_t\left [\|\tilde{\nabla}_{z} \mathcal{L}(z^t)\|_2^2\; \middle | \;u^{t+1}, \omega^{t+1}, v^{t+1} \right] \\ &\leq(\delta_{z}^t)^2 \left( \sigma^2 + \|\nabla_z \mathcal{L}(z^t)\|_2^2 \right) \leq (\delta_{z}^t)^2 \left( \sigma^2 + C_1^2 \right).
\end{align*}
Applying total expectation, summing over all $t$, and using \eqref{eq:step_size_condition} establish
\begin{align}
  \sum_{t=1}^{\infty} \mathbb{E} [ \left\|\omega^{t+1} - \omega^{t}\right\|_2^2] \leq  \left( \sigma^2 + C_1^2 \right)\sum_{t=1}^{\infty}(\delta_{\omega}^t)^2  < \infty, \\
    \sum_{t=1}^{\infty} \mathbb{E} [ \left\|z^{t+1} - z^{t}\right\|_2^2]  \leq \left( \sigma^2 + C_1^2\right)  \sum_{t=1}^{\infty}(\delta_{z}^t)^2 < \infty.
\end{align}
By Lemma \ref{lemma:sum_converge}, we have 
\begin{align}
\lim_{t \rightarrow \infty}  \left\|\omega^{t+1} - \omega^{t}\right\|_2^2 &=0 \; \text{a.s.} \implies \lim_{t \rightarrow \infty} \left(\omega^{t+1} - \omega^t \right)= 0\; \text{a.s.}, \label{eq:omega_diff_converge}\\
\lim_{t \rightarrow \infty} \left\|z^{t+1} - z^{t}\right\|_2^2 &=0 \; \text{a.s.} \implies \lim_{t \rightarrow \infty} \left( z^{t+1} - z^t \right)= 0 \; \text{a.s.} \label{eq:z_diff_converge}
\end{align}
For the sake of brevity, we will omit ``a.s." in the rest of the proof.
Earlier, in the proof of Lemma \ref{lemma:bounded_sequence}, we obtained
\begin{align}
    \lim_{t \rightarrow \infty} v^t - \nabla z^t &= 0.\label{eq:v_equal_nabla_z2}
\end{align}
By \eqref{eq:liminf1}-\eqref{eq:liminf2} in Proposition {\ref{prop:finite_length}}, there exists a subsequence $\{Z^{t_k}\}_{k=1}^{\infty}$ such that
\begin{align}
    \lim_{k \rightarrow \infty} \mathbb{E}[\|\nabla_{\omega} \mathcal{L}(\omega^{t_k})\|_2^2] &= 0 \label{eq:liminf_omega},\\
    \lim_{k \rightarrow \infty} \mathbb{E}[\|\nabla_{z} \mathcal{L}(z^{t_k})\|_2^2] &= 0 \label{eq:lim2}.
\end{align}
Because $\{(Z^{t_k}, \Omega^{t_k})\}_{k=1}^{\infty}$ is bounded, there exists a subsequence $\{(Z^{t_{k_{\ell}}}, \Omega^{t_{k_{\ell}}})\}_{\ell=1}^{\infty}$ such that $\displaystyle \lim_{\ell \rightarrow \infty}(Z^{t_{k_{\ell}}}, \Omega^{t_{k_{\ell}}})  = (Z^{\star}, \Omega^{\star})$. In addition, \eqref{eq:omega_diff_converge}-\eqref{eq:lim2} hold for this subsequence.


Now we show that $(Z^{\star},\Omega^{\star})$ is a KKT point. From \eqref{eq:v_equal_nabla_z2}, we have
\begin{align}\label{eq:v_star_converge}
    v^{\star} &= \lim_{\ell \rightarrow \infty} v^{t_{k_{\ell}}} = \lim_{\ell \rightarrow \infty} \nabla z^{t_{k_{\ell}}} = \nabla z^{\star}. 
\end{align}
Since $\sum_{t=1}^{\infty} \|\Lambda^{t+1} - \Lambda^t\|_2^2 < \infty$ by the assumption, it follows that $\displaystyle \lim_{t \rightarrow \infty} \left( \Lambda^{t+1} - \Lambda^t \right)= 0$. For each $j = 1, \ldots, N$, there exists $L_1 > 0$ for $\epsilon > 0$ such that $\ell \geq L_1$ implies the following:
\begin{align}
\|u^{\star}- u^{t_{k_{\ell}}}\|_2 &< \epsilon, \label{eq:u_epsilon}\\
\|\Lambda^{t_{k_{\ell}}} - \Lambda^{t_{k_{\ell}}-1}\|_2 &< \beta_1\epsilon, \label{eq:lambda_epsilon}\\
\|\omega^{\star} \circ S_j z^{\star} - \omega^{t_{k_{\ell}}} \circ S_j z^{t_{k_{\ell}}} \|_2 &< \epsilon. \label{eq:f_epsilon}
\end{align}
    Note that \eqref{eq:u_epsilon} is due to $\displaystyle \lim_{\ell \rightarrow \infty} u^{t_{k_{\ell}}}  = u^{\star}$; \eqref{eq:lambda_epsilon} is due to  $\displaystyle \lim_{t \rightarrow \infty}  \left( \Lambda^{t+1} - \Lambda^{t} \right) = 0$; and \eqref{eq:f_epsilon} is due to $\displaystyle \lim_{\ell \rightarrow \infty} \omega^{t_{k_{\ell}}} \circ S_j z^{t_{k_{\ell}}} = \omega^{\star}\circ S_j z^{\star}$. Then, by Lemma \ref{lemma:sampling}, there exists $\ell' \geq L_1$ such that $j \in n^{t_{k_{\ell'}}-1}$, so we have
\begin{align*}
\|u_j^{\star} - \mathcal{F}(\omega^{\star} \circ S_j z^{\star})\|_2 &\leq \|u_j^{\star} - u_j^{t_{k_{\ell'}}} \|_2 + \|u_j^{t_{k_{\ell'}}} - \mathcal{F}(\omega^{\star} \circ S_j z^{\star}) \|_2 \\
&= \|u_j^{\star} - u_j^{t_{k_{\ell'}}} \|_2 + \left \|\frac{1}{\beta_1} \left(\Lambda_j^{t_{k_{\ell'}}} -\Lambda_j^{t_{k_{\ell'}}-1} \right) + \mathcal{F}(\omega^{t_{k_{\ell'}}} \circ S_j z^{t_{k_{\ell'}}}) - \mathcal{F}(\omega^{\star} \circ S_j z^{\star}) \right\|_2 \\
& \leq \|u_j^{\star} - u_j^{t_{k_{\ell'}}} \|_2  + \frac{1}{\beta_1} \left\| \Lambda_j^{t_{k_{\ell'}}} -\Lambda_j^{t_{k_{\ell'}}-1}\right\|_2 + \left\|\omega^{t_{k_{\ell'}}} \circ S_j z^{t_{k_{\ell'}}} - \omega^{\star} \circ S_j z^{\star}  \right\|_2 \\
&< 3 \epsilon
\end{align*}
after applying \eqref{eq:u_update3} and \eqref{eq:u_epsilon}-\eqref{eq:f_epsilon}. Because $\epsilon > 0$ is chosen arbitrarily, it follows that
\begin{align}\label{eq:u_star_equal_f}
    u_j^{\star} =  \mathcal{F}(\omega^{\star} \circ S_j z^{\star}),\; \forall j=1,\ldots, N.
\end{align}

Before proving the rest of the conditions, we need to show that $\displaystyle \lim_{\ell \rightarrow \infty} v^{t_{k_{\ell}}+1} = v^{\star}$ and $\displaystyle \lim_{\ell \rightarrow \infty} u_j^{t_{k_{\ell}}+1} = u_j^{\star}$ for each $j = 1, \ldots, N$. By \eqref{eq:z_diff_converge}, \eqref{eq:v_equal_nabla_z2}, and \eqref{eq:v_star_converge}, we have
\begin{gather}
\begin{aligned}\label{eq:v_converge}
    \lim_{\ell \rightarrow \infty} v^{t_{k_{\ell}}+1} &=    \lim_{\ell \rightarrow \infty} \left( v^{t_{k_{\ell}}+1}- \nabla z^{t_{k_{\ell}}+1} + \nabla z^{t_{k_{\ell}}+1}  -  \nabla z^{t_{k_{\ell}}} +  \nabla z^{t_{k_{\ell}}} \right)   \\
    &=\lim_{\ell \rightarrow \infty} \left( v^{t_{k_{\ell}}+1}- \nabla z^{t_{k_{\ell}}+1} \right) +\lim_{\ell \rightarrow \infty} \left(  \nabla z^{t_{k_{\ell}}+1}  -  \nabla z^{t_{k_{\ell}}} \right) + \lim_{\ell \rightarrow \infty} \nabla z^{t_{k_{\ell}}} \\
&=   \lim_{\ell \rightarrow \infty} v^{t_{k_{\ell}}}=v^{\star}.
\end{aligned}
\end{gather}
Since  $\{(\omega^t, z^t)\}_{t=1}^{\infty}$ is bounded, then $\|\omega^t\|_{\infty}, \|z^t\|_{\infty} \leq C_2$ for all $t \in \mathbb{N}$ for some constant $C_2>0$. For $\epsilon > 0$, there exists $L_2 > 0$ such that $ \ell \geq L_2$ implies the following:
\begin{align}
\|\omega^{t_{k_{\ell}}+1} - \omega^{t_{k_{\ell}}}\|_2 &< \frac{1}{C_2} \epsilon \label{eq:successive_omega} \\
\|z^{t_{k_{\ell}}+1} - z^{t_{k_{\ell}}}\|_2 &< \frac{1}{ C_2} \epsilon \label{eq:successive_z} \\
\|\Lambda^{t_{k_{\ell}}+1} - \Lambda^{t_{k_{\ell}}}\|_2 &< \beta_1\epsilon. \label{eq:sucessive_Lambda}
\end{align}
Note that \eqref{eq:successive_omega}-\eqref{eq:successive_z} are due to \eqref{eq:omega_diff_converge}-\eqref{eq:z_diff_converge} and \eqref{eq:sucessive_Lambda} is due to $\displaystyle \lim_{t \rightarrow \infty} \left( \Lambda^{t+1} - \Lambda^t \right) = 0$.  By definition of $S_j$, we have $\|S_j z\|_2 \leq \|z\|_2$ for all $z \in \mathbb{C}^{n^2}$.   For each $j=1, \ldots, N$, when $\ell \geq \max\{L_1,L_2\}$, we establish that if $j \not \in n^{t_{k_{\ell}}}$, then
\begin{align}
\|u_j^{t_{k_{\ell}}+1} - u_j^{\star}\|_2 = \|u_j^{t_{k_{\ell}}}-u_j^{\star}\|_2 < \epsilon < 4\epsilon
\end{align}
 by \eqref{eq:u_update3} and \eqref{eq:u_epsilon}, and if $j \in n^{t_{k_{\ell}}}$, then
 \begin{align*}
 &\|u_j^{t_{k_{\ell}}+1} - u_j^{\star}\|_2 = \left\|   \frac{1}{\beta_1} \left( \Lambda_j^{t_{k_{\ell}}+1}- \Lambda_j^{t_{k_{\ell}}} \right) + \mathcal{F}(\omega^{t_{k_{\ell}}+1} \circ S_j z^{t_{k_{\ell}}+1})-u_j^{\star}\right\|_2\\ &\leq \frac{1}{\beta_1} \|\Lambda_j^{t_{k_{\ell}}+1} - \Lambda_j^{t_{k_{\ell}}}\|_2 + \| \mathcal{F} \left(\omega^{t_{k_{\ell}}+1} \circ S_j z^{t_{k_{\ell}}+1} - \omega^{t_{k_{\ell}}} \circ S_j z^{t_{k_{\ell}}}\right) \|_2 + \left\| \mathcal{F} \left(\omega^{t_{k_{\ell}}} \circ S_j z^{t_{k_{\ell}}}\right) -  \mathcal{F}(\omega^{\star} \circ S_j z^{\star})\right\|_2 \\
 &< \epsilon + \left\| \omega^{t_{k_{\ell}}+1} \circ S_j z^{t_{k_{\ell}}+1} - \omega^{t_{k_{\ell}}} \circ S_j z^{t_{k_{\ell}}}\right\|_2+\left\| \omega^{t_{k_{\ell}}} \circ S_j z^{t_{k_{\ell}}} - \omega^{\star} \circ S_j z^{\star}\right\|_2\\ &\leq 2 \epsilon +  \left\| \omega^{t_{k_{\ell}}+1} \circ S_j z^{t_{k_{\ell}}+1} - \omega^{t_{k_{\ell}}+1} \circ S_j z^{t_{k_{\ell}}}\right\|_2 + \left\|\omega^{t_{k_{\ell}}+1} \circ S_j z^{t_{k_{\ell}}} -\omega^{t_{k_{\ell}}} \circ S_j z^{t_{k_{\ell}}}\right\|_2 \\
 &\leq 2 \epsilon + \left\| \omega^{t_{k_{\ell}}+1} \right\|_{\infty}\left\|z^{t_{k_{\ell}}+1} - z^{t_{k_{\ell}}}  \right\|_2 + \|z^{t_{k_{\ell}}}\|_{\infty} \| \omega^{t_{k_{\ell}}+1} - \omega^{t_{k_{\ell}}} \|_2 \\
 &<4\epsilon 
 \end{align*}
 by \eqref{eq:u_update3}, \eqref{eq:f_epsilon}, \eqref{eq:u_star_equal_f}, and \eqref{eq:successive_omega}-\eqref{eq:sucessive_Lambda}. Because $\epsilon>0$ is chosen arbitrarily, we have
 \begin{align}\label{eq:u_t_k_l_equal_u_star}
     \displaystyle \lim_{\ell \rightarrow \infty} u_j^{t_{k_{\ell}}+1} = u_j^{\star}.
 \end{align}

 Next we prove \eqref{eq:u_kkt}. At iteration $t_{k_{\ell}}$ for each $j=1, \ldots, N$, the first-order optimality condition of \eqref{eq:u_sub} is
\begin{align*}
     - \Lambda_j^{t_{k_{\ell}}} - \beta_1 \left(u_j^{t_{k_{\ell}}+1} - \mathcal{F}(\omega^{t_{k_{\ell}}} \circ S_j z^{t_{k_{\ell}}}) \right) &\in  \begin{cases}\partial \left(\displaystyle \frac{1}{2}\left\| |u_j^{t_{k_{\ell}}+1}| - \sqrt{d_j} \right\|^2_2 \right), &\text{ if AGM},\\
      \partial \left(\displaystyle \frac{1}{2} \left \langle |u_j^{t_{k_{\ell}}+1}|^2 - d_j \log |u_j^{t_{k_{\ell}}+1}|^2, \mathbf{1}\right\rangle \right), & \text{ if IPM}
     \end{cases}\\ &= \begin{cases} \partial \left|u_j^{t_{k_{\ell}}+1} \right| \circ \left( |u_j^{t_{k_{\ell}}+1}| - \sqrt{d_j} \right), & \text{ if AGM,}\\
      \partial \left|u_j^{t_{k_{\ell}}+1} \right| \circ \left( |u_j^{t_{k_{\ell}}+1}| - \displaystyle\frac{d_j}{|u_j^{t_{k_{\ell}}+1}|} \right), & \text{ if IPM}. 
     \end{cases}
\end{align*}
For the IPM case, because $\displaystyle\lim_{x \rightarrow 0^+} x - d \log x = + \infty$ for any $d > 0$, it follows that $(u_j^{\star})_i \neq 0$ for all $i$, so we do  need to worry about $u_j^{t_{k_{\ell}}+1}= 0$ for any iteration $t_{k_{\ell}}+1$.
By \eqref{eq:u_star_equal_f} and \eqref{eq:u_t_k_l_equal_u_star}, we see that
\begin{align*}
   \lim_{\ell \rightarrow \infty}- \Lambda_j^{t_{k_{\ell}}} - \beta_1 \left(u_j^{t_{k_{\ell}}+1} - \mathcal{F}(\omega^{t_{k_{\ell}}} \circ S_j z^{t_{k_{\ell}}}) \right) = -\Lambda_j^{\star}
\end{align*}
and
\begin{align*}
&\lim_{\ell \rightarrow \infty} \frac{1}{2}\left\| |u_j^{t_{k_{\ell}}+1}| - \sqrt{d_j} \right\|^2_2 = \frac{1}{2}\left\| |u_j^{\star}| - \sqrt{d_j} \right\|^2_2,\\
&\lim_{\ell \rightarrow \infty} 
      \displaystyle \frac{1}{2} \left \langle |u_j^{t_{k_{\ell}}+1}|^2 - d_j \log |u_j^{t_{k_{\ell}}+1}|^2, \mathbf{1}\right\rangle  =\displaystyle \frac{1}{2} \left \langle |u_j^{\star}|^2 - d_j \log |u_j^{\star}|^2, \mathbf{1}\right\rangle. 
\end{align*}
By closedness of limiting subdifferential, i.e., \eqref{eq:subdiff_closed}, we establish that 
\begin{align*}
0 \in \begin{cases}
    \partial|u_j^{\star}| \circ ( |u_j^{\star}| - \sqrt{d_j}) + \Lambda_j^{\star}, & \text{if AGM}, \\
   \partial|u_j^{\star}| \circ\left( |u_j^{\star}| - \displaystyle \frac{d_j}{|u_j^{\star}|} \right) + \Lambda_j^{\star}, &\text{if IPM.}
    \end{cases}
\end{align*}

At iteration $t_{k_{\ell}}$, the first-order optimality condition of \eqref{eq:v_sub} is
\begin{align}
    -\frac{y^{t_{k_{\ell}}}}{\lambda} - \frac{\beta_2}{\lambda} \left(v^{t_{k_{\ell}}+1} - \nabla z^{t_{k_{\ell}}} \right) \in \partial \left(\|v^{t_{k_{\ell}}+1}\|_1 - \alpha  \|v^{t_{k_{\ell}}+1}\|_{2,1} \right).
\end{align}
By \eqref{eq:v_star_converge} and \eqref{eq:v_converge}, we have
\begin{align*}
    &\lim_{\ell \rightarrow \infty} -\frac{y^{t_{k_{\ell}}}}{\lambda} - \frac{\beta_2}{\lambda} \left(v^{t_{k_{\ell}}+1} - \nabla z^{t_{k_{\ell}}} \right) = -\frac{y^{\star}}{\lambda}.
\end{align*}
By \eqref{eq:v_converge}, 
    $ \displaystyle \lim_{\ell \rightarrow \infty} \|v^{t_{k_{\ell}}+1}\|_1 - \alpha \|v^{t_{k_{\ell}}+1}\|_{2,1} =  \|v^{\star}\|_1 - \alpha \|v^{\star}\|_{2,1}.$
Altogether, by closedness of limiting subdifferential, we establish \eqref{eq:v_kkt}.

Lastly, we prove \eqref{eq:omega_kkt}-\eqref{eq:z_kkt}. At iteration $t_{k_{\ell}}$, we have
\begin{align*}
    \nabla_{\omega} \mathcal{L}(\omega^{t_{k_{\ell}}}) =  -\sum_{j=1}^N &\left\{\beta_1 (S_j z^{t_{k_{\ell}}})^* \circ \left[ \mathcal{F}^{-1}\left( u_{j}^{{t_{k_{\ell}}}+1}  + \frac{\Lambda_j^{t_{k_{\ell}}}}{\beta_1} \right) - \omega^{t_{k_{\ell}}} \circ S_j z^{t_{k_{\ell}}} \right] \right\},\\
    \nabla_{z} \mathcal{L}(z^{t_{k_{\ell}}}) =
    -\sum_{j=1}^N &\Bigg\{\beta_1 S_j^{\top}(D_{\omega^{t_{k_{\ell}}+1}})^*\left[ \mathcal{F}^{-1}\left(u_j^{t_{k_{\ell}}+1} + \frac{\Lambda_j^{t_{k_{\ell}}}}{\beta_1} \right) -\omega^{t_{k_{\ell}}+1} \circ S_j z^{t_{k_{\ell}}} \right]\\
    &+\beta_2 \left[ \nabla^{\top} \left(v^{t_{k_{\ell}}+1} + \frac{y^{t_{k_{\ell}}}}{\beta_2} \right)+ \Delta z^{t_{k_{\ell}}} \right] \Bigg\}.
\end{align*}
Taking the limit and using \eqref{eq:omega_diff_converge}, \eqref{eq:v_converge},\eqref{eq:u_t_k_l_equal_u_star}, we establish
\begin{align*}
    \nabla_{\omega} \mathcal{L}(Z^{\star}, \Omega^{\star}) &= \lim_{\ell \rightarrow \infty} \nabla_{\omega} \mathcal{L}(\omega^{t_{k_{\ell}}}),\\
    \nabla_{z} \mathcal{L}(Z^{\star}, \Omega^{\star}) &= \lim_{\ell \rightarrow \infty} \nabla_{z} \mathcal{L}(z^{t_{k_{\ell}}}).
\end{align*}
Since $\{\nabla_{\omega} \mathcal{L}(\omega^t), \nabla_z \mathcal{L}(z^t)\}_{t=1}^{\infty}$ is bounded, we have
\begin{align*}
\mathbb{E}\left[\|\nabla_{\omega} \mathcal{L}(Z^{\star}, \Omega^{\star}) \|_2^2\right] &= \lim_{\ell \rightarrow \infty}  \mathbb{E}[\|\nabla_{\omega} \mathcal{L}(\omega^{t_{k_{\ell}}})\|_2^2] = 0,\\
\mathbb{E}\left[\|\nabla_{z} \mathcal{L}(Z^{\star}, \Omega^{\star}) \|_2^2\right] &= \lim_{\ell \rightarrow \infty}  \mathbb{E}[\|\nabla_{z} \mathcal{L}(z^{t_{k_{\ell}}})\|_2^2] = 0
\end{align*}
by bounded convergence theorem and \eqref{eq:liminf_omega}-\eqref{eq:lim2}.
It follows that $\nabla_{\omega} \mathcal{L}(Z^{\star}, \Omega^{\star}),\nabla_{z} \mathcal{L}(Z^{\star}, \Omega^{\star})  =0.$ 
Altogether, $(Z^{\star}, \Omega^{\star})$ is a KKT point a.s.\
\end{proof}
\section*{References}
\bibliographystyle{siam}
\bibliography{references}

\end{document}